\documentclass[11pt,a4paper]{article}

\usepackage{amssymb,amsthm,amsmath,amsfonts}
\usepackage[pdftex]{graphicx}
\usepackage{graphicx}
\usepackage[margin=1in]{geometry}
\usepackage{color}
\usepackage[toc,page]{appendix}
\usepackage{hyperref}
\usepackage[export]{adjustbox}
\usepackage{caption}
\usepackage[skip=1ex, belowskip=2ex]{subcaption}
\usepackage{tikz,amsmath}
\usetikzlibrary{shapes.geometric, arrows,shadows}
\tikzstyle{process} = [rectangle,copy shadow={fill=black,shadow xshift=-0.5ex,shadow yshift=-0.5ex}, minimum width=3em, minimum height=2em, text centered, draw=black, fill=gray!10,]

\tikzstyle{arrow} = [thick,->,>=stealth]
\numberwithin{equation}{section}
\newtheorem{theorem}{Theorem}[section]
\newtheorem{lemma}[theorem]{Lemma}
\newtheorem{proposition}[theorem]{Proposition}

\newtheorem{corollary}[theorem]{Corollary}
\newtheorem{remark}[theorem]{Remark}
\theoremstyle{definition}

\theoremstyle{remark}

\newtheorem{assume}{Assumption}

\usepackage{hyperref}
\usepackage{comment}

\title{On the convergence rates of moment-SOS hierarchies approximation of truncated moment sequences}
\author{Hoang Anh Tran, Kim-Chuan Toh}

\newcommand{\bx}{\mathbf{x}}
\newcommand{\bu}{\mathbf{u}}
\newcommand{\by}{\mathbf{y}}
\newcommand{\iy}{\mathit{y}}
\newcommand{\bz}{\mathbf{z}}
\newcommand{\RR}{\mathbb{R}}
\newcommand{\NN}{\mathbb{N}}
\newcommand{\bv}{\mathbf{v}}
\newcommand{\CP}{\mathcal{P}}
\newcommand{\CQ}{\mathcal{Q}}
\newcommand{\CX}{\mathcal{X}}
\newcommand{\CT}{\mathcal{T}}
\newcommand{\bM}{\mathbf{M}}

\newcommand{\fmin}{f_{\min}}
\newcommand{\fmax}{f_{\max}}
\newcommand{\mlb}{\mathrm{mlb}}
\newcommand{\lb}{\mathrm{lb}}
\newcommand{\ub}{\mathrm{ub}}
\newcommand{\CM}{\mathcal{M}}
\newcommand{\cf}{\mathbf{f}}
\newcommand{\bd}{\mathbf{d}}
\newcommand{\BC}{\mathbf{C}}

\newcommand{\CY}{\mathcal{Y}}
\newcommand{\blambda}{\mathbf{\lambda}}

\newcommand{\BK}{\mathbf{K}}

\newcommand{\CR}{\mathcal{R}}
\newcommand{\Lo}{\textit{\L}}
\newcommand{\bw}{\mathbf{w}}
\newcommand{\conv}{\mathrm{conv}}

\def\degg#1{\lceil#1\rceil}
\def\bb{\mathbf{b}}

\newcommand{\interior}[1]{%
  {\kern0pt#1}^{\mathrm{o}}%
}

\begin{document}
\maketitle

\begin{abstract}
    The moment-SOS hierarchy is a widely applicable framework to address polynomial optimization problems over basic semi-algebraic sets based on positivity certificates of polynomial. Recent works show that the convergence rate of this hierarchy over certain simple sets, namely, the unit ball, hypercube, and standard simplex, is of the order $\mathrm{O}(1/r^2)$, where $r$ denotes the level of the moment-SOS hierarchy. This paper aims to provide a comprehensive understanding of the convergence rate of the moment-SOS hierarchy by estimating the Hausdorff distance between the set of truncated pseudo-moment sequences and the set of truncated moment sequences specified by Tchakaloff’s theorem. Our results provide a connection between the convergence rate of the moment-SOS hierarchy and the \L{}ojasiewicz exponent $\Lo$ of the domain under the compactness assumption, where we establish
    the convergence rate of $\mathrm{O}(1/r^\Lo)$.
    Consequently, we obtain the convergence rate of $\mathrm{O}(1/r)$ for polytopes and sets satisfying the constraint qualification condition, $\mathrm{O}(1/\sqrt{r})$ for domains that either satisfy the Polyak-Łojasiewicz condition or are defined by locally strongly convex polynomials. We also obtain the convergence rate of $\mathrm{O}(1/r^2)$ for general polynomials over a
    sphere.
\end{abstract}

\section{Introduction}
Consider the problem of minimizing a polynomial $f \in \RR[\bx]$ over a compact basic semi-algebraic set $\CX \subset \RR^n$:
\begin{equation}\tag{POP}\label{POP}
    \fmin = \underset{\bx \in \CX}{\min}\ f(\bx).
\end{equation}
The semi-algebraic set $\CX$ is defined by polynomial inequalities and equalities as follows:  
\begin{equation}\label{domain X}
\CX := \left\{ \bx \in \RR^n: \ g_j(\bx) \geq 0 \ \forall j \in [m],\;
h_i(x) = 0\ \forall\; i\in [p] \right\},
\end{equation}
where each $g_j$ and $h_i$ is a polynomial in $\RR[\bx]$. The class of polynomial optimization problems \eqref{POP} has wide applications in various fields, we refer to \cite{lasserre2009moments} for an overview on the existing techniques and applications. There are 2 types of moment-SOS hierarchies to address \eqref{POP}:  one approximates \eqref{POP} from below, and the other approximates from above, which we outline next.

\subsection{Hierarchies of lower bounds} 
The moment-SOS hierarchy of lower bounds as described in e.g. \cite{lasserre2001global, lasserre2009moments, lasserre2011new},  consists of a sequence of semidefinite programming (SDP) relaxations of \eqref{POP}. At the $r$-th level, the moment-SOS hierarchy approximates \eqref{POP} through an SDP, whose constraints are defined by $r$-truncated moment and localizing matrices or the sum-of-squares (SOS) representation of a positive polynomial. These relaxations form a sequence of lower bounds for $\fmin$, whose convergence is guaranteed by positivity certificates such as Putinar's Positivstellensatz and Schmüdgen's Positivstellensatz (see, e.g., \cite{anjos2011handbook}, \cite{lasserre2015introduction}).

The moment-SOS hierarchy can be categorized into two primary formulations based on the type of SDP relaxations: the primal formulation, known as the moment hierarchy, generates SDPs based on a generalized moment problem; and the dual formulation, known as the SOS hierarchy, generates SDPs based on the SOS representations of positive polynomials. Furthermore, the choice of positivity certificates influences the structure of these hierarchies. The most commonly employed Positivstellensatz are Schmüdgen’s and  Putinar’s 
Positivstellensatz. Consequently, four distinct types of hierarchies are derived: the Schmüdgen-type moment hierarchy \eqref{hierarchy: moment Schmudgen}, the Schmüdgen-type SOS hierarchy \eqref{hierarchy: SOS Schmudegen}, the Putinar-type moment hierarchy \eqref{hierarchy: moment Putinar}, and the Putinar-type SOS hierarchy \eqref{hierarchy: SOS Putinar}.
In terms of convergence, the Schmüdgen-type hierarchies have  faster convergence to the optimal value, but they are much more expensive in terms of computational complexity than the Putinar-type hierarchies. 

When the domain $\CX$ is a simple set -- specifically, unit ball, hypercube, and standard simplex, the existing works \cite{slot2111sum}, \cite{laurent2023effective} and \cite{S.o.S-on-simplex} have developed a method utilizing the Christoffel-Darboux (CD) kernel to approximate a positive polynomial by an SOS polynomial. This method leads to an explicit convergence rate of ${\rm O}({1/r^2})$ for Schmüdgen-type moment-SOS hierarchies of lower bounds. For the hypersphere $S^{n-1}$, the same convergence rate of $\mathrm{O}(1/r^2)$ is shown in \cite{fang2021sum} for homogeneous polynomial objective functions. When $\CX$ is the binary hypercube $\{0,1\}^n$, the convergence rate of $\mathrm{O}(1/r^2)$ is also available. Moreover, it it known from \cite{fawzi2016}, \cite{sakaue2017} that the corresponding Putinar-type moment-SOS hierarchy on $\{0,1\}^n$ is exact when $r \geq (n+d-1)/2$.

For a general compact semi-algebraic set $\CX$, general methods have been proposed to obtain the convergence rate of $\mathrm{O}(1/r^c)$ for the moment-SOS hierarchy in the work \cite{schmudgen-complexity}, where $c$ is a constant depending on $\CX$. Furthermore, improved versions of these convergence rates are shown in \cite{moment-approximation} and \cite{Baldi_2025} for the Putinar-type moment-SOS hierarchy. In particular, \cite{Baldi_2025} proved the convergence rate of $\mathrm{O}(1/r^{1/10})$ under the constraint qualification condition (CQC). Other works studying the convergence rates of the moment-SOS hierarchies of lower bounds with weaker results include \cite{de2010error}, \cite{S.o.S-on-simplex}, and \cite{Putinar-complexity}.

\subsection{Hierarchies of upper bounds}\label{sec: hierarchy of upper bounds} 
Lasserre’s approach begins by fixing a reference probability measure on the domain \(\CX\) and then relaxing \eqref{POP} into a convex optimization problem over the set of probability measures whose density functions are non-negative polynomials on \(\CX\) (see, e.g., \cite{lasserre2011new}). This formulation is further relaxed by replacing the set of non-negative polynomials on \(\CX\) with sums-of-squares (SOS) polynomials, the preordering, and the quadratic module, respectively. These relaxations lead to a semidefinite programming (SDP) formulation whose size depends polynomially on the number of variables $n$ and the degree bound $2r$ of the density polynomial. We note that this method requires the choice of a reference measure and its moment sequence on $\CX$.

 Using the CD kernel, it is known in \cite{slot2111sum}, \cite{S.o.S-on-simplex} that the convergence rates of the moment-SOS hierarchies of upper bounds on simple sets are $\mathrm{O}(1/r^2)$. The same convergence rate is also obtained for the minimization of a homogeneous polynomial over the 
 hypersphere $S^{n-1}$ in \cite{fang2021sum}. However, we should mention that this approach relies on the explicit formula of the CD kernel, which has been successfully calculated only for the above mentioned simple sets. For a more general domain: a compact full-dimensional
semi-algebraic set $\CX$ equipped with the Lebesgue measure, the convergence rate of $\mathrm{O}(\log^2 r/r^2)$ is proved for all types of 
hierarchies of upper bounds (see e.g., \cite{Slot_upper_bound}).

\subsection*{Contribution} 
In this paper, we propose an entirely different approach to analyze the convergence rate of the Schmüdgen-type moment-SOS hierarchy of lower bounds and the hierarchy of upper bounds as follows: rather than estimating the SOS representations of the objective function, we consider the error of truncated pseudo-moment sequences, which are the feasible solutions of either the SDP relaxation stated in \eqref{hierarchy: moment Schmudgen} or 
\eqref{hierarchy: reduced moment}. By  "error", we mean the minimum distance between the set of truncated pseudo-moment sequences and the set of truncated moment sequences supported on $\CX$. 
Since the problem \eqref{POP} is equivalent to the generalized moment problem \eqref{moment problem}, whose feasible solutions are truncated moment sequences, we can treat the feasible set of the SDP relaxation in each level of the moment hierarchy as an outer spectrahedral approximation of the set of truncated moment sequences, denoted by $\CM_k(\CX)$, where $k$ is the truncation order. Hence, to estimate the error of the moment hierarchy, we analyze the Hausdorff distance between these outer spectrahedral approximations and $\CM_k(\CX)$. We consider the upper bound on this distance as an "error" of a truncated pseudo-moment sequence in the sense of how far we can move a truncated pseudo-moment sequence to a truncated moment sequence, which then delivers the tightness of the SDP relaxations within the moment hierarchy by Lemma~\ref{lemma: distance to convergence rate}. 

Because the SDP relaxations in the SOS hierarchy are the dual of the SDP relaxations in the moment sequence, for which the strong duality holds under the Archimedean condition ( see e.g., \cite{josz2016strong}), this leads to an identical convergence rate between the moment and SOS hierarchies. In addition, we construct a new certificate denoted by $\CR(\CX)$ which is weaker than the Schmüdgen certificate, and potentially leads to a reduced version of the Schmüdgen-type hierarchy without changing the theoretical convergence rate. 
The reduction of the moment-SOS hierarchy for a real algebraic variety in \cite{S.o.S-on-real-variety} shares a similar construction, but our results provide precise convergence rates under the Archimedean condition. The connection of the error estimation of the truncated pseudo-moment sequences with the Łojasiewicz inequality directly implies the convergence rates in various special cases such as strongly convex sets, sets satisfying the Polyak-Łojasiewicz condition \eqref{def LP} or constraint qualification condition (CQC), polytopes, and spheres. In conclusion, the main results of this paper are summarized in Table~\ref{tab: results}.

\begin{table}[tbhp]
\footnotesize
  
\begin{center}
\renewcommand{\arraystretch}{2.0}
  \begin{tabular}{c c c c c} 
   Domain $\CX$ (Archimedean) & \bf Certificate & \bf Error & \bf Convergence rate &\bf Theorem/Corollary\\ \hline
   Unit ball & $\CR(\CX),\CQ(\CX),\CT(\CX)$ &  $\mathrm{O}(1/r^2)$ &  $\mathrm{O}(1/r^2)$ & \eqref{thm: error on product of simple sets}, \eqref{cor: rate on product}
   \\
   Standard simplex & $\CR(\CX),\CT(\CX)$ &  $\mathrm{O}(1/r^2)$ &  $\mathrm{O}(1/r^2)$ & \eqref{thm: error on product of simple sets}, \eqref{cor: rate on product}\\
   Product of simple sets & $\CR(\CX),\CT(\CX)$ &  $\mathrm{O}(1/r^2)$ &  $\mathrm{O}(1/r^2)$ & \eqref{thm: error on product of simple sets}, \eqref{cor: rate on product}\\
   Compact  & $\CR(\CX),\CT(\CX)$ &  $\mathrm{O}(1/r^{\Lo})$ &  $\mathrm{O}(1/r^{\Lo})$ & \eqref{thm: error on general}, \eqref{cor: convergence rate for compact semi-algebraic set}\\
   Polyak-Łojasiewicz condition & $\CR(\CX),\CT(\CX)$ &  $\mathrm{O}(1/\sqrt{r})$ &  $\mathrm{O}(1/\sqrt{r})$ & \eqref{thm: error and rate under the strong convexity}\\
   Strongly convex  & $\CR(\CX),\CT(\CX)$ &  $\mathrm{O}(1/\sqrt{r})$ &  $\mathrm{O}(1/\sqrt{r})$ &  \eqref{thm: error and rate under the strong convexity}\\
   Polytope & $\CR(\CX),\CT(\CX)$ &  $\mathrm{O}(1/r)$ &  $\mathrm{O}(1/r)$ & \eqref{cor: error and rate for polytope} 
   \\
   Sphere & $\CR(\CX), \CQ(\CX),\CT(\CX)$ &  $\mathrm{O}(1/r^2)$ &  $\mathrm{O}(1/r^2)$ & \eqref{thm: error on sphere} \eqref{cor: rate on sphere}\\
   Under CQC & $\CR(\CX),\CT(\CX)$ &  $\mathrm{O}(1/r)$ &  $\mathrm{O}(1/r)$ & \eqref{cor: error and rate under CQC} \eqref{cor: rate on sphere}
   \\ \hline
  \end{tabular}
\end{center}
\caption{Error on pseudo-moment sequence approximation and the convergence rate in terms of the Łojasiewicz exponent $\Lo$.}\label{tab: results}
\end{table}

For the convergence rate of the hierarchy of upper bounds, we use the same method as in 
Section 5 of \cite{slot2111sum}, where we use a generalization of the CD kernel to bound the error of the optimal value $\fmin$ and the optimal value of the SDP relaxation in the $r$-th level of the hierarchy of upper bounds when the domain is a product of  simple sets.  The result is stated in Theorem~\ref{thm: convergence rate for hierarchy of upper bounds}.

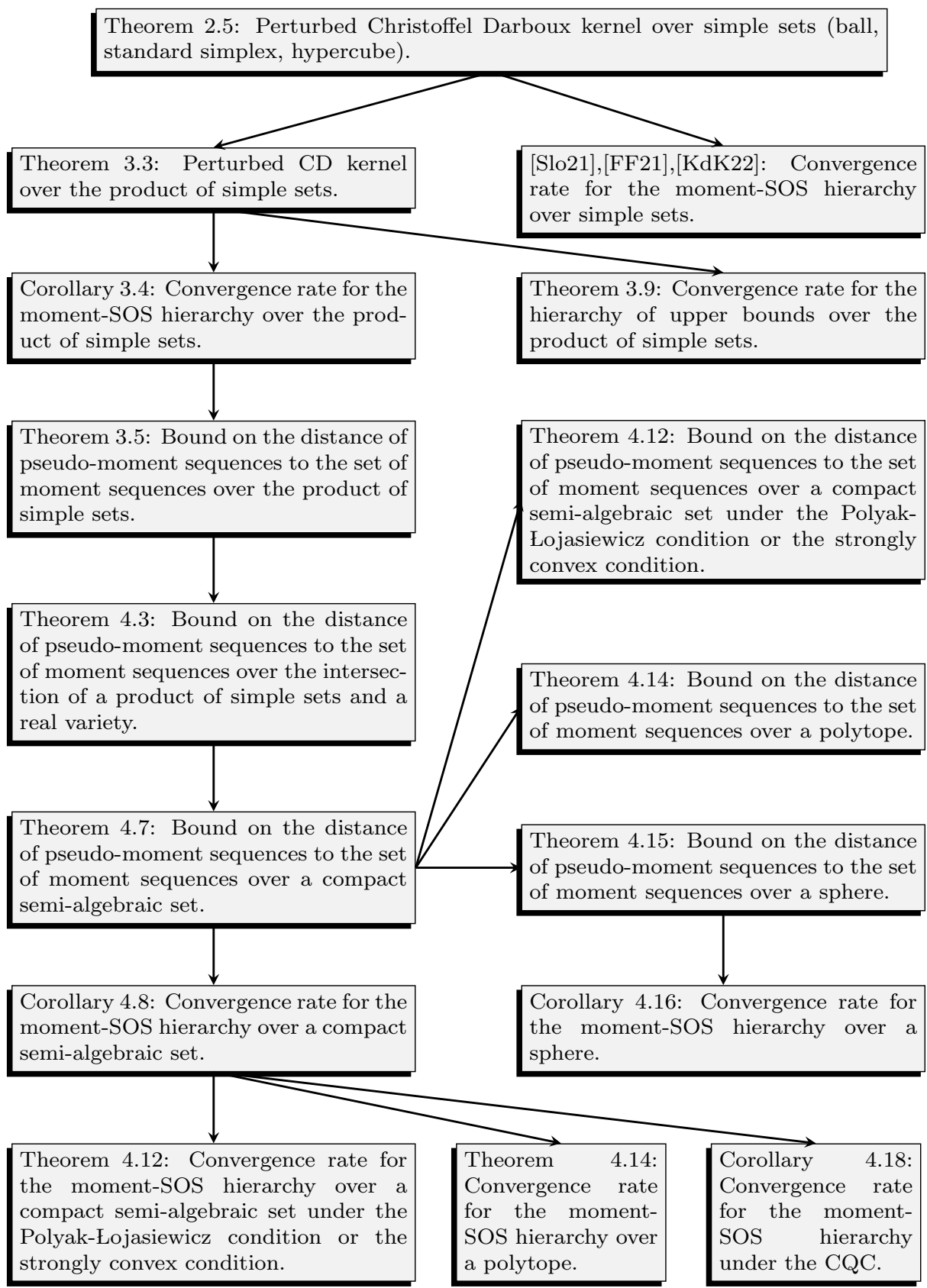
\begin{figure}[!htbp]
\centering
\hfill
\resizebox{\textwidth}{!}{
{\scriptsize \begin{tikzpicture}[node distance=2em]

\node (S1) [process]{\parbox[t][][t]{10cm}{ Theorem~\ref{thm: CD kernel on simple sets}: Perturbed Christoffel Darboux kernel over simple sets (ball, standard simplex, hypercube).} };

\node (S2) [process, below of=S1, yshift=-4em, xshift = -12em]{\parbox[t][][t]{5cm}{Theorem~\ref{thm: CD on product}: Perturbed  CD kernel over the product of simple sets.}};

\node (S3) [process, below of=S1, yshift=-4.5em, xshift = 10em]{\parbox[t][][t]{5cm}{\cite{slot2111sum},\cite{fang2021sum},\cite{S.o.S-on-simplex}: Convergence rate for the moment-SOS hierarchy over simple sets.}};

\node (S4) [process, below of=S2, yshift=-4em]{\parbox[t][][t]{5cm}{Corollary~\ref{cor: rate on product}: Convergence rate for the moment-SOS hierarchy over the product of simple sets.}};

\node (S5) [process, below of=S4, yshift=-5em]{\parbox[t][][t]{5cm}{Theorem~\ref{thm: error on product of simple sets}: Bound on the distance of pseudo-moment sequences to the set of moment sequences over the product of simple sets.}};

\node (S6) [process, below of=S3, yshift=-3.5em]{\parbox[t][][t]{5cm}{Theorem~\ref{thm: convergence rate for hierarchy of upper bounds}: Convergence rate for the hierarchy of upper bounds over the product of simple sets.}};

\node (S7) [process, below of=S5, yshift=-6.5em]{\parbox[t][][t]{5cm}{Theorem~\ref{thm: error on variety}: Bound on the distance of pseudo-moment sequences to the set of moment sequences over the intersection of a product of simple sets and a real variety.}};

\node (S8) [process, below of=S7, yshift=-6.5em]{\parbox[t][][t]{5cm}{Theorem~\ref{thm: error on general}: Bound on the distance of pseudo-moment sequences to the set of moment sequences over a compact semi-algebraic 
set.}};

\node (S9) [process, below of=S8, yshift=-5em]{\parbox[t][][t]{5cm}{Corollary~\ref{cor: convergence rate for compact semi-algebraic set}: Convergence rate for the moment-SOS hierarchy over a compact semi-algebraic set.}};

\node (S10) [process, below of=S6, yshift=-6em]{\parbox[t][][t]{5cm}{Theorem~\ref{thm: error and rate under the strong convexity}: Bound on the distance of pseudo-moment sequences to the set of moment sequences over a compact semi-algebraic set under the Polyak-{\L}ojasiewicz condition or the strongly convex condition.}};

\node (S11) [process, below of=S9, yshift=-6em]{\parbox[t][][t]{5cm}{Theorem~\ref{thm: error and rate under the strong convexity}: Convergence rate for the moment-SOS hierarchy over a compact semi-algebraic set under the Polyak-{\L}ojasiewicz condition or the strongly convex condition.}};

\node (S12) [process, below of=S10, yshift=-7em]{\parbox[t][][t]{5cm}{Theorem~\ref{cor: error and rate for polytope}: Bound on the distance of pseudo-moment sequences to the set of moment sequences over a polytope.}};

\node (S13) [process, below of=S12, yshift=-5em]{\parbox[t][][t]{5cm}{Theorem~\ref{thm: error on sphere}: Bound on the distance of pseudo-moment sequences to the set of moment sequences over a sphere.}};

\node (S14) [process, below of=S13, yshift=-5em]{\parbox[t][][t]{5cm}{Corollary~\ref{cor: rate on sphere}: Convergence rate for the moment-SOS hierarchy over a sphere.}};

\node (S15) [process, below of=S14,xshift = -7em, yshift=-6em]{\parbox[t][][t]{2.5cm}{Theorem~\ref{cor: error and rate for polytope}: Convergence rate for the moment-SOS hierarchy over a polytope.}};

\node (S16) [process, below of=S14,xshift = 4em, yshift=-6em]{\parbox[t][][t]{2.5cm}{Corollary~\ref{cor: error and rate under CQC}: Convergence rate for the moment-SOS hierarchy under the CQC.}};

\draw [arrow] (S1.south) -- node[anchor=east]{}(S2.north);
\draw [arrow] (S1.south) -- node[anchor=west]{}(S3.north);
\draw [arrow] (S2) -- (S4);
\draw [arrow] (S4) -- (S5);
\draw [arrow] (S2.south) -- (S6.north);
\draw [arrow] (S5) -- (S7);
\draw [arrow] (S7) -- (S8);
\draw [arrow] (S8) -- (S9);
\draw [arrow] (S8.east) -- (S10.west);
\draw [arrow] (S9.south) -- (S11.north);
\draw [arrow] (S8.east) -- (S12.west);
\draw [arrow] (S8.east) -- (S13.west);
\draw [arrow] (S9.south) -- (S15.north);
\draw [arrow] (S9.south) -- (S16.north);
\draw [arrow] (S13.south) -- (S14.north);

\end{tikzpicture}}}
\caption{Flow chart for the results established in this paper.}
\label{Fig: order of proofs}
\end{figure}

The results in this paper are systematically presented and closely interconnected. Their relationships are outlined in the flowchart in Figure~\ref{Fig: order of proofs} for clarity.

\section{Preliminaries}
\subsection{Notation and SOS polynomial} 
\label{subsec-notation}

We denote a closed ball in an Euclidean space with center at the origin
and radius $R$ by $\mathbb{B}_R$. We use $\|\bx\|$ to denote the Euclidean norm of a vector $\bx\in\RR^n$. The distance between a point $\bx$ and a set $\mathcal{A}$ in an Euclidean space is defined as $\bd(\bx,\mathcal{A}) = \inf \{\|\by-\bx\| : \by\in \mathcal{A}\}$, and 
the Hausdorff distance between two sets $\mathcal{A}, \mathcal{B}$ 
is defined by $\bd(\mathcal{A},\mathcal{B}) = \sup\{ \bd(\bx,\mathcal{B}) : \bx \in \mathcal{A}\}$. For any integer $m\in \NN$, $[m]:=\{1,\ldots,m\}.$

We use $\bx=(x_1,\dots,x_n)$ to denote a vector of variables and $\RR[\bx]$ as the ring of polynomials in $\bx$. Let $\alpha =(\alpha_1,\dots,\alpha_n)$ be a multi-index with length $|\alpha| = \sum_{i=1}^n \alpha_i$. The set of multi-index of length at most $r$ is denoted by $\NN^n_r = \{\alpha \in \NN^n: \ |\alpha| \leq r\}$. We let $\overline{\NN}^n_{r}$ to be the subset of $\NN^n_r$ whose elements have length exactly $r$. The monomials in $\bx$ are written in the form $\bx^{\alpha}= x_1^{\alpha_1}\cdots x_n^{\alpha_n}$. For any polynomial $f(\bx) = \sum_{\alpha}f_{\alpha}\bx^{\alpha} \in \RR[\bx]$, we define the norm $\|f\|_{1} = \sum_{\alpha}|f_{\alpha}|$, and  $\lceil f \rceil := \lceil (\mbox{deg}\ f) /2 \rceil$.

We denote the basis vector containing all standard monomials in $\bx$ and the 
$r$-truncated basis vector of all monomials of degree up to $r$, respectively, by 
\begin{displaymath}
\bv(\bx) = (\bx^{\alpha})_{\alpha \in \NN^n}, \quad \text{and} \quad \bv_r(\bx) = (\bx^{\alpha})_{\alpha \in \NN_r^n}.
\end{displaymath}
The dimension of the basis vector $\bv_r$ is $s(n,r)=\binom{n+r}{n}$. 
Using the monomials basis, any polynomial $f$ of degree $d$ can be expressed in the form:
\begin{displaymath}
    f(\bx) = \sum_{\alpha \in \NN^n_d}f_{\alpha} \bx^{\alpha}= \langle \mathbf{f},\bv_d(\bx) \rangle,\ \mathbf{f} \in \RR^{s(n,d)}.
\end{displaymath}

A polynomial $f$ is called a sum-of-squares (SOS) if there exists a finite number of polynomials $f_1,\dots,f_N$ such that $f(\bx)= \sum_{i=1}^N f_i(\bx)^2$. We denote the set of SOS polynomials and its subset of SOS polynomials of degree at most $2r$ by $\Sigma[\bx]$ and $\Sigma[\bx]_{2r}$, respectively.

Recall the basic semi-algebraic set $\CX$ in \eqref{domain X}.
For an index set $ J \subset [m]$, we define $g_J = \prod_{j \in J}g_j$, and
 $g_{\emptyset}=1.$
The truncated preordering and quadratic module of $\CX$ are defined, respectively, by 
\begin{eqnarray} 
    \CT(\CX)_{2r} &=& \Biggl\{q(\bx) = \sum_{j=1}^N\sigma_{J_j}(\bx)g_{J_j}(\bx) + \sum_{i=1}^p\tau_i(\bx)h_i(\bx) :\ 
     \label{eq-TX} \\ 
    &&\quad \exists N \in \NN,\; J_j\subset [m], \;\sigma_{J_j} \in \Sigma[\bx]_{2(r-\lceil g_{J_j} \rceil) } \ \forall j \in [N],\ \tau_i \in \RR[\bx]_{2(r-\lceil h_i \rceil)}\ \forall i \in [p] \ \Biggr\},
\nonumber\\
        \CQ(\CX)_{2r}&=&\Biggl\{q(\bx) = \sigma_0(\bx) +\sum_{j=1}^m\sigma_j(\bx)g_j(\bx) 
        +\sum_{i=1}^p\tau_i(\bx)h_i(\bx):
        \label{eq-QX} \\ 
    &&\quad \sigma_0 \in \Sigma[\bx]_{2r},\;
     \; \sigma_j \in \Sigma[\bx]_{2(r-\lceil g_j \rceil)}\ \forall\; j\in[m],\;
     \tau_i \in \RR[\bx]_{2(r-\lceil h_i \rceil)}\ \forall i \in [p]\Biggr\}.
    \nonumber
\end{eqnarray}
In the above definitions, the conditions $\degg{g_{J_j}} \leq r$ for all $j\in[N]$, $\degg{g_j}\leq r$ for all $j\in [m]$, 
and $\degg{h_i}\leq r$ for all $i\in[p]$ are assumed. 
The preordering $\CT(\CX)$  of $\CX$ is defined by removing the degree constraints on $\sigma_{J_j}\,\forall\, j\in[N]$ and $\tau_i\,\forall\, i\in[p]$ in $\CT(\CX)_{2r}$.
Similarly, the quadratic module $\CQ(\CX)$ of $\CX$ is defined by removing 
the degree constraints on 
$\sigma_0$, $\sigma_j\, \forall\, j\in [m]$ and 
$\tau_i \,\forall\, i\in[p]$ in $\CQ(\CX)_{2r}$.  

It is clear that $\CT(\CX)_{2r}$ and $\CQ(\CX)_{2r}$ are subsets of the set of non-negative polynomials over $\CX$. Furthermore, checking the membership of a polynomial in $\CT(\CX)_{2r}$ and $\CQ(\CX)_{2r}$ can be verified by an SDP. Hence, these sets are the relaxation of the set of non-negative polynomials corresponding to the Schmüdgen Positivstellensatz for a compact semi-algebraic set (see e.g., \cite[pp. 283--313]{schmudgen2017moment})
and the Putinar Positivstellensatz for an Archimedean semi-algebraic set (see e.g., \cite{putinar1993positive}), respectively. 

\begin{theorem}[Schmüdgen Positivstellensatz]\label{thm: Schmüdgen}
    Let $\CX$ be the semi-algebraic set in \eqref{domain X}.
    We assume that $\CX$ is compact. If $f$ is a positive polynomial on $\CX$, then $f \in \CT(\CX)$.
\end{theorem}
\begin{theorem}[Putinar Positivstellensatz]\label{thm: Putinar}
    Let $\CX$ be the semi-algebraic set in \eqref{domain X}. We assume that the Archimedean condition holds, i.e, there exists a positive number $R$ such that $R -\|\bx\|^2 \in \CQ(\CX)$. If $f$ is a positive polynomial on $\CX$, then $f \in \CQ(\CX).$
\end{theorem}

We now define a novel reduced version of 
$\CT(\CX)_{2r}$ as follows: For any $r\in\NN$, the reduced version of $\CT(\CX)_{2r}$  is defined as  
\begin{eqnarray} 
    \CR(\CX)_{2r} &:=& \Biggl\{q(\bx) = \sum_{j=1}^N\sigma_{J_j}(\bx)g_{J_j}(\bx) + \sum_{i=1}^p\tau_ih_i^2(\bx) \in 
    \CT(\CX): \; \tau_i \in \RR_{\geq 0} \ \forall i \in [p],\;
    \nonumber \\    
   &&\quad \exists N \in \NN,\; J_j\subset [m],\; \sigma_{J_j} \in \Sigma[\bx]_{2(r-\lceil g_{J_j} \rceil) } \ \forall j \in [N]
     \Biggr\}.
     \label{eq-RX}
\end{eqnarray}
It is clear that for any positive integer $r$, $\CR(\CX)_{2r} \subset \CT(\CX)_{2r}$.
In this paper, we analyze the error of the truncated pseudo-moment sequences associated 
with $\CR(\CX)_{2r}$ instead of $\CT(\CX)_{2r}$.
\subsection{The moment-SOS hierarchy}
\label{sec-2.2}

Let $\iy = (y_{\alpha})_{\alpha \in \NN^n}$ be a real sequence indexed by the
vector of monomials $\bv(\bx)$. We define the Riesz linear functional $\ell_{\iy}: \RR[\bx] \to \RR$ as follows: 
\begin{displaymath}
    f(\bx)= \sum_{\alpha \in \NN^n}f_{\alpha}x^{\alpha} \quad \mapsto \quad \ell_{\iy}(f) = \sum_{\alpha \in \NN^n}f_{\alpha}y_{\alpha}.
\end{displaymath}
The Riesz linear functional plays a central role in determining whether a sequence $\iy$ is a moment sequence for a Borel measure (see the Riesz-Haviland Theorem, e.g., \cite[Theorem 3.1]{lasserre2009moments}). We utilize $\ell_{\iy}$ to set up the moment matrix and localizing matrix as follows: given an infinite sequence $\iy$ as above, the moment matrix $\bM(\iy)$ with rows and columns indexed by $\bv(\bx)$ is defined by
\begin{displaymath}
    \bM(\iy)(\alpha, \beta) = \ell_{\iy}(\bx^{\alpha+\beta})=y_{\alpha +\beta}, \quad \forall \alpha, \beta \in \NN^n.
\end{displaymath}
For a given $r \in \NN$, the $r$-truncated moment matrix, denoted by $\bM_r(\iy)$, is the submatrix of $\bM(\iy)$ obtained by extracting the rows and columns of $\bM(\iy)$
indexed by $\bv_r(\bx)$.
Similarly, for a given polynomial $g \in \RR[\bx]$, the localizing matrix $\bM(g\iy)$ associated with $\iy$ and $g$ is defined by 
\begin{displaymath}
    \bM(g\iy)(\alpha,\beta)=\ell_{\iy}(g(\bx)\bx^{\alpha+\beta})
    = \sum_{\gamma} g_\gamma y_{\gamma+\alpha+\beta}, \quad \forall \alpha, \beta \in \NN^n.
\end{displaymath}
The $r$-truncated localizing matrix is similarly constructed by extracting all the rows and columns indexed by $\bv_r(\bx)$ from the the localizing matrix $\bM(g\iy)$. 

We now revisit the moment and SOS hierarchies used to solve the problem \eqref{POP} with $\CX$ defined as in \eqref{domain X}. These hierarchies come in two forms: one based 
on Schmüdgen’s Positivstellensatz in Theorem~\ref{thm: Schmüdgen}, and the other on
Putinar’s Positivstellensatz 
in Theorem~\ref{thm: Putinar}. For any $r \in \NN$ such that $r \geq 
\max\{\lceil f \rceil,\lceil g_1 \rceil,\ldots,\lceil g_m \rceil\} $, we define the hierarchies as follows:
\begin{align}\label{hierarchy: moment Schmudgen}
   & \mlb(f,\CT(\CX))_r = \inf \Big\{ \ell_{\iy}(f) = \sum_{\alpha \in \NN^n_{2r}}f_{\alpha}y_{\alpha} \,:\, \iy \in \CM(\CT(\CX)_{2r})\Big\}
    \\
    \text{where \quad} & \CM(\CT(\CX)_{2r}) := \Big\{ 
    \iy\in 
    \RR^{s(n,2r)}\,:\,\ y_0 =1,\ \bM_r(\iy) \succeq 0,\ \bM_{r - \lceil h_i \rceil}(h_i \iy) = 0 \; \forall i \in [p] \nonumber\\ 
    &\hspace{30mm}  \bM_{r-\lceil g_J \rceil}(g_J\iy) \succeq 0\;\;
     \forall\; J \subset [m] \;\mbox{such that} \; \lceil g_J \rceil\leq r \Big\}
    \nonumber.
\end{align}
The elements of $\CM(\CT(\CX)_{2r})$ are called pseudo-moment sequences.
This forms the Schmüdgen-type moment hierarchy, whose optimal values generate a sequence of lower bounds for $\fmin$. One can see that $\mlb(f,\CT(\CX))_r$ is an SDP, whose dual problem is given by 
\begin{equation}\label{hierarchy: SOS Schmudegen}
    \lb(f,\CT(\CX))_r = \sup\{c \in \RR:\ f(\bx) -c \in \CT(\CX)_{2r}\}.
\end{equation}
The hierarchy \eqref{hierarchy: SOS Schmudegen} is called the Schmüdgen-type SOS hierarchy, whose convergence to the optimal value $\fmin$ is guaranteed by the Schmüdgen’s Positivstellensatz
in Theorem ~\ref{thm: Schmüdgen}. Whence, we obtain that 
\begin{displaymath}
 \lb(f,\CT(\CX))_r \leq \mlb(f,\CT(\CX))_r \quad \forall r \in \NN, \quad \lim_{r \to \infty}\lb(f,\CT(\CX))_r=\lim_{r \to \infty}\mlb(f,\CT(\CX))_r=\fmin.
\end{displaymath}

In the same manner, when the Archimedean condition is met, we have the 
Putinar-type version of the moment-SOS hierarchy as follows: 
\begin{align}\label{hierarchy: moment Putinar}
    & \mlb(f,\CQ(\CX))_r = \inf\Big\{ \ell_{\iy}(f) = \sum_{\alpha \in \NN^n_{2r}}f_{\alpha}y_{\alpha}\,:\, \iy \in {\cal M}(\CQ(\CX)_{2r}) \Big\}
    \\
    \text{where \quad}& {\cal M}(\CQ(\CX)_{2r}) :=\Big\{ \iy \in \RR^{s(n,2r)}\,:\, y_0 =1,\ \bM_r(\iy) \succeq 0, \;
      \bM_{r-\lceil g_i \rceil}(g_i\iy) \succeq 0 \ \forall\; i \in [m]\Big\}.
    \nonumber
\end{align}
These SDP relaxations form the Putinar-type moment hierarchy, whose 
dual problems form the Putinar-type SOS hierarchy defined by
\begin{equation}\label{hierarchy: SOS Putinar}
    \lb(f,\CQ(\CX))_r = \sup\{c \in \RR:\ f(\bx) -c \in \CQ(\CX)_{2r}\}.
\end{equation}
Under the Archimedean condition, the strong duality between the primal SDP 
\eqref{hierarchy: moment Putinar} and dual SDP \eqref{hierarchy: SOS Putinar} in the same level of the Putinar-type hierarchy holds, i.e, 
$\lb(f,\CQ(\CX))_r=\mlb(f,\CQ(\CX))_r$ (see e.g., \cite{josz2016strong}). The convergence of the Putinar-type hierarchy is based on the Putinar’s Positivstellensatz in Theorem~\ref{thm: Putinar}.

The moment hierarchy associated with $\CR(\CX)_{2r}$ in \eqref{eq-RX} is defined by
\begin{align}\label{hierarchy: reduced moment}
   & \mlb(f,\CR(\CX))_r = \inf \Big\{ \ell_{\iy}(f) = \sum_{\alpha \in \NN^n_{2r}}f_{\alpha}y_{\alpha} \,:\, \iy \in \CM(\CR(\CX)_{2r})\Big\}
    \\
    \text{where \quad} & \CM(\CR(\CX)_{2r}) := \Big\{ 
    \iy\in 
    \RR^{s(n,2r)}\,:\,\ y_0 =1,\ \bM_r(\iy) \succeq 0, \nonumber\\ 
    &\hspace{30mm}  \ell_{\iy}(h_i^2(\bx))=0\ \forall i \in [p],\; \bM_{r-\lceil g_J \rceil}(g_J\iy) \succeq 0\;\;
     \forall\; J \subset [m] \;\mbox{s.t} \; \lceil g_J \rceil\leq r \Big\}
    \nonumber,
\end{align}
whose  dual problem is 
\begin{equation}\label{hierarchy: reduced SOS}
    \lb(f,\CR(\CX))_r = \sup\{c \in \RR:\ f(\bx) -c \in \CR(\CX)_{2r}\}.
\end{equation}
The elements of $\CM(\CR(\CX)_{2r})$ are also called pseudo-moment sequences.
It is straightforward from the definitions that the following inequalities hold:
\begin{displaymath}
    \lb(f,\CR(\CX))_r \leq  \mlb(f,\CR(\CX))_r \leq  \mlb(f,\CT(\CX))_r, \quad \lb(f,\CR(\CX))_r \leq \lb(f,\CT(\CX))_r.
\end{displaymath}
\begin{remark}\label{remark: Archimedian condition}
    In this paper, we aim to determine the asymptotic convergence rate of the moment-SOS hierarchy for a compact semi-algebraic set $\CX$. Since $\CX$ is bounded, there exists a constant $R$ such that the ball $\mathbb{B}_R$ (with radius $R$ and center at the origin) contains $\CX$. Without loss of generality, we can assume that $R- \|\bx\|^2$ is positive over $\CX$. The Schm\"udgen Positivstellensatz theorem implies that there exists a positive integer $t$ such that $R- \|\bx\|^2 \in \CT(\CX)_{2t}$. Therefore, if we set $\CT(\CX')_{2r}$ to be the preordering of order $2r$ of $\CX$ 
    with the additional constraint $R- \|\bx\|^2 \geq 0$, then we have 
\begin{displaymath}
     \CT(\CX')_{2r}\subset \CT(\CX)_{2r+2t} \quad \forall\; r \in \NN.
\end{displaymath}
This means that the asymptotic convergence rates of $\mlb(f,\CT(\CX))_r$ and $\lb(f,\CT(\CX))_r$ are equal to that of $\mlb(f,\CT(\CX'))_r$ and $\lb(f,\CT(\CX'))_r$, respectively. Therefore, without loss of generality,
throughout this paper, most results will be stated under the assumption that the ball constraint $R^2 - \|\bx\|^2 \geq 0$ for some suitable $R$ is added to the description~\eqref{domain X} of $\CX$. 
\end{remark}

\subsection{Hausdorff distances}

Let $\CM(\CX)$ denote the set of all moment sequences associated with a probability measure on $\CX$, and denote the set of all probability measures on $\CX$ by $\CP(\CX)$.
Next, we change the point of view for the moment hierarchy as follows: the problem \eqref{POP} admits an equivalent formulation for any integer $k \geq \degg{f}$ as follows:
\begin{equation}\label{moment problem}
     \fmin = \inf\left\{\int_{\CX}fd\mu: \mu \in \CP(\CX) \right\}=\inf\left\{\sum_{\alpha \in \NN^n_k} f_{\alpha}y_{\alpha}=\langle \cf,\iy \rangle:\ \iy \in \CM_k(\CX) \right\}.
\end{equation}
Here, $\cf = (f_{\alpha})_{\alpha \in \NN^n_k} \in \RR^{s(n,k)}$. Then, \eqref{POP} is a linear optimization problem on a convex set $\CM_k(\CX) \subset \RR^{s(n,k)}$, where $\CM_k(\CX)$ denotes the set of $k$-truncated moment sequences of $\CM(\CX)$.

For $k\leq 2r$, let  $\pi_k :\RR^{s(n,2r)} \to \RR^{s(n,k)}$ denote the  projection onto the first $s(n,k)$ coordinates. Then the primal SDP problems in the moment hierarchies \eqref{hierarchy: moment Schmudgen},  \eqref{hierarchy: reduced moment}, and \eqref{hierarchy: moment Putinar}  can be written as the following alternatives:
\begin{align*}
    &\mlb(f,\CT(\CX))_r = \inf\left\{\sum_{\alpha \in \NN^n_k} f_{\alpha}y_{\alpha}=\langle \cf_k,\iy \rangle:\ \iy \in \CM_k(\CT(\CX)_{2r}) \right\},\\
    &\mlb(f,\CR(\CX))_r = \inf\left\{\sum_{\alpha \in \NN^n_k} f_{\alpha}y_{\alpha}=\langle \cf_k,\iy \rangle:\ \iy \in \CM_k(\CR(\CX)_{2r}) \right\},\\
    &\mlb(f,\CQ(\CX))_r = \inf\left\{\sum_{\alpha \in \NN^n_k} f_{\alpha}y_{\alpha}=\langle \cf_k,\iy \rangle:\ \iy \in \CM_k(\CQ(\CX)_{2r}) \right\},
\end{align*}
where 
\begin{eqnarray*}
\CM_k(\CT(\CX)_{2r}) &=& \{ \pi_k(y)\in \RR^{s(n,k)} \,:\, y\in \CM(\CT(\CX)_{2r}) \},
\\
\CM_k(\CR(\CX)_{2r}) &=& \{ \pi_k(y)\in \RR^{s(n,k)} \,:\, y\in \CM(\CR(\CX)_{2r}) \},
\\
\CM_k(\CQ(\CX)_{2r}) &=& \{ \pi_k(y)\in \RR^{s(n,k)} \,:\, y\in \CM(\CQ(\CX)_{2r}) \}.
\end{eqnarray*}
We have that $\CM_k(\CT(\CX)_{2r})$ and $\CM_k(\CQ(\CX)_{2r})$ are outer convex approximations of $\CM_k(\CX)$, i.e.,
\begin{displaymath}
    \CM_k(\CX) \subset \CM_k(\CT(\CX)_{2r}) \subset \CM_k(\CQ(\CX)_{2r}).
\end{displaymath}
To study the error of truncated pseudo-moment sequences, we analyze the bound on the following Hausdorff distances:
\begin{eqnarray*}\label{hausdorff distance}
   \qquad \bd_k(\CT(\CX)_{2r}) &:=&\bd(\CM_k(\CT(\CX)_{2r}),\CM_k(\CX))
   =\max\{\bd(\iy,\CM_k(\CX)):\ \iy \in \CM_k(\CT(\CX)_{2r})\},\\
\qquad \bd_k(\CR(\CX)_{2r}) &:=& \bd(\CM_k(\CR(\CX)_{2r}),\CM_k(\CX))
=\max\{\bd(\iy,\CM_k(\CX)):\ \iy \in \CM_k(\CR(\CX)_{2r})\},\\
  \qquad  \bd_k(\CQ(\CX)_{2r}) &:=& \bd(\CM_k(\CQ(\CX)_{2r}),\CM_k(\CX))
  =\max\{\bd(\iy,\CM_k(\CX)):\ \iy \in \CM_k(\CQ(\CX)_{2r})\},
\end{eqnarray*}
which can be used to establish the convergence rates of the associated moment hierarchies as stated in the following lemma.

\begin{lemma}\label{lemma: distance to convergence rate}
    Let $\CX$ be a compact semi-algebraic set and $k \geq {\rm deg}(f)$. Then the errors of the $r$-level of the moment hierarchies are bounded proportionally to the Hausdorff distances as follows:
    \begin{eqnarray*}
        \qquad \fmin - \mlb(f,\CT(\CX))_r &\leq& \|f\|_1\,\bd_k(\CT(\CX)_{2r}),\\
        \qquad \fmin - \mlb(f,\CR(\CX))_r &\leq&\|f\|_1\,\bd_k(\CR(\CX)_{2r}),\\
        \qquad \fmin - \mlb(f,\CQ(\CX))_r &\leq&\|f\|_1\,\bd_k(\CQ(\CX)_{2r}).
    \end{eqnarray*}
    Combining with Remark~\ref{remark: Archimedian condition}, the convergence rates of the pairs of moment-SOS hierarchies \eqref{hierarchy: moment Schmudgen}--\eqref{hierarchy: SOS Schmudegen}, 
    \eqref{hierarchy: reduced moment}--\eqref{hierarchy: reduced SOS}, and 
    \eqref{hierarchy: moment Putinar}--\eqref{hierarchy: SOS Putinar} are the same as the rates of the Hausdorff distances $\bd_k(\CT(\CX)_{2r})$,
    $\bd_k(\CR(\CX)_{2r})$, and $\bd_k(\CQ(\CX)_{2r})$, respectively. 
\end{lemma}
\begin{proof}
We only prove the result for the preordering $\CT(\CX)_{2r}$ since it is similar for the quadratic module $\CQ(\CX)_{2r}$ and $\CR(\CX)_{2r}$. Notice that the problem \eqref{POP} admits an equivalent formulation defined by $k$-truncated moment sequences as follows:
    \begin{displaymath}
        \fmin = \inf \left\{\sum_{\alpha \in \NN^n_k}f_{\alpha}\iy_{\alpha}=\langle \mathbf{f},\iy \rangle:\ \iy \in \CM_k(\CX) \right\}.
    \end{displaymath}
    Since $\CX$ is compact, so is $\CM_k(\CX)$. For any $\iy \in \CM_k(\CT(\CX)_{2r})$, there exists its projection $\overline{\iy} \in \CM_k(\CX)$ such that $\|\iy - \overline{\iy}\| = \bd_k(\CT(\CX)_{2r})$. Hence the Cauchy–Schwarz inequality implies that 
    \begin{displaymath}
       \left| \langle \mathbf{f},\iy \rangle- \langle \mathbf{f},\overline{\iy} \rangle\right| \leq \|f\|_1\bd_k(\CT(\CX)_{2r})\ \Rightarrow \  \fmin - \mlb(f,\CT(\CX))_r \leq \|f\|_1\,\bd_k(\CT(\CX)_{2r}).
    \end{displaymath}
    This completes the proof.
\end{proof}

\subsection{The hierarchies of upper bounds} Consider the alternative form of \eqref{POP} defined as follows:
\begin{equation}\label{form of upper bounds}
    \fmin = \inf_{\nu \in \CM_+(\CX)} \left\{\int_{\CX}f(\bx)d\nu(\bx)\; : \; \int_{\CX}d\nu(\bx) =1 \right\},
\end{equation}
where $\CM_{+}(\CX)$ denotes the set of positive measures supported on $\CX$. The idea of Lasserre is to relax $\CM_+(\CX)$ into the set of measures that are absolutely continuous with respect to a fixed reference measure $\mu$ supported on $\CX$. It then continues to 
inner approximate the set of density functions by the quadratic module and the preordering of $\CX$ for different level $r \in \NN$:
\begin{align*}
    \ub(f, \CQ(\CX))_r &:= \inf_{q \in \CQ(\CX)_{2r}}\left\{\int_{\CX} 
    f(\bx)q(\bx)d\mu(\bx) \;:\; \int_{\CX}q(\bx)d\mu(\bx)=1 \right\},\\ 
    \ub(f, \CT(\CX))_r &:= \inf_{q \in \CT(\CX)_{2r}}\left\{\int_{\CX}
    f(\bx)q(\bx)d\mu(\bx) \;:\; \int_{\CX}q(\bx)d\mu(\bx)=1 \right\}.
\end{align*}
These hierarchies of upper bounds (for $f_{\min}$) 
are called the Putinar-type hierarchy of upper bounds, and the Schm\"udgen-type hierarchy of upper bounds, respectively. We have summarized the works studying the convergence rates of these hierarchies of upper bounds in Section~\ref{sec: hierarchy of upper bounds}

\subsection{Christoffel-Darboux kernel}

Despite the convergence of the moment-SOS hierarchies based on Putinar and Schmüdgen Positivstellensatzs, the convergence rate of these hierarchies are challenging to analyze in general. However, when the domain is simple, recent investigation on the convergence rate has achieved substantial success by using the corresponding Christoffel-Darboux kernel (CD kernel) to show the convergence rate of $\mathrm{O}(1/r^2)$. In this paper, our methodology uses the CD kernel for the simple sets, unit ball $B_n$ and standard simplex $\Delta_n$, combining with the Łojasiewicz inequality to study the convergence rate of the moment-SOS hierarchy in general. This section summarizes the technique of the CD kernel that we use throughout this paper.

Let $\CX$ be a compact subset of $\RR^n$, and $\mu$ be a probability measure whose support is exactly $\CX$. The measure $\mu$ defines an inner product in the ring of polynomials $\RR[\bx]$ as follows:
\begin{displaymath}
    \langle p, q \rangle_{\mu} = \int_{\CX}p(\bx)q(\bx)d\mu(\bx), \quad \forall\; p,\; q \in \RR[\bx].
\end{displaymath}

Let $\{P_{\alpha}:\ \alpha \in \NN^n\}$ be an orthonormal basis with respect to the inner product $\langle \cdot, \cdot \rangle_{\mu}$, where ${\rm deg} (P_{\alpha}) = |\alpha|$. 
 The $\lambda$-perturbed CD kernel of degree $2r$ associated with the measure $\mu$ and weight vector $\lambda=(\lambda_i)_{i=0,\dots,2r}$ is defined as follows:
 \begin{displaymath}
     C_{2r}[\CX,\mu,\lambda](\bx,\by):= \sum_{i=0}^{2r}\lambda_iC^{(i)}[\CX,\mu](\bx,\by),
     \quad\mbox{with} \;\; C^{(i)}[\CX,\mu](\bx,\by) := \sum_{|\alpha|=i}P_{\alpha}(\bx)P_{\alpha}(\by). 
 \end{displaymath}  
For any $ 2r\geq k$, we define the linear operator $\mathbf{C}_{2r}$ associated with the polynomial kernel $C_{2r}(\bx,\by,\lambda)$ by 
 \begin{displaymath}
     \mathbf{C}_{2r}[\CX,\mu,\lambda]:\; \RR[\bx]_k \to \RR[\bx]_k,\quad \mathbf{C}_{2r}[\CX,\mu,\lambda]p(\bx) = \int_{\CX}
     C_{2r}[\CX,\mu,\lambda](\bx,\by)p(\by)d\mu(\by).
 \end{displaymath}
 The measures used to construct the CD kernels on simple sets are given in 
 Table~\ref{tab:measure}. The following theorem summarizes the properties
 of the operator $\mathbf{C}_{2r}[\CX,\mu,\lambda]$ for the unit ball and the standard simplex established in \cite{slot2111sum}. We refer to \cite{fang2021sum,laurent2023effective} for the analogue of this result for the unit sphere and the hypercube, respectively.

\begin{table}[tbhp]

\footnotesize
\begin{center}
\renewcommand{\arraystretch}{2.0}
\begin{tabular}{c c c}
\textbf{domain $\CX$} & \textbf{measure $\mu$} & \textbf{reference} \\
\hline
Unit ball $B_n$ & $(1-\|\bx\|^2)^{-1/2}d\bx$ & \cite{slot2111sum} \\
\par
Standard simplex $\Delta_n$ &$(1-|\bx|)^{-1/2}\prod_i(1-x_i)^{-1/2}d\bx$& \cite{slot2111sum}\\
\par
Unit sphere $S^n$ & Haar measure on $SO(n)$  & \cite{fang2021sum} \\
\par
Hypercube $\mathbb{B}^n$ & $\prod_i(1-x_i)^{-1/2}d\bx$ & \cite{laurent2023effective} \\
\hline
\end{tabular}
\end{center}
\caption{Measures $\mu$ used on simple sets.}
\label{tab:measure}
\end{table}

 \begin{theorem}[\cite{slot2111sum}]
 \label{thm: CD kernel on simple sets} 
 Let $\CX$ be either the unit ball $B_n$ or the standard simplex $\Delta_n$. For any $k \in \NN$ such that $r \geq 2(n+1)k$, there exist $1/2 \leq \lambda_i \leq 1 \ \forall i=\{0,\dots,2r\}$ 
 such that $C_{2r}[\CX,\mu,\lambda](\cdot,\by) \in \CT(\CX)_{2r}$ for any fixed $\by \in \CX$, and the associated operator $\mathbf{C}_{2r}[\CX,\mu,\lambda]$ is an invertible linear operator on $\RR[\bx]_k$ satisfying the following properties:
 \begin{align*}
    &\mathbf{C}_{2r}[\CX,\mu,\lambda](1) =1,\\
    &\mathbf{C}_{2r}[\CX,\mu,\lambda]f \in \CT(\CX)_{2r} \quad \forall\; f \in \RR[\bx]_k \;\mbox{such that}\; f(\bx) \geq 0 \ \forall \bx \in \CX,\\
    &\|\mathbf{C}_{2r}[\CX,\mu,\lambda]^{-1}f-f\|_{\CX} \leq \gamma(\CX,k)\dfrac{c(n,k)}{r^2}\|f\|_{\CX} \quad \forall\; f\in \RR[\bx]_k,
 \end{align*}
 where $c(n,k) = 2(n+1)^2k^2$, $\|f\|_{\CX}=\max_{\bx \in \CX}|f(\bx)|$,
 and $\gamma(\CX,k)$ is the harmonic constant depending on $\CX$ (which depends polynomially on $n$ for fixed $k$, polynomially on $k$ for fixed $n$).
 \end{theorem}

In the following lemma, we restate the result from \cite{slot2111sum} regarding to the parameter $\lambda$ in the perturbed CD kernel for the unit ball and the standard simplex.
    
 \begin{lemma}[\cite{slot2111sum}]\label{lemma: recall results}
     Let $\CX \subset \RR^n$ be either the unit ball $B_n$ or the standard simplex $\Delta_n$  with $\mu$ as the corresponding measure listed in Table~\ref{tab:measure}. We fix $\{P_{\alpha}(\bx)\; :\; \alpha \in \NN^n\}$ to be an orthonormal basis of $\RR[\bx]$ with respect to the inner product $\langle \cdot, \cdot \rangle_{\mu}$ induced by $\mu$ on $\CX$. Then for any $k \in \NN$ such that $r \geq 2(n+1)k$, there exists a sequence of positive numbers $\lambda = (\lambda_j)_{0 \leq j \leq 2r}$ such that
     \begin{itemize}
        \item  $\lambda_0 =1$,
         \item $C_{2r}[\CX,\mu,\lambda](\cdot,\by):= \sum_{j=0}^{2r}\lambda_j C^{(j)}[\CX,\mu](\cdot,\by) 
         \in \CT(\CX)_{2r}\quad \forall\; \by \in \CX$,
         \item $\lambda_j \in [1/2,1]\quad \forall\; 0 \leq j \leq 2r$ and\; $\sum_{j=0}^k\left|1-\dfrac{1}{\lambda_j} \right| \leq 
         \dfrac{c(n,k)}{r^2}$.
     \end{itemize}
     Here, $c(n,k) = 2(n+1)^2k^2$. In addition, the spectrum of $\BC_{2r}[\CX,\mu,\lambda]$ consists of the eigenvalues $\lambda_j$'s for 
     $0 \leq j \leq 2r$, whose eigen-space is spanned by the set of polynomials $\{P_{\alpha}(\bx)\; :\;\alpha \in \NN^n,\; |\alpha|= j\}$.
 \end{lemma}

\section{Approximation of moment sequences on a product of simple sets}

 While the SOS hierarchy approximates the optimal value $\fmin$ based on the SOS representation of the function $f(\bx)-\fmin$, the moment hierarchy 
 approximates the truncated moment sequences on a compact $\CX$ 
 by a spectrahedron. Thus, the latter approach can be analyzed independently of the objective function. In particular, the problem \eqref{POP} admits an alternative formulation as follows: for any 
 $k \geq {\rm deg}(f)$, 
 \begin{equation}\label{form of moment sequence}
    \fmin \;=\; \inf_{\iy \in \CM(\CX)} \ell_{\iy}(f) \; \left(=\sum_{\alpha \in \NN^n}f_{\alpha}\iy_{\alpha} \right) \;=\; \inf_{\iy \in \CM_k(\CX)} \ell_{\iy}(f),
\end{equation}
where $\CM_k(\CX)$ is the set of $k$-truncated moment sequences associated with the probability measures on $\CX$. The methodology for our main estimation is developed based on results for the CD kernels on $B_n$ and $\Delta_n$. The next section generalizes Theorem~\ref{thm: CD kernel on simple sets} for a product of simple sets consisting of either the unit ball or standard simplex. 
\begin{remark}
In fact we can generalize Theorem~\ref{thm: CD kernel on simple sets} 
for a product of simple sets consisting of the unit ball, standard simplex, and hypercube. However, since the version of the CD kernel for the hypercube (see e.g., \cite{laurent2023effective}) is slightly different from that of the other two sets, we only present the results on bounding the error of truncated pseudo-moment sequences and the convergence rate of the moment-SOS hierarchy for a product of unit balls and standard simplexes 
for  simplicity.
\end{remark}

\subsection{Christoffel-Darboux kernel on a product of simple sets}

In this section, we consider the domain $\CX = \prod_{i=1}^m\CX_i= \{\bx = (\bx^{(1)},\dots,\bx^{(m)}):\ \bx^{(i)} \in \CX_i \; \forall i \in [m]\}$, where each set $\CX_i \subset \RR^{n_i}$ is either the unit ball $B_{n_i}$ or the standard simplex $\Delta_{n_i}$. Based on Table~\ref{tab:measure}, we fix the probability measures $\mu_i$ on $\CX_i$, and set $\mu = \otimes_{i=1}^m \mu_i$. For these fixed measures, let $\{P^{(i)}_{\alpha}:\ \alpha \in \NN^{n_i}\}$ be an orthonormal basis on $\CX_i$ with respect to $\mu_i$. 
Let $n=\sum_{i=1}^m n_i$.
Then $\{P_{\alpha}(\bx):=\prod_{i=1}^mP^{(i)}_{\alpha_i}(\bx^{(i)})\; :\; \alpha= (\alpha_1,\dots,\alpha_m) \in \NN^n,\;\alpha_i \in \NN^{n_i} \; \forall i \in [m]\}$ is an orthonormal basis on $\CX$ with respect to $\mu$, i.e., for any $\alpha, \beta \in \NN^n$,
 \begin{eqnarray*}
    && \left \langle P_{\alpha}(\bx), P_{\beta}(\bx) \right \rangle_{\mu}=\left\langle \prod_{i=1}^mP^{(i)}_{\alpha_i}(\bx^{(i)}),\prod_{i=1}^mP^{(i)}_{\beta_i}(\bx^{(i)}) \right\rangle_{\mu}= \int_{\CX}\prod_{i=1}^mP^{(i)}_{\alpha_i}(\bx^{(i)})\prod_{i=1}^mP^{(i)}_{\beta_i}(\bx^{(i)})d\mu(\bx)
     \\
     &=&\prod_{i=1}^m\int_{\CX_i} P^{(i)}_{\alpha_i}(\bx^{(i)})P^{(i)}_{\beta_i}(\bx^{(i)})d\mu_i(\bx^{(i)})= \prod_{i=1}^m \left\langle P^{(i)}_{\alpha_i}(\bx^{(i)}),P^{(i)}_{\beta_i}(\bx^{(i)}) \right\rangle_{\mu_i}.
 \end{eqnarray*}
 
 For any $i \in [m]$ and $r \in \NN$, we denote the CD kernel on $\CX_i$ and $\CX$ by $C_{2r}[\CX_i,\mu_i]$ and $C_{2r}[\CX,\mu]$, whose associated operators are denoted by $\BC_{2r}[\CX_i,\mu_i]$ and $\BC_{2r}[\CX,\mu]$, respectively. Then, the CD kernel of degree $2r$ associated with $\mu $ is given by 
 \begin{align*}
     &C_{2r}[\CX,\mu](\bx,\by)= \sum_{j=0}^{2r}C^{(j)}[\CX,\mu](\bx,\by)
     =\sum_{j=0}^{2r}\sum_{\alpha \in \overline{\NN}^n_{j}}P_{\alpha}(\bx)P_{\alpha}(\by)
     \\
     =&\sum_{ j_1+\dots+j_m \leq 2r}
     \sum_{ \alpha_i \in \overline{\NN}^{n_i}_{j_i}, i\in [m]}
     \prod_{i=1}^mP^{(i)}_{\alpha_i}(\bx^{(i)})
     P^{(i)}_{\alpha_i}(\by^{(i)})\\
     =& \sum_{ j_1+\dots+j_m \leq 2r}\prod_{i=1}^mC^{(j_i)}[\CX_i,\mu_i](\bx^{(i)},\by^{(i)}).
 \end{align*}
 
 For each of the domains $\CX_i$ and $\CX$, the operator $\BC_{2r}$ associated with $C_{2r}$ reproduces the space of polynomials of degree at most $2r$. Similar to the technique in  \cite{slot2111sum}, we 
 modify the CD kernel partially as follows: for any sequence $\lambda= (\lambda^{(1)},\dots,\lambda^{(m)})$, where $\blambda^{(i)} := (\lambda^{(i)}_j)_{0 \leq j \leq 2r}$, we consider the  following kernel:
 \begin{equation} \label{eq-CDK}
     C_{2r}[\CX,\mu, \lambda](\bx,\by)= \sum_{ j_1+\dots+j_m \leq 2r}\; \prod_{i=1}^m\lambda^{(i)}_{j_i} C^{(j_i)}[\CX_i,\mu_i](\bx^{(i)},\by^{(i)}).
 \end{equation}
 This so-called perturbed CD kernel has a useful property that is stated in the next lemma.
 \begin{lemma}\label{lem: equivalent kernel}
     Define the following polynomial kernel
     \begin{equation*}
         K(\bx,\by) = \prod_{i=1}^mC_{2r}[\CX_i,\mu_i,\lambda^{(i)}](\bx^{(i)},\by^{(i)})
         =\prod_{i=1}^m\left(\sum_{j=0}^{2r}\lambda^{(i)}_jC^{(j)}[\CX_i,\mu_i](\bx^{(i)},\by^{(i)}) \right).
     \end{equation*}
     Then the linear operator $\BK$ associated with $K$ is identical to the operator $\BC_{2r}[\CX,\mu,\lambda]$ on $\RR[\bx]_{2r}$.
\end{lemma}
    \begin{proof}
     The orthogonality of $\{P^{(i)}_{\alpha^{(i)}}:\ \alpha^{ (i)} \in \NN^{n_i}\}$ implies that for any $P_{\alpha}(\bx)$ satisfying $\alpha=(\alpha^{(1)},\dots,\alpha^{(m)}) \in \NN^n,\ |\alpha^{(i)}|=j_i,\ \sum_{i=1}^mj_i \leq 2r$, the following identities hold:
         \begin{align*}
            &\BK P_{\alpha}(\bx)\;=\; \int_{\CX}\left(\prod_{i=1}^mC_{2r}[\CX_i,\mu_i,\lambda^{(i)}](\bx^{(i)},\overline{\bx}^{(i)})\right) \cdot \left( \prod_{i=1}^mP^{(i)}_{\alpha^{(i)}}(\overline{\bx}^{(i)})\right)d\mu(\overline{\bx})\\
            =\;& \prod_{i=1}^m\left(\int_{\CX_i}C_{2r}[\CX_i,\mu_i,\lambda^{(i)}](\bx^{(i)},\overline{\bx}^{(i)})P^{(i)}_{\alpha^{(i)}}(\overline{\bx}^{(i)} )d\mu_i(\overline{\bx}^{(i)})\right)\\
            =\;& \prod_{i=1}^m \left(\lambda^{(i)}_{j_i}P^{(i)}_{\alpha^{(i)}}(\bx^{(i)}) \right)= \left(\prod_{i=1}^m \lambda^{(i)}_{j_i}\right)P_{\alpha}(\bx),
         \end{align*} 
         and 
         \begin{align*}
             &\mathbf{C}_{2r}[\CX,\mu,\lambda]P_{\alpha}(\bx)
             =\int_{\CX}C_{2r}[\CX,\mu,\lambda](\bx,\overline{\bx})
             P_{\alpha}(\overline{\bx})d\mu(\overline{\bx})\\
             =&\sum_{j^\prime_1+\dots+j^\prime_m\leq 2r}
             \int_{\CX}
             \left( \prod_{i=1}^m \lambda^{(i)}_{j_i^\prime}C^{(j^\prime_i)}[\CX_i,\mu_i](\bx^{(i)},\overline{\bx}^{(i)}) 
             P^{(i)}_{\alpha^{(i)}}(\overline{\bx}^{(i)}) \right)
             d\mu(\overline{\bx}) 
             \quad 
             \\
             =&\sum_{j^\prime_1+\dots+j_m ^\prime\leq 2r}\prod_{i=1}^m\left(\int_{\CX_i}\lambda^{(i)}_{j_i^\prime}C^{(j^\prime_i)}[\CX_i,\mu_i] (\bx^{(i)},\overline{\bx}^{(i)})P^{(i)}_{\alpha^{(i)}}(\overline{\bx}^{(i)})d\mu_i(\overline{\bx}^{(i)})\right)\\
             =&\prod_{i=1}^m\left(\int_{\CX_i}\lambda^{(i)}_{j_i}C^{(j_i)}[\CX_i,\mu_i] (\bx^{(i)},\overline{\bx}^{(i)})P^{(i)}_{\alpha^{(i)}}(\overline{\bx}^{(i)})d\mu_i(\overline{\bx}^{(i)})\right)\quad \mbox{(since $|\alpha^{(i)}| = j_i$)} 
             \\
             =& \prod_{i=1}^m \left(\lambda^{(i)}_{j_i}P^{(i)}_{\alpha^{(i)}}(\bx^{(i)}) \right)= \left(\prod_{i=1}^m \lambda^{(i)}_{j_i}\right)P_{\alpha}(\bx). 
         \end{align*}
         Hence, both linear operators $\BK$ and $\BC_{2r}[\CX,\mu,\lambda]$ share the same spectrum on $\RR[\bx]_{2r}$ with the following properties:
     \begin{itemize}
         \item For any $0 \leq j_1,\dots,j_m \leq 2r$ such that $\sum_{i=1}^mj_i \leq 2r$, $\prod_{i=1}^m\lambda^{(i)}_{j_i}$ is an eigenvalue with multiplicity equal to the number of polynomials of the form $P_{\alpha}(\bx)$ with $\alpha=(\alpha^{(1)},\dots,\alpha^{(m)})$ such that $ |\alpha^{(i)}| =j_i$, i.e., the multiplicity of $\prod_{i=1}^m\lambda^{(i)}_{j_i}$ is $\frac{|\alpha|!}{\prod_{i=1}^m|\alpha^{(i)}|!}$.
         \item The definitions of $\BK$ and $\BC_{2r}[\CX,\mu,\lambda]$ imply that 
         the eigenspace for $\prod_{i=1}^m\lambda^{(i)}_{j_i}$ is 
         \begin{equation}
             S_{j_1,\dots, j_m}:= {\rm span}\Big(\Big\{P_{\alpha}(\bx):\ \alpha=(\alpha^{(1)},\dots,\alpha^{(m)}) \in \NN^n,\ |\alpha^{(i)}|=j_i \; \forall i \in [m] \Big\}\Big).
             \label{eq-S}
         \end{equation}
    \end{itemize}
    This completes the proof.
    \end{proof}
    
Lemma~\ref{lem: equivalent kernel} also shows that 
\begin{displaymath}
    \RR[\bx]_{2r} = \bigoplus_{j_1 + \dots + j_m \leq 2r}S_{j_1,\dots,j_m}.
\end{displaymath}
Then we can decompose any polynomial $p \in \RR[\bx]_{2r}$ as 
    \begin{displaymath}
         p(\bx) = \sum_{j_1+\dots+j_m  \leq 2r} p_{j_1,\dots,j_m}(\bx),\quad p_{j_1,\dots,j_m} \in S_{j_1,\dots,j_m}\;\; 
         \forall\; j_1+\dots+j_m \leq 2r.  
    \end{displaymath}
    By the compactness of $\CX$, we can define the following harmonic constant for $\CX$ and any $k\in\NN$: 
    \begin{eqnarray}
        \Lambda(\CX,k):= \underset{p \in \RR[\bx]_{k}}{\max}\; \underset{0 \leq j_1+\dots+j_m \leq k}{\max}\ \dfrac{\|p_{j_1,\dots,j_m}\|_{\CX}}{\|p\|_{\CX}}.
        \label{eq-Lambda}
    \end{eqnarray}
    The constant $\Lambda(\CX,k)$ depends only on $m,n$ and $k$ (see e.g., \cite{fang2021sum,slot2111sum}) and plays the role of the harmonic constant 
    in Theorem~\ref{thm: CD kernel on simple sets}. The quantitative analysis on bounding this harmonic constant 
    is presented in Appendix~\ref{appendix: harmonic constant}.

 We  now extend Theorem~\ref{thm: CD kernel on simple sets} to any product of unit balls and standard simplexes in the following theorem.
 
 \begin{theorem}\label{thm: CD on product}
     Let $\CX = \prod_{i=1}^m \CX_i$ where each $\CX_i\subset \mathbb{R}^{n_i}$
     is either the unit ball $B_{n_i}$ or the standard simplex $\Delta_{n_i}$. We fix a positive integer $k$ and consider any $r > 2(\max\{n_1,\dots,n_m\}+1)k$. Let $\lambda= (\lambda^{(1)},\dots,\lambda^{(m)})$, 
     where $\lambda^{(i)}$ is  the sequence associated with the simple set $\CX_i$ satisfying the conditions as in Lemma~\ref{lemma: recall results}. Then the perturbed CD kernel $C_{2r}[\CX,\mu,\lambda]$ and its associated operator $\BC_{2r}[\CX,\mu,\lambda]$ on $\RR[\bx]_k$ satisfies the following conditions:
     \begin{align}
        &\mathbf{C}_{2r}[\CX,\mu,\lambda](1) =1,\tag{P1}\label{P1}\\
        &\mathbf{C}_{2r}[\CX,\mu,\lambda]f \in \CT(\CX)_{2mr} \quad \forall f \in \RR[\bx]_k \; \mbox{such that}\;
        f(\bx) \geq 0 \ \forall \bx \in \CX,\tag{P2}\label{P2}\\
        &\|\mathbf{C}_{2r}^{-1}f-f\|_{\CX} \leq 
          2^{m-1}\binom{k+m-1}{m-1}\left(\sum_{i=1}^m c(n_i,k)\right)\frac{\Lambda(\CX,k) \|f\|_{\CX}}{r^2}, \ \forall \; f\in \RR[\bx]_k,
        \tag{P3}\label{P3}
     \end{align}
     where $\Lambda(\CX,k)$ is the harmonic constant defined in 
     \eqref{eq-Lambda}.
 \end{theorem}
 \begin{proof}
    For any $r > 2(\max\{n_1,\ldots,n_m\}+1)k$, we define $\BK$ to be the linear operator on $\RR[\bx]_k$ associated to the kernel $K$ defined in Lemma \ref{lem: equivalent kernel}. Since $2r \geq k$, Lemma~\ref{lem: equivalent kernel} implies that $\BK$ and $\BC_{2r}[\CX,\mu,\lambda]$ are identical on $\RR[\bx]_{k}$. Thus, it suffices to prove that the linear operator $\BK: \RR[\bx]_k \to \RR[\bx]_k$ satisfies \eqref{P1},\eqref{P2} and \eqref{P3}. We recall from Lemma~\ref{lemma: recall results} that $\lambda^{(i)}_0 =1\; \forall i \in [m]$, 
    and $1 = \prod_{i=1}^m P^{(i)}_0(\bx^{(i)})$ is the eigen-polynomial of the eigenvalue $\prod_{i=1}^m\lambda^{(i)}_0=1$. This implies that 
         \begin{displaymath}
             \BK(1) = 1\quad \mbox{and hence} \quad \eqref{P1} \text{ is satisfied.}
         \end{displaymath}
         We next prove \eqref{P2} by utilizing Tchakaloff's theorem (see e.g., \cite{Tchakaloff}) to show the existence of a cubature rule for the integration of polynomials, i.e., there exist $\{(\bx_j,w_j): 1 \leq j \leq N := k + {\rm deg} (K) \} \subset \CX \times \RR_{>0}$ 
         such that for any non-negative polynomial $f \in \RR[\bx]_k$ on $\CX$, we have that 
         \begin{displaymath}
             \BK f(\bx) = \int_{\CX}\BK(\bx,\overline{\bx})f(\overline{\bx})d\mu (\overline{\bx})=\sum_{j=1}^NK(\bx,\bx_j)w_jf(\bx_j).
         \end{displaymath}
    Recall that $K(\bx,\bx_j)= \prod_{i=1}^mC_{2r}[\CX_i,\mu_i,\lambda^{(i)}](\bx^{(i)},\bx^{(i)}_j)$. Lemma~\ref{lemma: recall results} gives $\lambda^{(i)}$'s satisfying that
     \begin{displaymath}
         C_{2r}[\CX_i,\mu_i,\lambda^{(i)}](\bx^{(i)},\bx^{(i)}_j) \in \CT(\CX_i)_{2r}\quad \forall\; i \in [m],\;\; j\in[N].
     \end{displaymath}
     Therefore, $K(\bx,\bx_j) = \prod_{i=1}^mC_{2r}[\CX_i,\mu_i,\lambda^{(i)}](\bx^{(i)},\bx^{(i)}_j) \in \CT(\CX)_{2mr}$.
     Combining this with the condition $w_jf(\bx_j) \geq 0\ \forall i \in [N]$, we obtain property \eqref{P2} as
     follows: 
     \begin{displaymath}
         \BK f(\bx) =\sum_{j=1}^N K(\bx,\bx_j)w_j f(\bx_j) \in \CT(\CX)_{2mr}. 
     \end{displaymath}
     
     Finally, we prove property (P3). Since $\lambda^{(i)}_j$'s are positive numbers, the linear operator $\BK: \RR[\bx]_{2r} \to \RR[\bx]_{2r}$ has all its eigenvalues being positive. Thus, the inverse $\BK^{-1}$ exists and the spectrum is
     \begin{displaymath}
         \left\{\left(\prod_{i=1}^m\lambda^{(i)}_{j_i}\right)^{-1}\;:\: j_1+\dots j_m \leq 2r \right\}.
     \end{displaymath}
    Then, we have 
    \begin{eqnarray}\label{equation 1}
        |\BK^{-1}f(\bx) -f(\bx)| &=& 
        \left|\sum_{j_1+\dots+j_m \leq k}\left(1-\dfrac{1}{\prod_{i=1}^m\lambda^{(i)}_{j_i}} \right)f_{j_1,\dots,j_m}(\bx) \right| \nonumber \\ 
        &\leq& \left( \sum_{j_1\dots+j_m \leq k}\left|1-\dfrac{1}{\prod_{i=1}^m\lambda^{(i)}_{j_i}} \right| \right) \Lambda(\CX,k) \|f\|_{\CX}.
    \end{eqnarray}
     Taking maximum on both sides of \eqref{equation 1} over $\bx \in \CX$, we obtain that 
    \begin{displaymath}
        \|\BK^{-1}f -f\|_{\CX} \leq \Lambda(\CX,k) \|f\|_{\CX} \sum_{j_1+\dots+j_m \leq k}\left|1-\dfrac{1}{\prod_{i=1}^m\lambda^{(i)}_{j_i}} \right|.
    \end{displaymath}
    Based on Lemma~\ref{lemma: recall results}, we know that for all $i \in [m]$ and $0 \leq j \leq 2r $, we have $\lambda^{(i)}_j \in [1/2,1]$. Then 
    \begin{eqnarray}
        && \hspace{-7mm}
        \sum_{j_1+\dots+j_m \leq k}\left|1-\dfrac{1}{\prod_{i=1}^m\lambda^{(i)}_{j_i}} \right|= \sum_{j_1+\dots+j_m \leq k}\frac{1- \prod_{i=1}^m\lambda^{(i)}_{j_i}}{\prod_{i=1}^m\lambda^{(i)}_{j_i}}\leq \sum_{j_1+\dots+j_m \leq k}\frac{\sum_{i=1}^m (1- \lambda^{(i)}_{j_i})}{\prod_{i=1}^m\lambda^{(i)}_{j_i}}
        \nonumber \\
        &\leq& 2^{m-1}\sum_{j_1+\dots+j_m \leq k}\sum_{i=1}^m\Big|1-\frac{1}{\lambda^{(i)}_{j_i}} \Big|\leq 2^{m-1}\binom{k+m-1}{m-1}\sum_{i=1}^m\sum_{j=1}^k\Big|1-\frac{1}{\lambda^{(i)}_{j}} \Big|
        \nonumber \\
        &\leq &  2^{m-1}\binom{k+m-1}{m-1}\left(\sum_{i=1}^m c(n_i,k)\right)\frac{1}{r^2}.
        \label{eq-lambda-bound}
    \end{eqnarray}
    Here, the coefficient $2^{m-1}$ comes from the fact that $\lambda^{(i)}_j \in [1/2,1]$, the coefficient $\binom{k+m-1}{m-1}$ is the result of counting the tuple $(j_1,\dots,j_{i-1},j_{i+1},\dots,j_m)$ such that $\sum_{i=1}^mj_i \leq k$, and the last inequality is based on Lemma~\ref{lemma: recall results}.
    In short, we can obtain the following inequality:
    \begin{displaymath}
        \|\BK^{-1}f -f\|_{\CX} \leq 2^{m-1}\binom{k+m-1}{m-1}\left(\sum_{i=1}^m c(n_i,k)\right)\frac{\Lambda(\CX,k) \|f\|_{\CX}}{r^2},
    \end{displaymath}
    which proved \eqref{P3} and this completes the proof. 
 \end{proof}
 
 \begin{corollary}\label{cor: rate on product}
    Consider the polynomial optimization problem 
    \begin{displaymath}
    \fmin := \underset{\bx \in \CX}{\min}\ f(\bx),
    \end{displaymath}
    where $\CX$ is a product of unit balls and standard simplexes as in Theorem 
    \ref{thm: CD on product}. Let  $k = {\rm deg} (f)$ and
    $\fmax = \max \{f(\bx): \bx\in\CX\}$. Then we have that 
    the convergence rate of the Schmüdgen-type moment-SOS hierarchies \eqref{hierarchy: SOS Schmudegen} and \eqref{hierarchy: moment Schmudgen} is $\mathrm{O}(1/r^2)$ with the following bound for any $r \geq 2m(\max\{n_1,\ldots,n_m\}+1)k+m$, 
    \begin{eqnarray*}
       0 &\leq& \fmin - \mlb(f,\CT(\CX))_{r} \leq \fmin - \lb(f,\CT(\CX))_{r}
       \;\leq\;
       \Gamma(\CX,k)\frac{(\fmax-\fmin)}{(r-m)^2},
    \end{eqnarray*}
    where the constant
$\Gamma(\CX,k)= m^2\, 2^{m-1} \binom{k+m-1}{m-1}\left(\sum_{i=1}^m c(n_i,k)\right)\Lambda(\CX,k)$
is dependent on the domain $\CX$ and the degree $k$.
 \end{corollary}
 \begin{proof}
     We only need to prove the last inequality. Let $r' = \lfloor r/m \rfloor$. Then $r' > 2(\max\{n_1,\dots,n_m\}+1)k$. Consider the non-negative polynomial $f(\bx) -\fmin$ on $\CX$. By applying Theorem~\ref{thm: CD on product} to $f(\bx)-\fmin$, we can choose the parameter $\lambda= (\lambda^{(1)},\dots,\lambda^{(m)})$ such that the polynomial
     \begin{displaymath}
         \BC_{2r'}^{-1}f(\bx)-\fmin+ \varepsilon \geq 0\;\; \forall \; \bx \in \CX, 
    \end{displaymath}
    where $\varepsilon:=  2^{m-1}\binom{k+m-1}{m-1}\big(\sum_{i=1}^m c(n_i,k)\big)\frac{\Lambda(\CX,k)(\fmax-\fmin)}{(r')^2}$. Hence, \eqref{P2} and \eqref{P1} in Proposition~\ref{thm: CD on product} imply that 
     \begin{align*}
         &\BC_{2r'}(\BC_{2r'}^{-1}f(\bx)-\fmin+ \varepsilon) 
         \in \CT(\CX)_{2mr'} \subset \CT(\CX)_{2r}
         \\
         \Rightarrow\;\; & f(\bx) -\fmin + \varepsilon \in \CT(\CX)_{2r}\\
         \Rightarrow\;\; & \fmin - \lb(f,\CT(\CX))_r \leq  
         \varepsilon = 2^{m-1}\binom{k+m-1}{m-1}\left(\sum_{i=1}^m c(n_i,k)\right)\Lambda(\CX,k)\frac{\fmax-\fmin}{(r')^2}
         \\
         &\hspace{33mm}
         \leq 2^{m-1}\binom{k+m-1}{m-1}\left(\sum_{i=1}^m c(n_i,k)\right)\Lambda(\CX,k)\frac{(\fmax-\fmin)m^2}{(r-m)^2}.
     \end{align*}
From here, we get the required result. 
 \end{proof}

\subsection{Tightness of  truncated pseudo-moment sequences on a product of simple sets}

In this section, we use the results developed in the previous section to evaluate the tightness of $\CM_k(\CT(\CX)_{2r})$ in approximating 
$\CM_k(\CX)$. In particular, we show that the Hausdorff distance $\bd_k(\CT(\CX)_{2r})$ is $\mathrm{O}(1/r^2)$ for a product of simple sets $\CX$. 
We remind the reader that for simplicity, the simple sets considered in this paper are either a unit ball or a standard simplex. (With appropriate modifications in the analysis, we can allow the simple sets to also include the hypercube). 

We fix a positive integer $k$ (throughout this section, we always assume $k=2l$ to be even for convenience). The set $\CM_k(\CX)$ consists of $k$-truncated moment sequences $\iy:=(\iy_{\alpha})_{\alpha \in \NN^n_k}$.
We know that $\CM_k(\CT(\CX)_{2r})$ is an outer approximation of $\CM_k(\CX)$. As a  convention, we always assume that $2r \geq k$ so that $s(n,2r)\geq s(n,k)$, and the projection onto the  first  $s(n,k)$ coordinates of the
feasible vectors of $\CM(\CT(\CX)_{2r})$ is well-defined. For any 
$k$--truncated moment sequence $\iy \in \CM_k(\CX)$, there exists a probability measure $\kappa \in \CP(\CX)$ corresponding to a moment sequence, whose projection onto the first $s(n,k)$ coordinates is $\iy$. Moreover, because of Tchakaloff's theorem (see e.g.,\cite{Tchakerloff's_theorem}), there exist at most $N:=s(n,k)$ points $\{\bx_j: j\in[N]\} \subset \mbox{supp}(\kappa)$ and corresponding positive weight $\{w_j: j\in [N]\}$ satisfying $\sum_{j=1}^N w_j=1$ such that 
for any $f \in \RR[\bx]_k$, the integral of $f$ over $\CX$ with respect to $\kappa$ can be calculated as follows: 
\begin{displaymath}
    \int_{\CX}f(\bx) d\kappa(\bx) = \sum_{j=1}^{N}w_j f(\bx_j).
\end{displaymath}
This implies that when considering an element $y\in \CM_k(\CX)$, its corresponding measure can always be assumed to be a discrete measure with support contained in $\CX$. Thus, the 
$l$--truncated moment matrix admits the following decomposition 
$\forall \,\ g \in \RR[\bx]$ such that 
$2t + {\rm deg} (g)\leq k=2l$:
\begin{equation*}
    \bM_l(\iy) = \sum_{j=1}^Nw_j \bv_l(\bx_j)\bv_l(\bx_j)^{\top},\ \mbox{ and } \ \bM_t(g\iy) = \sum_{j=1}^Nw_jg(\bx_j)\bv_t(\bx_j)\bv_t(\bx_j)^{\top}.
\end{equation*}

The following theorem represents the main result of this section on the error estimation of the truncated pseudo-moment sequences in $\CM(\CT(\CX)_{2r}).$

\begin{theorem}\label{thm: error on product of simple sets}
    Let $\CX\subset \RR^n$ be a compact set that is a product $\Pi_{i=1}^m\CX_i$ where each $\CX_i\subset \RR^{n_i}$ is either the unit ball $B_{n_i}$ or standard simplex $\Delta_{n_i}$. We assume that there exists $R$ such that $\CX \subset \mathbb{B}_R$, and the inequality, $R^2 -\|\bx\|^2 \geq0$, is included in the definition of $\CX$. For a fixed $k=2l$ and $r \geq 2m(\max\{n_1,\dots,n_m\}+1)k+m$, the Hausdorff distance $\bd_k(\CT(\CX)_{2r})$ admits the following upper bound:
    \begin{equation}\label{error on moment sequence on product of simple sets}
        \bd_k(\CT(\CX)_{2r}) \leq \Gamma(\CX,k)\dfrac{2\gamma(R,n,k)}{(r-m)^2}.
    \end{equation}
    Here, the parameter $\gamma(R,n,k)$ is the radius of the ball centered at the origin that
    contains $\CM_k(\CT(\mathbb{B}_R)_{2r})$, and it depends polynomially on $n$ and $k$.
\end{theorem}
\begin{proof}
    Since both $\CM_k(\CT(\CX)_{2r})$ and $\CM_k(\CX)$ are compact (this is elaborated in Remark~\ref{rem: compactness of pseudo-moment sequence} below), we can let $\overline{\iy} \in \CM_k(\CT(\CX)_{2r})$ be a $k-$truncated pseudo-moment sequence such that 
    \begin{equation*}
        \bd_k(\CT(\CX)_{2r}) := \bd( \CM_k(\CT(\CX)_{2r}),\CM_k(\CX)) = \bd(\overline{\iy}, \CM_k(\CX)).
    \end{equation*}
    Because $\CX$ is compact and we have argued above by Tchakaloff's theorem that every $k$--truncated moment sequence is a convex combination of $\{\bv_k(\bx_1),\ldots,\bv_k(\bx_N)\} $ with $\bx_i\in\CX$ for all $i\in [N]$ and $N = s(n,k)$, we know that $\CM_k(\CX)$ is a compact convex set. Hence, there exists the unique projection of $\overline{\iy}$ on $\CM_k(\CX)$, denoted by $\Tilde{\iy}$, such that 
    \begin{displaymath}
      \bd_k(\CT(\CX)_{2r})^2 = \|\Tilde{\iy} - \overline{\iy}\|^2 =   \min\{\|\iy - \overline{\iy}\|^2\,:\, \iy \in \CM_k(\CX)\}.
    \end{displaymath}
    
    According to the first-order optimality condition, $\langle \Tilde{\iy}-\overline{\iy}, \iy - \Tilde{\iy} \rangle \geq 0$ for any $\iy \in \CM_k(\CX)$, and we obtain
    \begin{equation*}
        L(\iy):=\langle \Tilde{\iy}-\overline{\iy}, \iy - \overline{\iy} \rangle=\langle \Tilde{\iy}-\overline{\iy}, \iy - \Tilde{\iy} \rangle +\langle \Tilde{\iy}-\overline{\iy}, \Tilde{\iy} - \overline{\iy}\rangle\geq \bd_k(\CT(\CX)_{2r})^2 \;\; \forall\; \iy \in \CM_k(\CX).
    \end{equation*}
        Hence, $\Tilde{\iy}$ is also a minimizer of $L(\iy)$ on $\CM_k(\CX)$ with the minimum value $\bd_k(\CT(\CX)_{2r})^2$. 
        
        Consider the following problem
    \begin{align}\label{SDP: dual form}
        \underset{\iy \in \CM_k(\CX)}{\min}\ L(\iy) = \sum_{\alpha \in \NN^n_k}(\Tilde{\iy}_{\alpha}-\overline{\iy}_{\alpha})\iy_{\alpha}-\langle \Tilde{\iy}-\overline{\iy}, \overline{\iy} \rangle.
    \end{align}
    This problem is actually the equivalent form of the following POP:
    \begin{equation}\label{sub POP}
        \min  \Bigl\{\ f(\bx) = \sum_{\alpha \in \NN^n_k}(\Tilde{\iy}_{\alpha}-\overline{\iy}_{\alpha})\bx^{\alpha}-\langle \Tilde{\iy}-\overline{\iy}, \overline{\iy}
        \rangle
        \,:\, \bx \in \CX\Bigr\}.
    \end{equation}
    In addition, there exists a positive number $\gamma(R,n,k)$ such that for any $2r \geq k$, all the $k-$truncated pseudo-moment sequences of $\CM_k(\CT(\mathbb{B}_R)_{2r})$
    are contained in the ball centered at the origin with radius ${\gamma(R,n,k)}$ (see Remark~\ref{rem: compactness of pseudo-moment sequence}). Thus for any $\iy \in \CM_k(\CX)\subset \CM_k(\CT(\CX)_{2r}) \subset \CM_k(\CT(\mathbb{B}_R)_{2r})$, the Cauchy–Schwarz inequality implies that 
    \begin{equation*}
        L(\iy)=\langle \Tilde{\iy}-\overline{\iy}, \iy - \overline{\iy} \rangle \leq \|\Tilde{\iy}-\overline{\iy}\|
        \| \iy - \overline{\iy}\|
        \leq 
        \bd_k(\CT(\CX)_{2r})(2 \gamma(R,n,k)).
    \end{equation*}
    
    Consider the problem \eqref{sub POP} with the conditions $\fmin = \min_{\iy \in \CM_k(\CX)} L(\iy)=\bd_k(\CT(\CX)_{2r})^2$ and $\fmax \leq 2\gamma(R,n,k) \bd_k(\CT(\CX)_{2r})$.
    We apply the $r$-th level of the Schmüdgen-type moment hierarchy \eqref{hierarchy: moment Schmudgen} for the problem \eqref{sub POP} to get
    \begin{displaymath}
        \mlb(f,\CT(\CX))_r=\underset{\iy \in \CM_k(\CT(\CX)_{2r})}{\min}\ L(\iy) = \langle \Tilde{\iy}-\overline{\iy}, \iy - \overline{\iy} \rangle,
    \end{displaymath}
    whose error can be upper bounded by using Corollary~\ref{cor: rate on product} under the degree condition $r \geq 2m(\max\{n_1,\dots,n_m\}+1)k+m$ as
    follows: 
    \begin{equation}\label{inequality 2}
         \fmin- \mlb(f,\CT(\CX))_r\leq \Gamma(\CX,k)\dfrac{\fmax-\fmin}{(r-m)^2} \leq \Gamma(\CX,k)\dfrac{2\gamma(R,n,k) \bd_k(\CT(\CX)_{2r})}{(r-m)^2}.
    \end{equation}

 
    Since $\overline{\iy} \in \CM_k(\CT(\CX)_{2r})$, and $L(\overline{\iy}) =0$, we get $\mlb(f,\CT(\CX))_{r} \leq 0$. From this, the inequality \eqref{inequality 2} implies that 
    \begin{displaymath}
        \bd_k(\CT(\CX)_{2r})^2 = \fmin \leq \Gamma(\CX,k)\dfrac{2\gamma(R,n,k) \bd_k(\CT(\CX)_{2r})}{(r-m)^2}\quad \Rightarrow \quad \bd_k(\CT(\CX)_{2r}) \leq \Gamma(\CX,k)\dfrac{2\gamma(R,n,k) }{(r-m)^2}.
    \end{displaymath}
    This completes the proof.
\end{proof}
\begin{remark}\label{rem: compactness of pseudo-moment sequence}
    From now on, we will frequently use the parameter $\gamma(R,n,k)$, where $n$ denotes the dimension of variable $\bx$, $k$ denotes the truncation order of the moment sequences, and $R$ is the radius such that $\CX \subset \mathbb{B}_R$. Since we assume that the inequality $R- \|\bx\|^2\geq 0$ is included in the description of $\CX$, therefore for any positive even integer $k=2l$, we always have 
    $\CM(\CT(\CX)_k) \subset \CM(\CT(\mathbb{B}_R)_k)$.
    Thus, $\gamma(R,n,k)$ can be used to bound the Euclidean norm of the pseudo-moment sequences in $\CM(\CT(\CX)_k)$.
    The explicit expression of $\gamma(R,n,k)$ is given in 
    the following lemma.
    \begin{lemma}[\cite{josz2016strong}, Lemma 3]\label{lemm: radius of preordering}
    For $k=2l,$
        $\CM(\CT(\mathbb{B}_R)_k)$ is contained in the Euclidean ball 
        centered at the origin with radius 
        \begin{displaymath}
            \gamma(R,n,k):=\sqrt{\binom{n+l}{n}}\sum_{i=0}^lR^{2i}\,.
        \end{displaymath}
    \end{lemma}
\end{remark}

\begin{remark}
   Theorem~\ref{thm: error on product of simple sets} can be extended to any product $\CX$ of the ball $\mathbb{B}_R$ and the simplex $\Delta_K^n := \{ \bx \in \RR^n:\; x_i \geq 0 \ \forall i \in [n],\; \sum_{i=1}^n x_i \leq K\}$ for any positive numbers $R$ and $K$. Additionally, the CD kernel over these sets are similarly defined as in Table~\ref{tab:measure} by scaling. Furthermore, the convergence rate and the error on the truncated pseudo-moment sequences can be obtained
   based on Theorem~\ref{thm: error on product of simple sets}  via an invertible linear transformation, which is presented in Appendix~\ref{appendix: linear transformation}. From now on, we refer the term "simple sets" to the ball $\mathbb{B}_R$ and the simplex $\Delta^n_K$.
   For any product $\CX$ of simple sets $\CX_i$, Theorem~\ref{thm: error on product of simple sets} remains valid as 
   \begin{equation*}
        \bd_k(\CT(\CX)_{2r}) \leq \Gamma(\CX,k)\dfrac{2\gamma(R,n,k)}{r^2}
    \end{equation*}
    where $\Gamma(\CX,k)$ is specified in Remark~\ref{rem: scalling}.
\end{remark}

\subsection{Convergence rate of the Schm\"udgen-type hierarchy of upper-bounds over a product of simple sets}

Theorem~\ref{thm: CD on product} is an extension of the result established in the paper \cite{slot2111sum}. Here we reuse the analysis in \cite[Section 5]{slot2111sum} to establish the same convergence rate of 
$\mathrm{O}(1/r^2)$ for the Schm\"udgen-type hierarchy of upper-bounds over the product $\CX$ of unit balls and standard simplexes. In other words, this section is devoted to proving the following theorem.

\begin{theorem}\label{thm: convergence rate for hierarchy of upper bounds}
    Let $\CX = \prod_{i=1}^m\CX_i$, where each $\CX_i \subset \RR^{n_i}$ is either the unit ball $B_{n_i}$ or the standard simplex $\Delta_{n_i}$.
       We choose the reference measure $\mu =\otimes_{i=1}^m\mu_i$ for $\CX$, where $\mu_i$ is the reference measure on $\CX_i$ as in Table~\ref{tab:measure}. 
       Let $n =\sum_{i=1}^m n_i$. For any $i \in [m]$, let $\iy^{(i)}=\big(\iy^{(i)}_{\alpha^{(i)}\in \NN^{n_i}}\big)$ be the moment sequence with respect to $\mu_i$. 
    Then the moment sequence $\iy=(\iy_{\alpha\in \NN^n})$ with respect to $\mu$ is defined as follows: for any $\alpha= (\alpha^{(1)},\dots,\alpha^{(m)}) \in \NN^n$ with $\alpha^{(i)} \in \NN^n_i\; \forall i \in [m]$, the $\alpha-$component of $\iy$ is defined by 
    \begin{displaymath}
        \iy_{\alpha}= \prod_{i=1}^m \iy^{(i)}_{\alpha^{(i)}}.
    \end{displaymath}
    In addition, for a given $f(\bx)\in \RR[\bx]$ of degree $k$, consider the POP 
    \begin{equation*}
    \fmin := \min_{\bx \in \CX}f(\bx).
    \end{equation*}
   Then for any integer $r > 2(\max\{n_1,\ldots,n_m\}+1)k$,
   the following inequality holds: 
    \begin{displaymath}
        \ub(f,\CT(\CX))_{mr}-\fmin \leq \binom{k+m-1}{m-1}\left(\sum_{i=1}^m c(n_i,k)\right)\frac{\Lambda(\CX,k)\|f\|_{\CX}}{r^2}.
    \end{displaymath}
    Consequently, the convergence rate of the Schm\"udgen-type hierarchy of upper-bounds over the product of unit balls and standard simplexes is
     $\mathrm{O}(1/r^2)$.
\end{theorem}

\begin{proof}
To analyze the convergence rate of the Schm\"udgen-type hierarchy of upper-bounds, we use the same evaluation method in \cite[Section 5]{slot2111sum}. For a fixed positive integer $r$, let $\bx^*$ be one of the minimizers of $f$ over $\CX$, whose existence is due to the compactness of $\CX$. We recall the perturbed CD kernel $C_{2r}[\CX,\mu,\lambda]$ from Theorem~\ref{thm: CD on product} and set 
\begin{equation*}
    \sigma(\bx) = C_{2r}[\CX,\mu,\lambda](\bx,\bx^*).
\end{equation*}
In the proof of Theorem~\ref{thm: CD on product}, we have pointed out that $\sigma(\bx):= C_{2r}[\CX,\mu,\lambda](\bx,\bx^*)\in \CT(\CX)_{2mr}$. Whence, we have 
\begin{equation*}
    \fmin \leq \ub(f,\CT(\CX))_{mr} \leq \int_{\CX}f(\bx)\sigma(\bx)d\mu(\bx)=\BC_{2r}[\CX,\mu,\lambda]f(\bx^*).
\end{equation*}
Note that $f(\bx)$ admits the following decomposition
\begin{displaymath}
         f(\bx) = \sum_{j_1+\dots+j_m  \leq k}f_{j_1,\dots,j_m}(\bx),\quad f_{j_1,\dots,j_m} \in S_{j_1,\dots,j_m}\ \forall j_1+\dots+j_m \leq k.     
\end{displaymath}
Thus, the image of $f$ under $\BC_{2r}[\CX,\mu,\lambda]$ is 
\begin{equation*}
    \BC_{2r}[\CX,\mu,\lambda]f(\bx) = \sum_{j_1+\dots+j_m  \leq k} \Big(\prod_{i=1}^m\lambda^{(i)}_{j_i}\Big) 
    f_{j_1,\dots,j_m}(\bx).
\end{equation*}
Therefore we can bound 
\begin{align*}
    &\ub(f,\CT(\CX))_{mr}-\fmin \leq \BC_{2r}[\CX,\mu,\lambda]f(\bx^*)-f(\bx^*)\\
        \leq\; & \sum_{j_1+\dots+j_m  \leq k}
        \left|\left(1-\prod_{i=1}^m\lambda^{(i)}_{j_i}\right)f_{j_1,\dots,j_m}(\bx^*)\right|
        \quad 
        \\
        \leq\; & \left(\sum_{j_1+\dots+j_m  \leq k}\sum_{i=1}^m\left|1-\lambda^{(i)}_{j_i}\right|\right) \Lambda(\CX,k)\|f\|_{\CX}
        \quad \mbox{(by the definition of $\Lambda(\CX,k)$)}\\
        \leq\; & \binom{k+m-1}{m-1}\left(\sum_{i=1}^m\sum_{j=1}^k\left|1-\lambda^{(i)}_{j}\right|\right)
        \Lambda(\CX,k)\|f\|_{\CX} \\
        \leq\;&\binom{k+m-1}{m-1}\left(\sum_{i=1}^m\sum_{j=1}^k\Big|1-\frac{1}{\lambda^{(i)}_{j}} \Big|\right)
        \Lambda(\CX,k)\|f\|_{\CX}\quad \text{(since $\lambda^{(i)}_j \in [{1}/{2},1]$)}\\
        \leq\;& \binom{k+m-1}{m-1}\left(\sum_{i=1}^m c(n_i,k)\right)\frac{\Lambda(\CX,k)\|f\|_{\CX}}{r^2},
\end{align*}
where the last inequality follows from Lemma~\ref{lemma: recall results}, which leads to our desired bound. 
\end{proof}

\section{Error estimation of truncated pseudo-moment sequences on a compact semi-algebraic set}\label{section: error in algebraic set}

 We now extend the error estimation of truncated pseudo-moment sequences to  a compact domain $\CX$ defined as in \eqref{domain X}.
 The methodology starts by estimating the error of the truncated pseudo-moment sequence on a product of simple sets. We then leverage the positive semi-definiteness of localizing matrices, in conjunction with the Łojasiewicz inequality, to extend these results to any compact subset of the original simple set. Our approach progresses from specific cases to more general ones, starting with algebraic varieties.

 \subsection{Error estimation on the intersection of a real algebraic variety
 with a product of simple sets}
 
 We first consider the case where $\CX$ is the intersection of a real algebraic variety and a product of simple sets, i.e., 
 \begin{eqnarray} \label{eq-real-variety}
     \CX = \{\bx \in \RR^n: h_i(\bx) = 0 \ \forall i \in [p]\} \,\cap\, 
     \CY.
 \end{eqnarray}
 Here, the equalities $h_i(\bx) =0 \; \forall i \in [p]$ defines a real algebraic variety in $\RR^n$ and 
 $\CY=\Pi_{i=1}^s\CY_i$ is a product of $s$ simple sets $\CY_i$ in $\RR^{n_i}$. For convenience, we express $\CY$ using $m$ polynomial inequalities as follows:  
 \begin{eqnarray}
     \CY = \{\bx \in \RR^n: g_j(\bx) \geq 0 \ \forall j \in [m]\}.
     \label{eq-Y}
 \end{eqnarray}
Since $\CX$ is contained in a product of simple sets $\CY$, without loss of generality, we can add the polynomial inequality $g_0(\bx):= R^2 - \|\bx\|^2 \geq 0 $ for a suitable positive number $R$ to the description of $\CY$ without changing the domain. That is, $\CY \subset \mathbb{B}_R$. 

We recall that $\lceil g_j \rceil = \lceil {\rm deg} (g_j) /2 \rceil \; \forall j \in [m]$, 
and 
$ \lceil g_J \rceil= \lceil {\rm deg} (g_J)/2\rceil$ for an index set 
$J\subset \{0,1,\ldots,m\}$.
Define   
\begin{eqnarray} \label{eq-d}
    d:= \max\{\lceil h_i \rceil,\; \lceil g_j \rceil :\; i \in [p],\; 0 \leq j \leq m\}.
\end{eqnarray}
For any integers $r$ and $k$ such that 
$2r \geq k  \geq 4d$, we define
 \begin{multline*}
     \CM(\CR(\CX)_{2r})=   \Bigl\{\iy \in \RR^{s(n,2r)}\;:\; \iy_0 =1,\ \bM_r(\iy) \succeq 0,\; \ell_{\iy}(h_i^2(\bx))=0\ \forall i \in [p],\\
     \bM_{r-\lceil g_J \rceil}(g_J\iy) \succeq 0 \; \forall J \subset \{0,1,\dots,m\}\; \mbox{such that}\;\degg{g_J} \leq r\Bigr\},
 \end{multline*}
 and $\CM_k(\CR(\CX)_{2r}) = \{ \pi_k(y)\,:\, y\in \CM(\CR(\CX)_{2r})\}$. 
 We observe that 
 \begin{displaymath}
     \CM_k(\CX) \subset \CM_k(\CT(\CX)_{2r})  \subset \CM_k(\CR(\CX)_{2r}).
 \end{displaymath}
 Hence, the Hausdorff distances satisfy that 
 \begin{equation}\label{inequality of distances }
     \bd_k(\CT(\CX)_{2r})  \leq \bd(\CM_k(\CR(\CX)_{2r}),\CM_k(\CX)) =: \bd_k(\CR(\CX)_{2r}).
 \end{equation}
 Moreover, the outer approximation $\CM_k(\CR(\CX)_{2r})$ tends to $\CM_k(\CX)$ as $r \to +\infty$ in the sense of Hausdorff distance, as stated in the following proposition.
 \begin{proposition}
     For any $k \in \NN$, we have $\lim_{r \to \infty}\bd_k(\CR(\CX)_{2r}) =0$.
 \end{proposition}
 \begin{proof}
     We observe that 
     \begin{displaymath}
         \cdots \subset \CM_k(\CR(\CX)_{2r}) \subset \CM_k(\CR(\CX)_{2r+2}) \subset \CM_k(\CR(\CX)_{2r+4}) \subset \cdots.
     \end{displaymath}
     Thus it is sufficient to prove that $\bigcap_{r \in \NN} \CM_k(\CR(\CX)_{2r})= \CM_k(\CX)$.
     Indeed, let $\iy$ be a sequence belonging to $\CM_k(\CR(\CX)_{2r})$ for all $r$. Then there exists an infinite sequence $\overline{\iy}$ such that the first $s(n,k)$ coordinates of $\overline{\iy}$ is $\iy$ and it satisfies the following conditions:
     \begin{equation*}
         \overline{\iy}_0=1,\quad \bM(\overline{\iy}) \succeq 0,\quad 
          \bM_{r-\lceil g_J \rceil}(g_J\overline{\iy}) \succeq 0 \quad \forall J \subset \{0,1,\dots,m\},\; \forall r \in \NN,\quad \ell_{\overline{\iy}}(h_i^2)=0 \quad  \forall i \in [p].
     \end{equation*}
     
     According to \cite[Theorem 3.8]{lasserre2009moments}, the conditions $\bM(\overline{\iy}) \succeq 0$ and $ \bM_{r-\degg{g_J}}(g_J\overline{\iy}) \succeq 0 \; \forall J \subset \{0,1,\dots,m\},\ r \in \NN$, imply that there exists a probability measure $\mu$ supported on $\CY$
     that is represented by
     the moment sequence $\overline{\iy}$. Thus, the condition $\ell_{\overline{\iy}}(h_i^2)=0\;\forall i \in [p]$ is equivalent to 
     \begin{equation*}
         \int_{\mathbb{B}_R}h_i(\bx)^2 d\mu(\bx) = 0 \quad \forall i \in [p].
     \end{equation*}
     This implies that $\int_{\mathbb{B}_R\backslash\CX}h_i(\bx)^2 d\mu(\bx) = 0\; \forall\; i\in[p]$, and hence
     $\mbox{supp}(\mu) \subset \CX$. Therefore, $\mu \in \CP(\CX)$ and $\overline{\iy} \in \CM(\CX)$, which directly leads to $\iy \in \CM_k(\CX)$.
 \end{proof} 
 
 To estimate the error of pseudo-moment sequences on $\CX$, we introduce the Łojasiewicz inequality in the next lemma, which plays a key role in our estimation. 
 
 \begin{lemma}[\cite{Real_algebraic_geometry},Corollary 2.6.7]\label{lemma: Lojasiewicz inequality}
     Let $B$ be a compact semi-algebraic set, and $f$ and $g$ be two continuous 
     semi-algebraic functions from $B$ to $\RR$ such that $f^{-1}(0) \subset g^{-1}(0)$. Then there exist a Łojasiewicz constant $c > 0 $ and a Łojasiewicz exponent $\Lo > 0$ such that 
     \begin{displaymath}
         |g(\bx)| \leq c|f(\bx)|^{\Lo} \quad \forall\; \bx \in B.
     \end{displaymath}
 \end{lemma}
 Define the distance function:
 \begin{displaymath}
     \bd_{\CX}(\bx) = \bd(\bx,\CX), \quad \bx \in \CY.
 \end{displaymath}
 Since $\CX$ is compact, the set of minimizers $\pi_{\CX}(\bx)$ of the
 problem $\min\{ \|\by-\bx\|\,:\, \by\in \CX\}$
  is non-empty and compact. Moreover, the fact that $\CX$ is a basic semi-algebraic set implies that $\bd_{\CX}(\bx)$ is a continuous semi-algebraic function. We next define the function 
 \begin{displaymath}
     f(\bx) = \max\{|h_i(\bx)|\ :\ i \in [p] \}, \quad \bx\in\CY,
 \end{displaymath}
 which is also a continuous semi-algebraic function. Moreover,  the following relation holds true 
 \begin{displaymath}
     \bd_{\CX}^{-1}(0)= \CX = f^{-1}(0).
 \end{displaymath}
 Hence, applying Lemma~\ref{lemma: Lojasiewicz inequality} gives us the following
 inequality with Łojasiewicz constant $c $ and  exponent $\Lo $:
 \begin{equation}\label{Lojaciewicz inequality}
     \bd_{\CX}(\bx) \leq c\max\{|h_i(\bx)|\ :\ i \in [p] \}^{\Lo} \quad \forall\; \bx \in \CY.
 \end{equation}
 Without loss of generality, 
 we may assume that $\Lo \leq 1$ since otherwise, we can replace $\Lo$ by $1$ and multiply the Łojasiewicz constant by $\max_{\bx \in \CY}\max\{|h_i(\bx)|\ :\ i \in [p] \}^{L-1} <\infty$ 
 to obtain a new inequality with the Łojasiewicz exponent $1$.
 We now apply the inequality \eqref{Lojaciewicz inequality} to estimate the error of truncated pseudo-moment sequences in $\CM_k(\CR(\CX)_{2r})$.
 
\begin{theorem}\label{thm: error on variety}
    Let $\CX$ be a semi-algebraic set defined in \eqref{eq-real-variety}, $\CY=\Pi_{i=1}^s\CY_i$ is a product of $s$ simple sets $\CY_i\subset \RR^{n_i}$ as expressed in \eqref{eq-Y}. We assume that the inequality $g_0(\bx) = R^2 -\|\bx\|^2 \geq 0$ is added to the description of $\CX$. For any $k = 2l \in\NN$ and $r \in \NN$ such that $k\geq 4d$ and $r \geq 2s(\max\{n_1,\dots,n_s\}+1)k+s$,  where $d$ is defined in \eqref{eq-d}, the Hausdorff distance $\bd_k(\CR(\CX)_{2r}) := \bd(\CM_k(\CR(\CX)_{2r}),\CM_k(\CX)) $ admits an upper bound as follows:
    \begin{displaymath}
        \bd_k(\CR(\CX)_{2r}) \leq 
        \dfrac{2\gamma(R,n,k)\Gamma(\CY,k)}{(r-s)^2}+ cL(R,k)\left(\dfrac{2\gamma(R,n,k)\,
        \Gamma(\CY,k)}{(r-s)^2}\right)^{\Lo/2}
        \left(\sum_{i=1}^p\|h_i^2\|_{1}\right)^{\Lo/2}.
    \end{displaymath}
    Here, $L(R,k)$ is the Lipschitz number of $\bv_k(\bx)$ on the compact set $\mathbb{B}_R$, $\Gamma(\CY,k)$ is a parameter defined in Theorem~\ref{thm: error on product of simple sets}. In short, the error of  $k$-truncated pseudo-moment sequences in the $r$-th level of Schmüdgen-type moment hierarchy is $\mathrm{O}(1/r^{\Lo})$.
\end{theorem}
\begin{proof}
\begin{figure}[!ht]
\centering
\begin{tikzpicture}

    \begin{scope}[blend mode=multiply]
    \filldraw[fill=gray!10, draw=black] (1,0) ellipse (2.9cm and 1.5cm);
    \filldraw[fill=white!10, draw=black] (0,0) circle (1.2cm);
    \filldraw[fill=gray!10, draw=black] (0,0) ellipse (1.9cm and 2.8cm);
    \filldraw[fill=white!10, draw=black] (1,0) circle (3.5cm);
    \end{scope}
    
    \node at (0,0) {$\CM_k(\CX)$};
    \node at (2.7,0) {$\CM_k(\CY)$};
    \node at (0,1.8) {$\CM_k(\CR(\CX)_{2r})$};
    \node at (2.5,-2.1) {$\CM_k(\CT(\CY)_{2r})$};
    \node at (0,-2.5) (S1) {$\iy$};
    \node at (2.7,-0.8) (S2) {$\overline{\iy}$};
    \node at (0,-0.7) (S3) {$\widetilde{\iy}$};

    \draw[arrow] (S1.east) -- (S2.west);
    \draw[arrow] (S2.west) -- (S3.south);
    \draw[arrow] (S1.east) -- (S3.south);

\end{tikzpicture}
\caption{Projection of a sequence $y\in \CM_k(\CT(\CY)_{2r})$ onto $\CM_k(\CX).$}
\label{Fig: projection with real algebraic variety}
\end{figure}
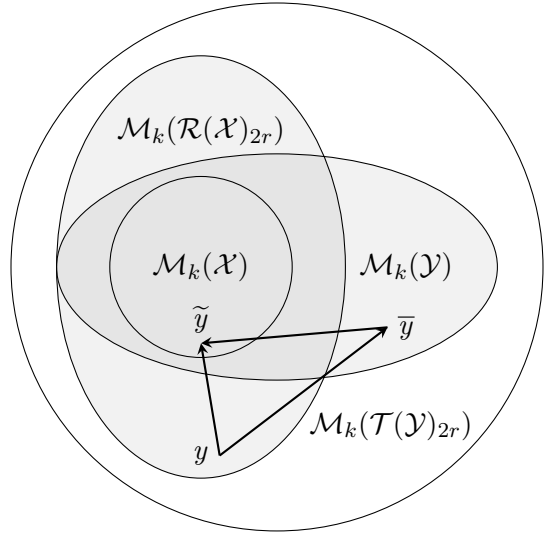
We prove the theorem by evaluating the distance between an arbitrary pseudo-moment sequence $\iy \in \CM_k(\CR(\CX)_{2r})$ and $\CM_k(\CX)$. First, we conduct two consecutive projections as follows: since $\CM_k(\CX) \subset \CM_k(\CY)$,
we first project $\iy$ onto $\CM_k(\CY)$ and denote its projection by $\overline{\iy}$. By Tchakaloff's theorem, $\overline{\iy}$ can be written as a convex combination of some $\bv_k(\overline{\bx}_j)'$s with points $\{\overline{\bx}_j : j\in[N]\} \subset \CY$. We continue projecting $\overline{\bx}_j$ onto $\CX$ to define a new moment sequence $\widetilde{\iy} \in \CM_k(\CX)$. These projections are demonstrated in Figure~\ref{Fig: projection with real algebraic variety} and elaborated in the later part of the proof. Then we can upper bound the distance of $\iy$ to $\CM_k(\CX)$ by 
the triangle inequality: 
\begin{equation}\label{eq: inequality for distance with variety}
    \bd(\iy,\CM_k(\CX)) \;\leq\; \|\iy - \widetilde{\iy}\| \;\leq\; \|\iy-\overline{\iy}\| + \|\overline{\iy}-\widetilde{\iy}\|.
\end{equation}
The idea of upper-bounding $\bd(\iy,\CM_k(\CX))$ by $\|\iy-\overline{\iy}\|$ and $\|\overline{\iy}-\widetilde{\iy}\|$ is based on the fact that we already have Theorem~\ref{thm: error on product of simple sets} and the Łojasiewicz inequality~\eqref{Lojaciewicz inequality} as tools to evaluate these terms, which is elaborated next. 

\textbf{Evaluating $\|\iy-\overline{\iy}\|$:}  according to Theorem~\ref{thm: error on product of simple sets} under the 
condition $r \geq 2s(\max\{n_1,\dots,n_s\}+1)k+s$, we have
    \begin{equation}\label{eq: inequality for distance with variety 1}
        \|\iy -\overline{\iy}\| \;\leq\; 
        \dfrac{2\gamma(R,n,k)\,\Gamma(\CY,k)}{(r-s)^2}\;=: \varepsilon
    \end{equation}
    
    \textbf{Evaluating $\|\overline{\iy}-\widetilde{\iy}\|$:} according to Tchakaloff's theorem, there exists at most $N = s(n,k)$ points $\{\overline{\bx}_j:j\in[N]\}$ in $\CY$ and positive real numbers $\{w_j : j\in[N]\}$ such that $\sum_{j=1}^N w_j =1$ and
    \begin{displaymath}
        \overline{\iy}= \sum_{j=1}^N w_j\bv_k(\overline{\bx}_j).
    \end{displaymath}
    Therefore, we can rewrite the equalities of the Riesz functional associated with $\iy$ and $\overline{\iy}$ as follows:
    \begin{equation*}
        \ell_{\iy}(h_i^2) = \langle \mathbf{h}_i^2, \iy \rangle= 0, 
        \quad \mbox{and} \quad  \ell_{\overline{\iy}}(h_i^2) =\langle \mathbf{h}_i^2, \overline{\iy} \rangle= \sum_{j=1}^N w_j \langle \mathbf{h}_i^2,\bv_k(\overline{\bx}_j)\rangle =
        \sum_{j=1}^N w_j h_i^2(\overline{\bx}_j),
    \end{equation*}
    where $\mathbf{h}_i^2$ denotes the vector of coefficients of $h_i^2(\bx)$. Since $\|\iy -\overline{\iy}\| \leq \varepsilon$, by the Cauchy–Schwarz 
    inequality, we have the following inequality for any $i \in [p]$:
    \begin{eqnarray*}
        &\varepsilon \|h_i^2\|_1 \geq \left| \langle \mathbf{h}^2_i, \iy - \overline{\iy} \rangle \right| = \Big| \ell_{\iy}(h_i^2) - \sum_{j=1}^N w_j {h}_i^2(\overline{\bx}_j)  \Big| 
        = \sum_{j=1}^N w_j h_i(\overline{\bx}_j)^2.
    \end{eqnarray*}
     Thus, applying the Cauchy–Schwarz inequality again, we obtain that
    \begin{eqnarray*}
     \varepsilon \sum_{i=1}^p \|h_i^2\|_{1} 
     &\geq& \sum_{j=1}^N w_j \sum_{i=1}^p h_i(\overline{\bx}_j)^2 \geq \sum_{j=1}^N w_j 
        \max\left\{ |h_i(\overline{\bx}_j)|:\ i \in [p]\right\}^2
        \\
        &=& \left(\sum_{j=1}^N w_j\right)\left(\sum_{j=1}^N w_j \max\left\{ h_i(\overline{\bx}_j):\ i \in [p]\right\}^2\right)
        \\
        & \geq& \left(\sum_{j=1}^N w_j \max\left\{ |h_i(\overline{\bx}_j)|:\ i \in [p]\right\} \right)^2
        \\
        \Rightarrow \ \varepsilon^{1/2} \Big(\sum_{i=1}^p \|h_i^2\|_{1}\Big)^{1/2}
        &\geq& \sum_{j=1}^N w_j \max\left\{ |h_i(\overline{\bx}_j)|:\ i \in [p]\right\} .
    \end{eqnarray*}

        Let ${\bx}_j\in\CX$ be the projection of $\overline{\bx}_j\in\CY$ onto $\CX$ for all $j \in [N]$. We define 
    \begin{equation*}
        \widetilde{\iy} \;=\; \sum_{j=1}^N w_j\bv_k({\bx}_j) \;\in \CM_k(\CX).
    \end{equation*}
    In addition, based on the Łojasiewicz inequality~\eqref{Lojaciewicz inequality}, we have
    \begin{displaymath}
        \|\bx_j -\overline{\bx}_j\| = \bd_{\CX}(\overline{\bx}_j) \leq c \max\left\{ |h_i(\overline{\bx}_j)|:\ i \in [p]\right\}^{\Lo}\;\;\forall\; j\in[N].
    \end{displaymath}
    We then use the common Lipschitz number $L(R,k)$ to get the inequality below:
    \begin{displaymath}
        \|\bv_k(\bx_j) - \bv_k(\overline{\bx}_j)\| \leq L(R,k) \|\bx_j -\overline{\bx}_j\|  \leq cL(R,k) \max\left\{ |h_i(\overline{\bx}_j)|:\ i \in [p]\right\}^{\Lo}\;\;\forall\; j\in[N].
    \end{displaymath}
    Hence, we obtain that 
    \begin{align*}
        &\|\overline{\iy}-\widetilde{\iy}\|\;=\;\left\|\sum_{j=1}^N 
        w_j\bv_k(\overline{\bx}_j) -\sum_{j=1}^N w_j\bv_k(\bx_j)\right\|\\ 
        \leq\;& \sum_{j=1}^N w_j\|\bv_k(\overline{\bx}_j) - \bv_k(\bx_j)\| 
        \;\leq \; cL(R,k)\sum_{j=1}^N w_j\max\left\{ |h_i(\overline{\bx}_j)|:\ i \in [p]\right\}^{\Lo}.
    \end{align*} 
    Since $\Lo \leq 1$, applying the Jensen's inequality gives us the following inequality:
        \begin{displaymath}
            \sum_{j=1}^N w_j\max\left\{ |h_i(\overline{\bx}_j)|:\ i \in [p]\right\}^{\Lo} \leq \left(\sum_{j=1}^N w_j\max\left\{ |h_i(\overline{\bx}_j)|:\ i \in [p]\right\}\right)^{\Lo}\leq  \varepsilon^{\Lo/2} \left(\sum_{i=1}^p \|h_i^2\|_1 \right) ^{\Lo/2}.
        \end{displaymath}
        Combining the above with the inequality \eqref{eq: inequality for distance with variety 1}, we obtain that 
        \begin{align*}
            & \bd_k(\CR(\CX)_{2r}) \;\leq\;\|y - \tilde{y}\| \leq \|y-\overline{y}\| + \|\overline{y}-\tilde{y}\|
            \;\leq \; \varepsilon + cL(R,k)\varepsilon^{\Lo/2} \left(\sum_{i=1}^p \|h_i^2\|_{1} \right) ^{\Lo/2}\\ 
            \leq\;&  \dfrac{2\gamma(R,n,k)\Gamma(\CY,k)}{r^2}+ cL(R,k)\left(\dfrac{2\gamma(R,n,k)\Gamma(\CY,k)}{(r-s)^2}\right)^{\Lo/2}\left(\sum_{i=1}^p\|h_i^2\|_{1}\right)^{\Lo/2}.
        \end{align*}
        This completes the proof. 
\end{proof}

\subsection{Error estimation on a compact semi-algebraic set}

We now extend our analysis to any compact semi-algebraic set $\CX$ defined by \eqref{domain X}, i.e.,
\begin{eqnarray*}
    \CX = \big\{\bx \in \RR^n: \  g_j(\bx) \geq 0 \ \forall j \in [m],\;
    h_i(\bx) = 0 \ \forall i \in [p]\big\}.
    \label{eq-CX}
\end{eqnarray*}
Without loss of generality, we can assume that the Archimedean condition is satisfied by adding the following inequality for suitable $R \geq 1$ into the description of $\CX$ without changing the domain:
\begin{displaymath}
    g_0(\bx):= R^2 - \|\bx\|^2 \geq 0.
\end{displaymath}

We note that the method used in Section~\ref{section: error in algebraic set} is actually valid for any intersection of a real algebraic variety and a product of simple sets. In order to reuse this method for a general semi-algebraic set, we perform a lifting of $\CX$ into $\RR^{n+m}$ by the following polynomial mapping:
\begin{displaymath}
    \varphi: \RR^n \to \RR^{n+m}: \quad \bx \mapsto \varphi(\bx):= (\bx,g_1(\bx),\dots,g_m(\bx)).
\end{displaymath}

Recall that $d = \max\{ \lceil  h_i \rceil, \  \lceil  g_j  \rceil \,:\ i \in [p],\ j \in [m]\}$. Then by a simple estimation, we obtain that 
\begin{eqnarray} \label{eq-K}
    \sum_{j=1}^mg_j(\bx) \leq R^{2d}\sum_{j=1}^m \|g_j\|_{1} =: K, \quad
    \forall\; \bx\in \CX \subset \mathbb{B}_R.
\end{eqnarray}
Hence, we have that 
\begin{displaymath}
    \varphi(\bx) \in \mathbb{B}_R \times \Delta^m_K \quad \forall\; \bx  \in \CX,
\end{displaymath}
where $\Delta^m_K := \left\{ \bu \in \RR^m:\ u_j \geq 0 \ \forall j \in [m],\ \sum_{j=1}^m u_j \leq K \right\}$ is an $m-$dimensional simplex. As an attempt to reuse Theorem~\ref{thm: error on variety}, we observe that the image of $\CX$ via the lifting $\varphi$ is the intersection of a real variety and a product of simple sets $\mathbb{B}_R \times \Delta^m_K$, i.e., 
\begin{eqnarray}\label{def of varphi(X)}
     &&\mathbb{B}_R \times \Delta_K^m \;=\;\Bigl\{ \bz = (\bx,\bu) \in \RR^{n}\times \RR^m \; :\; p_0(\bz) := R^2 - \|\bx\|^2 \geq 0,\; p_j(\bz):= u_j \geq 0 \ \forall j \in [m],\nonumber\\ 
     && \hspace{26mm} p_{m+1}(\bz):= K - \sum_{j=1}^m u_j\geq 0\Bigr\},\nonumber\\
    &&\varphi(\CX)= \Bigl\{\bz =  (\bx,\bu) \in \RR^{n}\times \RR^m\;:\; h_i(\bx) =0 \; \forall i \in [p],\;q_j(\bz):= u_j - g_j(\bx) = 0 \; \forall j \in [m],\nonumber\\
    &&\hspace{17mm}
    p_0(\bz) \geq 0,\;
    p_j(\bz)=u_j \geq 0 \; \forall j \in [m],\; p_{m+1}(\bz)=K - \sum_{j=1}^m u_j \geq 0
    \Bigr\}.
\end{eqnarray}
 We show in the following lemma that $\CM(\CT(\mathbb{B}_R \times \Delta_K^m)_{2r})$ is a compact set, and we can bound the corresponding Hausdorff distance 
 by Theorem~\ref{thm: error on product of simple sets}.
 \begin{lemma}\label{lem: bounded of M_k(B_R x Delta)}
     For any positive integer $k$, $\CM_k(\CT(\varphi(\CX))_{2r})$, $\CM_k(\CR(\varphi(\CX))_{2r})$, and $\CM_k(\CT(\mathbb{B}_R \times \Delta_K^m)_{2r})$ are compact sets contained in the Euclidean ball centered at the origin with radius $\gamma(n+m,R+mK^2,k)$. In addition, we have 
     \begin{displaymath}
         \bd_k(\CT(\mathbb{B}_R \times \Delta_K^m)_{2r}) \leq \Gamma(\mathbb{B}_R \times \Delta_K^m,k)\dfrac{2\gamma(n+m,R+mK^2,k)}{r^2}.
     \end{displaymath}
 \end{lemma}
 \begin{proof}
     We show that the constraint $\widehat{p}_0(\bz) :=R+ mK^2 - \|\bx\|^2 - \sum_{j=1}^m u_j^2 \geq 0$ can be added to the description of $\mathbb{B}_R \times \Delta_K^m$ without changing $\CM_k(\CT(\mathbb{B}_R \times \Delta_K^m)_{2r})$.
     Indeed, we can write
     \begin{multline*}
         K^2 - u_j^2= (K-u_j)(K+u_j)= \Bigl(K-\sum_{l=1}^m u_l +\sum_{l \neq j}u_l\Bigr)(K+u_j)
         \\
         = K\Bigl(K-\sum_{l=1}^m u_l\Bigr)+u_j\Bigl(K-\sum_{l=1}^m u_l\Bigr) 
         + K \sum_{l\neq j} u_l + \sum_{l\neq j} u_l u_j \in 
         \CT\big((\mathbb{B}_R \times \Delta_K^m)_{2r}\big)
         \quad \forall j \in [m].
     \end{multline*}
     This implies that we can write $\widehat{p}_0$ as a sum of elements of $\CT(\mathbb{B}_R \times \Delta_K^m)_{2r}$ of degree at most $2$. In other words, for any product polynomial $q$ of $p_0,\; p_j,\; j \in [m+1]$ with degree at most $2t$ and $\iy \in \CM(\CT(\mathbb{B}_R \times \Delta_K^m)_{2r})$, we have 
     \begin{eqnarray*}
         &&\widehat{p}_0q = p_0q+\sum_{j=1}^m\Biggl[K\Bigl(K-\sum_{l=1}^m  
         u_l\Bigr)+u_j\Bigl(K-\sum_{l=1}^m u_l\Bigr) 
         + K \sum_{l \neq j}u_l
         + \sum_{l\neq j}u_lu_j
         \Bigr]q
         \\
         &\Rightarrow&\; \bM_{r-t-1}(\widehat{p}_0q\iy)
         = \bM_{r-t-1}(p_0q\iy)+\sum_{j=1}^mK\bM_{r-t-1}\Biggl(\Bigl(K-\sum_{l=1}^m
         u_l\Bigr)q\iy\Biggr)\\
         &&\qquad +\bM_{r-t-1}\Biggl(u_j\Bigl(K-\sum_{l=1}^m u_l\Bigr)q\iy\Biggr) +\sum_{l \neq j} K \bM_{r-t-1} (u_l qy) 
         + \sum_{l\neq j} \bM_{r-t-1}\Biggl(u_l u_jq\iy\Biggr) \succeq 0.
     \end{eqnarray*}
     If we add $\widehat{p}_0(\bz) \geq 0$ into the description of $\mathbb{B}_R \times \Delta_K^m$, the new constraint added to $\CM(\CT(\mathbb{B}_R \times \Delta_K^m)_{2r})$
     then take the form of 
     \begin{equation*}
         \bM_{r-t-1}(\widehat{p}_0q\iy) \succeq 0.
     \end{equation*}
     However, as we have shown above, the constraint is already satisfied without adding $\widehat{p}_0(\bz) \geq 0$. Thus, adding it does not change 
     $\CM(\CT(\mathbb{B}_R \times \Delta_K^m)_{2r})$. Furthermore, applying Lemma~\ref{lemm: radius of preordering} gives us that $\CM_k(\CT(\mathbb{B}_R \times \Delta_K^m)_{2r})$ is a compact set contained in the Euclidean ball centered at the origin with radius $\gamma(n+m,R+mK^2,k)$. Applying Theorem~\ref{thm: error on product of simple sets} for the product of two simple sets, we obtain that
     \begin{equation*}
         \bd_k(\CT(\mathbb{B}_R \times \Delta_K^m)_{2r}) \leq \Gamma(\mathbb{B}_R \times \Delta_K^m,k)\dfrac{2\gamma(n+m,R+mK^2,k)}{(r-2)^2},
     \end{equation*}
     which is the upper bound we desire. 
 
  Notice that $\varphi(\CX) \subset \mathbb{B}_R \times \Delta_K^m$, and both $\CM_k(\CT(\varphi(\CX))_{2r})$ and $\CM_k(\CR(\varphi(\CX))_{2r})$ are closed subsets of $\CM_k(\CT(\mathbb{B}_R \times \Delta_K^m)_{2r})$. Hence, we can claim that these sets are also compact.
  \end{proof}
 Since $\varphi(\CX)$ is the intersection of the real algebraic variety 
 \begin{displaymath}
     \left\{\bz = (\bx,\bu) \in \RR^{n}\times \RR^m\; :\; h_i(\bx) =0 \; \forall i \in[p],\; u_j - g_j(\bx) =0 \; \forall j \in [m] \right\}
 \end{displaymath}
 and the product $\mathbb{B}_R \times \Delta^m_K$ of the simple sets $\mathbb{B}_R$ and $\Delta^m_K$, we can therefore  apply Theorem~\ref{thm: error on variety} to determine the tightness of the relaxation $\CM_k(\CT(\varphi(\CX))_{2r})$ of $\CM_k(\varphi(\CX))$. To convey that tightness back to the relaxation $\CM_k(\CT(\CX)_{2r})$ of $\CM_k(\CX)$, we need to examine the connection between the truncated pseudo-moment sequences on $\CX$ and those on $\varphi(\CX)$.
Indeed, the examination is conducted as follows: For any positive integers
$r$ and $k$, we set $t = \lfloor r/(2d) \rfloor$ and always assume in this section that $2t \geq k$. Recall the special superset $\CM_k(\CR(\CX)_{2r})$ 
of $\CM_k(\CX)$ defined in Section~\ref{sec-2.2}:
\begin{align*}
    \CM_k(\CR(\CX)_{2r}) =  \Bigl\{\pi_k(\iy): &\ \iy \in \RR^{s(n,2r)},\ \iy_0=1,\ \bM_r(\iy) \succeq 0,\ \ell_{\iy}(h_i^2 )=0 \ \forall i \in [p],\\
    \quad & \bM_{r-\lceil g_J \rceil}(g_J \iy) \succeq 0 \ \forall J \subset [m],\ \lceil g_J \rceil \leq r  \Bigr\}.
\end{align*}
Similar inclusions  as in Section~\ref{section: error in algebraic set} also hold, i.e., we obtain that  
 \begin{displaymath}
     \CM_k(\CX) \subset \CM_k(\CT(\CX)_{2r})  \subset \CM_k(\CR(\CX)_{2r}).
 \end{displaymath}
 Hence, the Hausdorff distances between these sets satisfy the  following inequality:  
 \begin{equation}
     \bd_k(\CT(\CX)_{2r})  \leq \bd_k(\CM_k(\CR(\CX)_{2r}),\CM_k(\CX)) =: \bd_k(\CR(\CX)_{2r}).
 \end{equation}

We next show that the lifting from $\bx \mapsto \varphi(\bx)$ induces a lifting of moment sequences from $\RR^{s(n,2r)}$ to $\RR^{s(n+m,2t)}$ defined as follows: 
 \begin{equation}\label{eq-ylift}
     \iy^{\varphi}_{(\alpha,\beta)} =\ell_{\iy^{\varphi}}(\bz^{(\alpha,\beta)}) := \ell_{\iy}(\bx^{\alpha}g(\bx)^{\beta}),\; \forall\; (\alpha,\beta) \in \NN^n\times\NN^{m},\; |\alpha| + |\beta| \leq 2t,    
 \end{equation}
 where $g(\bx)=(g_1(\bx),\ldots,g_m(\bx))$.
Here, $\bz^{(\alpha,\beta)}=(\bx,\bu)^{(\alpha,\beta)}=\bx^\alpha\bu^\beta.$ Equation \eqref{eq-ylift} is well-defined since we have ${\rm deg} (\bx^{\alpha}g(\bx)^{\beta}) \leq |\alpha| + 2d|\beta| \leq 2d(|\alpha| +|\beta|) \leq 4dt \leq 2r.$  
In particular, $\iy^{\varphi}_{(\alpha,0_m)}=\ell_y(\bx^\alpha) = \iy_\alpha$ for all $|\alpha|\leq 2t,$ and $\ell_{\iy^{\varphi}}(p(\bx)\bu^\beta)=\ell_{y}(p(\bx)g(\bx)^\beta)$ for any
$p\in \RR[\bx]$ such that ${\rm deg}(p)+|\beta|\leq 2t.$
 
 The following lemma shows that if $\iy$ is a truncated pseudo-moment sequence on $\CX$, then
 $\iy^\varphi$ 
 is a truncated pseudo-moment sequence on 
 $\varphi(\CX)$. We 
 adopt the notational convention that if $\iy \in \RR^{s(n,2r)}$, then  $\pi_k(y)$ denotes the projection onto the first $s(n,k)$  coordinates of $\iy$.
 
 \begin{lemma}\label{lemma: lift of moment sequence}
 We set $t = \lfloor r/(2d) \rfloor$. 
 Let $\iy \in \CM(\CR(\CX)_{2r})$ be a truncated pseudo-moment sequence, and $\iy^{\varphi}$ be the sequence defined as in
\eqref{eq-ylift}. If $t \geq 2d$, then $\iy^{\varphi} \in \CM(\CT(\varphi(\CX))_{2t})$.
 \end{lemma}
 \begin{proof}
 We recall the description of $\varphi(\CX)$:
 \begin{eqnarray*}
   && \varphi(\CX)= \Bigl\{\bz =  (\bx,\bu) \in \RR^{n}\times \RR^m\;:\; h_i(\bx) =0 \; \forall i \in [p],\; 
    q_j(\bz):=u_j - g_j(\bx) = 0 \; \forall j \in [m],
    \\ 
    &&\qquad\quad p_0(\bz):=R^2-\|\bx\|^2 \geq 0,\; p_j(\bz):=u_j \geq 0 \; \forall j \in [m],\; p_{m+1}(\bz):=K - \sum_{j=1}^mu_j \geq 0
    \Bigr\}.
\end{eqnarray*}
Then 
\begin{multline*}
    \CM(\CT(\varphi(\CX))_{2t})= \Bigl\{\iy^\varphi \in \RR^{s(n+m,2t)}:\; \ell_{\iy^\varphi}(h_i^2) =0\; \forall i \in [p],\; \ell_{\iy^\varphi}(q_j(\bz)^2)=0 \; \forall j \in [m],
    \\
    \iy^\varphi_0 =1,\; \bM_t(\iy^\varphi) \succeq 0,\; \bM_{t-\lceil p_J \rceil}(p_J\iy^\varphi) \succeq 0\; \forall\, J \subset \{0,1,\dots,m+1\},\; \lceil p_J \rceil \leq t
    \Bigr\}.
\end{multline*}
 Here, for any $J \subset \{0,1,\dots,m+1\}$, define $p_J(\bz) = \prod_{j \in J}p_j(\bz)$. The condition $t \geq 2d$ ensures that $\ell_{\iy^\varphi}(h_i^2)$ and $\ell_{\iy^\varphi}((u_j-g_j(\bx))^2)$ are well-defined for any $i \in [p],\ j \in [m]$.
 
     To prove that $\iy^{\varphi} \in \CM(\CT(\varphi(\CX))_{2t})$, we need to prove the followings:
     \begin{enumerate}
        \item $\iy^{\varphi}_0 =1,\; \ell_{\iy^{\varphi}}(h_i^2) =0\; \forall\; i \in [p],\; \ell_{\iy^{\varphi}}((u_j-g_j(\bx))^2)=0 \; \forall\; j \in [m]$,
         \item $\bM_t(\iy^{\varphi}) \succeq 0$,
         \item  $\bM_{t- \lceil p_J \rceil}(p_J \iy^{\varphi}) \succeq 0$
         for any $J\subset \{0,1,\ldots,m+1\}$ with $\degg{p_J}\leq t.$
     \end{enumerate}
     The first condition is straightforward from the definition of $\iy^{\varphi}$. Indeed, we have
     \begin{eqnarray*}
         &&\iy^{\varphi}_0 = \iy_0 =1, \\
         &&\ell_{\iy^{\varphi}}(h_i^2)=\ell_{\iy}(h_i^2) =0\;\; \forall\; i \in [p], \\ 
         &&\ell_{\iy^{\varphi}}((u_j-g_j(\bx))^2)=\ell_{\iy}((g_j(\bx)-g_j(\bx))^2)=0 \;\; \forall\; j \in [m].
     \end{eqnarray*}
     For the second condition, based on the definition of $\iy^{\varphi}$, the Riesz functional $\ell_{\iy^{\varphi}}$ satisfies that for any $p \in \RR[\bz]_t$, we have 
     \begin{displaymath}
         \ell_{\iy^{\varphi}}(p(\bz)) = \ell_{\iy}(p(\bx,g(\bx)))\quad \Rightarrow \quad \bM_t(\iy^{\varphi}) 
         = \ell_{\iy^{\varphi}}(\bv_t(\bz)
         \bv_t(\bz)^{\top})=\ell_{\iy}\left(\bv_t(\bx,g(\bx))\bv_t(\bx,g(\bx)) ^{\top} \right).
     \end{displaymath}
     Since $t = \lfloor r/(2d) \rfloor$, there exists an $s(n+m,t) \times s(n,r)$ matrix $T_1$ such that 
     \begin{displaymath}
         \bv_t(\bx,g(\bx)) = T_1 \bv_r(\bx) \quad \Rightarrow \quad \bM_t(\iy^{\varphi})=\ell_{\iy}\left(T_1\bv_r(\bx)\bv_r(\bx)^{\top}T_1^{\top}\right)=T_1 \bM_r(\iy) T_1^{\top} \succeq 0.
     \end{displaymath}
     
     The third condition is more complicated to verify.  For any positive integer $c$ and $\widehat{g} \in \RR[\bx]_{2c}$, we define the auxiliary functions $g_{m+1}= R^{2c}\|\widehat{g}\|_1-\widehat{g}$, and $g_0(\bx) = R^2 - \|\bx\|^2.$
     Let $J \subset \{0,1,\ldots,m,m+1\}$. 
     We first prove the following claim: If $\lceil g_J \rceil \leq r$, then $\bM_{r-\lceil g_J \rceil}(g_J\iy) \succeq 0$. 

     \textbf{Proof of the claim.} If $J \subset \{ 0,1\dots,m\}$, $\bM_{r-\lceil g_J \rceil}(g_J\iy)$ is one of the localizing matrix in the description of $\CM(\CR(\CX)_{2r})$. Hence, it is positive semidefinite. Otherwise, $J = J^\prime \cup \{m+1\} $ and $J^\prime \subset \{0,1\ldots,m\}$. We recall the inner product associated with $\iy$ in $\RR[\bx]_r$ defined by 
     \begin{displaymath}
         \langle p ,q \rangle_{\iy} = \mathbf{p}^{\top}\bM_r(\iy)\mathbf{q} \quad \forall\;
         p(\bx), q(\bx) \in \RR[\bx]_r.
     \end{displaymath}
    The inner product $\langle \cdot,\cdot \rangle_y$ possesses the following properties stated in the book\cite{lasserre2009moments}:
    \begin{align*}
        &\langle q_1,q_2q_3 \rangle_{\iy} = \langle q_1q_2,q_3 \rangle_{\iy} \quad \forall q_1, q_2, q_3 \in \RR[\bx],\; {\rm deg} (q_1) + {\rm deg} (q_2) \leq r,\; {\rm deg} (q_2) + {\rm deg} (q_3) \leq r;\\
        & \langle p, g_Jq\rangle_{\iy} = \mathbf{p}^{\top} \bM_{r-\lceil g_J \rceil}(g_J y)\mathbf{q} \quad \forall\; p, q \in 
        \RR[\bx]_{r-\degg{g_J}}.
    \end{align*}
    Therefore, proving the positive semidefiniteness of 
    $\bM_{r-\lceil g_J \rceil}(g_J\iy)$ is equivalent to proving that 
    \begin{equation}\label{eq sdp}
        \langle p, g_Jq\rangle_{\iy} = \langle p, {g_{m+1}} g_{J^\prime}q\rangle_{\iy} \;\geq\; 0 \quad \forall\; p, q \in \RR[\bx]_{r-\lceil g_J \rceil}.
    \end{equation}
    In what follows, we prove \eqref{eq sdp} by considering the form of $g_{m+1}$. 

    \textbf{Case 1: $\widehat{g}(\bx)= \bx^{2\alpha}$.}  For any $i \in [n]$ and $ q \in \RR[\bx]_{r-\lceil g_{J^\prime} \rceil-1}$, we have
\begin{align}\label{eq-local}
    &\langle q, (R^2-x_i^2)g_{J^\prime}q\rangle_{\iy} = \Biggl\langle q, \Bigl(R^2-\sum_{j=1}^nx_j^2\Bigr)g_{J^\prime}q\Biggr \rangle_{\iy} +\sum_{j \neq i}\langle q, x_j^2g_{J^\prime}q\rangle_{\iy}\nonumber\\
    =\;& \mathbf{q}^{\top}\bM_{r-\lceil g_{J^\prime} \rceil-1}\Biggl( \Bigl(R^2-\sum_{j=1}^nx_j^2\Bigr)g_{J^\prime} \iy\Biggr)\mathbf{q} + \sum_{j \neq i}\langle x_jq,g_{J^\prime}x_jq\rangle_{\iy}\nonumber\\
    =\;& \mathbf{q}^{\top}\bM_{r-\lceil g_{J^\prime} \rceil-1}\Biggl( \Bigl(R^2-\sum_{j=1}^nx_j^2\Bigr)g_{J^\prime} \iy\Biggr)\mathbf{q} + \sum_{j \neq i}\widehat{\mathbf{q}}_j^{\top}\bM_{r-\lceil g_{J^\prime} \rceil}( g_{J^\prime} \iy)\widehat{\mathbf{q}}_j \geq 0 \nonumber \\
    \Rightarrow\;& \langle q, x_i^2g_{J^\prime}q\rangle_{\iy} \leq R^2\langle q, g_{J^\prime}q\rangle_{\iy}.
\end{align}
In the third equality above, $\widehat{\mathbf{q}}_j$ denotes the vector of coefficients of $x_jq(\bx)$. The last inequality is based on the positive semidefiniteness of the matrices $\bM_{r-\lceil g_{J^\prime} \rceil-1}\big( \big(R^2-\sum_{j=1}^nx_j^2\big)g_{J^\prime} \iy\big)$ and $\bM_{r-\lceil g_{J^\prime} \rceil}( g_{J^\prime} \iy)$. 
Now, for any index $i$ such that $\alpha_i \neq 0$, we can apply \eqref{eq-local} to obtain 
\begin{eqnarray*}
    &\left \langle q, \bx^{2\alpha}g_{J^\prime}q \right \rangle_{\iy} = \Bigl \langle qx_i^{\alpha_i-1}\prod_{j \neq i}^nx_j^{\alpha_j}, x_i^2g_{J^\prime}qx_i^{\alpha_i-1}\prod_{j \neq i}^nx_j^{\alpha_j} \Bigr \rangle_{\iy} \leq R^2\Bigl \langle q,x_i^{2\alpha_i-2}\prod_{j\neq i}^nx_j^{2\alpha_j}g_{J^\prime}q \Bigr \rangle_{\iy}. 
\end{eqnarray*}
By  repeating the above process for the non-zero components of $\alpha$,
we get 
$$\left \langle q, \bx^{2\alpha}g_{J^\prime}q \right \rangle_{\iy}
\leq R^{2|\alpha|}\left \langle q, g_{J^\prime}q \right \rangle_{\iy} \ \forall\; q \in \RR[\bx]_{r-\lceil g_{J^\prime} \rceil-|\alpha|}.
$$

\textbf{Case 2: $\widehat{g}(\bx)= \pm\bx^{\alpha}$.} We construct $\alpha^{(1)}$ and $\alpha^{(2)}$ as follows: without loss of generality, we assume that there is an index $j$ such that $\alpha_i$ is even for all $i >j$ and $\alpha_i$ is odd for all $i\leq j$. 
Then we set $ \alpha^{(1)}_i=\alpha^{(2)}_i =\alpha_ i $ for $i>j$. 
For index $i \leq j$, we set 
$\alpha^{(1)}_i=\alpha_i+(-1)^i$ and $\alpha^{(2)}_i=\alpha_i-(-1)^i$. Thus, we have  
$\alpha^{(1)} = 2\beta^{(1)}$, $\alpha^{(2)} = 2\beta^{(2)}$
with $\beta^{(1)},\beta^{(2)} \in \NN^n_r$ and $\alpha = \beta^{(1)}+\beta^{(2)}$.
By using the identity 
$(\bx^{\beta^{(1)}} \pm \bx^{\beta^{(2)}})^2 = 
\bx^{2\beta^{(1)}} + \bx^{2\beta^{(2)}} \pm  2\bx^{\alpha} $,
we have that
\begin{eqnarray*}
\Bigl\langle q, (\bx^{2\beta^{(1)}} + \bx^{2\beta^{(2)}} \pm  2\bx^{\alpha})g_{J^\prime}  q\Bigr\rangle \;=\; 
\Bigl\langle (\bx^{\beta^{(1)}} \pm \bx^{\beta^{(2)}}) q, g_{J^\prime} (\bx^{\beta^{(1)}} \pm \bx^{\beta^{(2)}}) q\Bigr\rangle_y\; \geq\; 0.
\end{eqnarray*}
Recall that $R \geq 1$. Then we have 
\begin{align*}
   &2 \left \langle q, \pm \bx^{\alpha}g_{J^\prime} q \right \rangle_{\iy} \leq \left \langle q, \bx^{2\beta^{(1)}}g_{J^\prime} q \right \rangle_{\iy}+ \left \langle q, \bx^{2\beta^{(2)}}g_{J^\prime} q \right \rangle_{\iy} \leq 2R^{2\lceil |\alpha|/2\rceil} \left \langle q, g_{J^\prime} q \right \rangle_{\iy}\\
   \Rightarrow\; & \left \langle q, \pm \bx^{\alpha}g_{J^\prime} q \right \rangle_{\iy} \leq R^{2\lceil |\alpha|/2\rceil} \left \langle q, g_{J^\prime} q \right \rangle_{\iy}.
\end{align*}

\textbf{Case 3: $\widehat{g}(\bx)= \sum_{|\alpha| \leq 2c} g_{\alpha}\bx^{\alpha}$.} 
Let $\varepsilon_{\alpha} ={\rm sign}(g_\alpha)$.
Combining what we have done so far, we obtain that 
\begin{equation*}
    \langle q, \widehat{g}g_{J^\prime}q \rangle_{\iy} \;=\; 
    \sum_{|\alpha| \leq 2c} |g_{\alpha}|\langle q, \varepsilon_{\alpha}\bx^{\alpha} g_{J^\prime}q \rangle_{\iy} \;\leq\; R^{2c}\|\widehat{g}\|_1\langle q, g_{J^\prime}q \rangle_{\iy} \quad  \forall\; q \in \RR[\bx]_{r-\degg{g_{J^\prime}}-c}.
\end{equation*}
This completes the proof of the claim, which gives $\bM_{r-\degg{g_J}}(g_J\iy)=\bM_{r-\degg{g_J}}\big((R^{2c}\|\widehat{g}\|_1-\widehat{g})g_{J^\prime}\iy\big) \succeq 0$.

We can now prove the third condition. For notational convenience 
in the proof, we define
the auxiliary function $g_{m+1}(\bx)= K- \sum_{j=1}^mg_j(\bx)$. 
Recall the notation $g(\bx) = (g_1(\bx),\ldots,g_m(\bx)).$
For any $J \subset \{0,1,\dots,m+1 \}$ satisfying that $ \lceil p_J \rceil \leq t$, it is clear that $p_J(\bx,g(\bx)) \in \RR[\bx]_{4td}\subset \RR[\bx]_{2r}$, and there exists an $s(n+m,\lceil p_J \rceil) \times s(n,r)$ matrix $T_2$ satisfying that 
\begin{displaymath}
    \bv_{t-\lceil p_J \rceil}\big((\bx,g(\bx))\big) = T_2 \bv_r(\bx)  \;\Rightarrow\;  \bM_{t- \lceil p_J \rceil}(p_J \iy^{\varphi})= T_2 \bM_{r- \degg{g_J}}(g_J \iy)T_2^{\top}\succeq 0.
\end{displaymath}
This completes the proof.
 \end{proof}
 
We next present a lemma that shows the properties of a projection from $\RR^{s(n+m,k)}$ onto $\RR^{s(n,k)}$, which projects a pseudo-moment sequence of higher dimension to one of lower dimension.
\begin{lemma}\label{lem: projection into lows dimension}
    For a positive integer $k=2l$, we define the projection $\psi_k:\; \RR^{s(n+m,k)} \to \RR^{s(n,k)}$ as follows:
    \begin{displaymath}
        \iy \in \RR^{s(n+m,k)} \mapsto  \psi_k(y) \;\; 
        \mbox{such that}\;\;
        (\psi_k(y))_{\alpha} = y_{(\alpha,0_m)}\;
        \forall\; 
        \alpha \in \NN^n_k.
    \end{displaymath}
    Then for any $2r\geq k$, the followings hold true.
    \begin{enumerate}
        \item If $\iy \in \CM_k(\CT(\mathbb{B}_R \times \Delta^m_K)_{2r})$, then $\psi_k(\iy) \in \CM_k(\CT(\mathbb{B}_R)_{2r})$.
        \item If $\iy \in \CM_k(\CT(\varphi(\CX))_{2r})$, then $\psi_k(\iy) \in \CM_k(\CT(\CX)_{2r})$.
        \item If $\iy \in \CM_k(\mathbb{B}_R \times \Delta^m_K)$, then $\psi_k(\iy) \in \CM_k(\mathbb{B}_R)$.
        \item If If $\iy \in \CM_k(\varphi(\CX))$, then $\psi_k(\iy) \in \CM_k(\CX)$.
    \end{enumerate}
\end{lemma}
\begin{proof} We only prove the first property since the others can
be proved similarly. 
    We first consider the case $k=2r$. In this case,
    \begin{displaymath}
        \CM_k(\CT(\mathbb{B}_R \times \Delta^m_K)_{2r})= \CM(\CT(\mathbb{B}_R \times \Delta^m_K)_{2r}). \quad 
    \end{displaymath}
    For any $\iy \in \CM_k(\CT(\mathbb{B}_R \times \Delta^m_K)_{2r})$, based on the definition of $\psi_k$, it is clear that $\bM_r(\psi_k(\iy))$ and $\bM_{r-1}\bigl((R^2-\sum_{j=1}^nx_j^2)\psi_k(\iy)\bigr)$ are  principle submatrices of $\bM_r(\iy)$ and $\bM_{r-1}\bigl((R^2-\sum_{j=1}^nx_j^2)\iy\bigr)$, respectively. Hence, we obtain
    \begin{displaymath}
       \bM_r(\psi_k(\iy)) \succeq 0, \quad \mbox{and} \quad  
       \bM_{r-1}\Big(\big(R^2-\sum_{j=1}^nx_j^2\big)\psi_k(\iy)\Big) \succeq 0 \; \Rightarrow \; \psi_k(\iy) \in \CM(\CT(\mathbb{B}_R)_{2r}).
    \end{displaymath}
    For the case $k \leq 2r$, since $\CM_k(\CT(\mathbb{B}_R)_{2r})$ 
    is the first $s(n,k)$--coordinate 
    projection of $\CM(\CT(\mathbb{B}_R)_{2r})$, 
    we obtain that $\psi_k(\iy) \in \CM_k(\CT(\mathbb{B}_R)_{2r})$.
\end{proof}

We use both Lemmas~\ref{lemma: lift of moment sequence} and~\ref{lem: projection into lows dimension} to prove the following theorem on bounding the Hausdorff distance $\bd_k(\CR(\CX)_{2r})$.
 \begin{theorem}\label{thm: error on general}
     Let $\CX \in \RR^n$ be a compact basic semi-algebraic set defined as in \eqref{eq-CX}. For any positive integers $k=2l$ and $r$ such that $t = \lfloor r/(2d) \rfloor$, $l \geq 2d$, and $t \geq 4(\max\{n,m\}+1)k+2$ the Hausdorff distance $\bd_k(\CR(\CX)_{2r})$ admits the following bound:
     \begin{multline*}
          \bd_k(\CR(\CX)_{2r}) \leq \Gamma(\mathbb{B}_R \times \Delta_K^m,k)2\gamma(n+m,R+mK^2,k)\dfrac{(2d)^2}{(r-4d)^2}\\+\frac{cL(R,k)(2d)^{\Lo}}{(r-4d)^{\Lo}}
            [\Gamma(\mathbb{B}_R \times \Delta_K^m,k)2\gamma(n+m,R+mK^2,k)]^{\Lo/2}\left(m+\sum_{j=1}^m\|g_j\|_1+\sum_{i=1}^p\|h_i\|_1\right)^{\Lo}.
     \end{multline*}
     Here, the parameter $L(R,k)$ is the Lipschitz number of $\bv_k(\bx)$ on the simple set $\mathbb{B}_R$. In short, the error for the set of truncated pseudo-moment sequences $\CM_k(\CR(\CX)_{2r})$
     is $\mathrm{O}(1/r^{\Lo})$. 
 \end{theorem}
 \begin{proof}
     Let $\iy \in \CM_k(\CR(\CX)_{2r})$ be an arbitrary truncated pseudo-moment sequence. We proceed to bound the distance from $\iy$ to $\CM_k(\CX)$ by finding an appropriate point $\widetilde{\iy} \in \CM_k(\CX)$ and a point $\overline{\iy} \in \CM_k(\mathbb{B}_R)$ such that
     \begin{equation}\label{eq: inequality of distance with compactness}
         \bd(\iy,\CM_k(\CX)) \leq \|\iy- \widetilde{\iy}\| \leq \|\iy- \overline{\iy}\| + \|\overline{\iy}- \widetilde{\iy}\|.
     \end{equation}
     The idea of finding $\overline{\iy}$ and $\widetilde{\iy}$ is through lifting $\CM_k(\CX)$ to $\CM(\CT(\varphi(\CX))_{2r})$, where $ \varphi(\CX)$ is the intersection of a real variety and the product $\mathbb{B}_R \times \Delta_K^m$ of simple sets. 
     We first find the corresponding points 
     $\overline{y}'$ and $\widetilde{y}'$ in the lifted spaces and then 
     obtain the desired points $\overline{y}$ and $\widetilde{y}$ by the 
     projection $\psi_k$ defined as in Lemma~\ref{lem: projection into lows dimension}.
     We can conduct the proof as in Theorem~\ref{thm: error on product of simple sets}. The detail is elaborated below.
     \begin{enumerate}
         \item Since $\iy \in \CM_k(\CR(\CX)_{2r})$, there exists $\widehat{\iy} \in \CM(\CR(\CX)_{2r})$ such that $\iy = \pi_k(\widehat{\iy})$. Lemma~\ref{lemma: lift of moment sequence} implies that 
         $\widehat{\iy}^{\varphi} \in \CM(\CT(\varphi(\CX))_{2t})$. Let $\overline{\iy}^\prime \in \CM_{2t}(\mathbb{B}_R \times \Delta_K^m)$ be the projection of $\widehat{\iy}^{\varphi}$ onto $\CM_{2t}(\mathbb{B}_R \times \Delta_K^m)$. We set
         \begin{equation*}
             \overline{\iy}= \psi_k(\pi_k(\overline{\iy}^\prime)) \in \CM_k(\mathbb{B}_R) \quad \text{(by Lemma~\ref{lem: projection into lows dimension})}.
         \end{equation*}
         \item By Tchakaloff's theorem, there exist at most $N = s(n+m,2t)$ points $\{\overline{\bz}_1,\dots,\overline{\bz}_N\} \subset \mathbb{B}_R \times \Delta_K^m$ and positive weights $\{w_s : s\in [N]\}$ satisfying that $\sum_{s=1}^Nw_s =1$ and 
         \begin{displaymath}
             \overline{\iy}^\prime = \sum_{s=1}^N w_s
             \bv_{2t}(\overline{\bz}_s).
         \end{displaymath}
         For any $s \in [N]$, $\overline{\bz}_s = (\overline{\bx}_s, \overline{u}_{s1},\dots,\overline{u}_{sm})$, we set ${\bx}_s$ to be the projection of $\overline{\bx}_s$ onto $\CX$. 
         Then we define 
         \begin{eqnarray*}
            \widetilde{\iy}^\prime &=&\sum_{s=1}^Nw_s \bv_{2t}\big( ({\bx}_s,g_1({\bx}_s),\dots,g_m({\bx}_s))\big)\in \CM_{2t}(\varphi(\CX)), 
            \\
            \widetilde{\iy} &=& \psi_k(\pi_k(\widetilde{\iy}^\prime))= \sum_{s=1}^Nw_s \bv_k({\bx}_s) \in \CM_k(\CX).
         \end{eqnarray*}
     \end{enumerate}
     Figure~\ref{fig-lifting} illustrates our idea. 
         After defining $\overline{\iy}$ and $\widetilde{\iy}$, we evaluate each term in \eqref{eq: inequality of distance with compactness}.
         
    \begin{figure}[!ht]
    \centering
    \begin{subfigure}{0.41\textwidth}
    \begin{tikzpicture}

    \begin{scope}[blend mode=multiply]
    \filldraw[fill=gray!10, draw=black] (1,0) ellipse (3.3cm and 1.3cm);
    \filldraw[fill=white!10, draw=black] (0,0) circle (1.2cm);
    \filldraw[fill=gray!10, draw=black] (0,0) ellipse (1.7cm and 2.6cm);
    \filldraw[fill=white!10, draw=black] (1,0) circle (3.5cm);
    \end{scope}
    
    \node at (0,0) {$\CM_{2t}(\varphi(\CX))$};
    \node at (3.0,0) {$\CM_{2t}(\mathbb{B}_R \times \Delta_K^m)$};
    \node at (0,1.6) {$\CM(\CR(\varphi(\CX))_{2t})$};
    \node at (1.0,-2.9) {$\CM(\CT(\mathbb{B}_R \times \Delta_K^m)_{2t})$};
    \node at (0,-2.3) (S1) {$\widehat{\iy}^{\varphi}$};
    \node at (2.7,-0.8) (S2) {$\overline{\iy}^\prime$};
    \node at (0,-0.5) (S3) {$\widetilde{\iy}^\prime$};

    \draw[arrow] (S1.east) -- (S2.west);
    \draw[arrow] (S2.west) -- (S3.south);
    \draw[arrow] (S1.east) -- (S3.south);
\end{tikzpicture}
\end{subfigure}
\hfill
\qquad\tikz[overlay,remember picture] {\draw[thick, ->] (-1.5,3.5) --node[above]{$\psi_k \circ \pi_k$}node[below]{$\substack{
    \iy= \psi_k(\pi_k(\widehat{\iy}^{\varphi}))\\
    \overline{\iy}= \psi_k(\pi_k(\overline{\iy}^\prime))\\
    \widetilde{\iy}= \psi_k(\pi_k(\widetilde{\iy}^\prime))}$} (0.7,3.5);}\qquad
\begin{subfigure}{0.41\textwidth}
    \begin{tikzpicture}

    \begin{scope}[blend mode=multiply]
    \filldraw[fill=gray!10, draw=black] (1,0) ellipse (3.3cm and 1.5cm);
    \filldraw[fill=white!10, draw=black] (0,0) circle (1.2cm);
    \filldraw[fill=gray!10, draw=black] (0,0) ellipse (1.7cm and 2.7cm);
    \filldraw[fill=white!10, draw=black] (1,0) circle (3.5cm);
    \end{scope}
    
    \node at (0,0) {$\CM_{k}(\CX)$};
    \node at (3.0,0) {$\CM_k(\mathbb{B}_R)$};
    \node at (0,1.8) {$\CM_k(\CR(\CX)_{2r})$};
    \node at (1.2,-3.0) {$\CM(\CT(\mathbb{B}_R)_{2r})$};
    \node at (0,-2.3) (S1) {$\iy$};
    \node at (2.7,-0.8) (S2) {$\overline{\iy}$};
    \node at (0,-0.7) (S3) {$\widetilde{\iy}$};

    \draw[arrow] (S1.east) -- (S2.west);
    \draw[arrow] (S2.west) -- (S3.south);
    \draw[arrow] (S1.east) -- (S3.south);
\end{tikzpicture}
\end{subfigure}
\caption{Lifting of $\CM_k(\CX)$ to $\CM_{2t}(\varphi(\CX))$.}
\label{fig-lifting}
\end{figure}
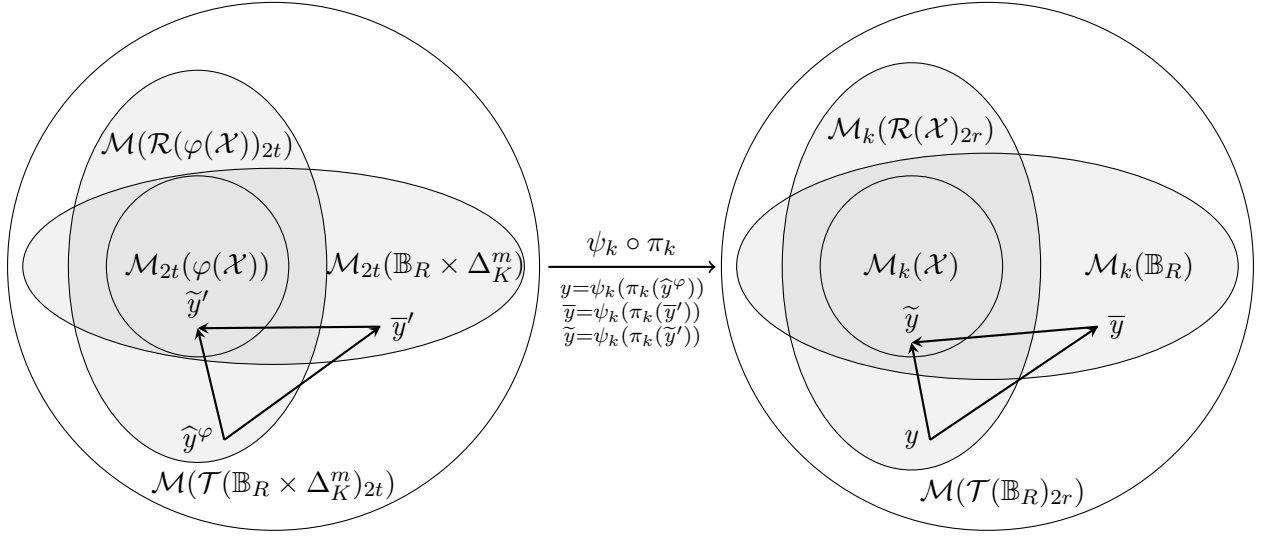

     \textbf{Evaluating $\|\iy -\overline{\iy}\|$:} 
     We apply Lemma~\ref{lem: bounded of M_k(B_R x Delta)} to $\bd_k(\CT(\mathbb{B}_R \times \Delta_K^m)_{2t})$ under the degree condition $t \geq 4(\max\{n,m\}+1)k+2$ for a product of two simple sets  to obtain that 
     \begin{multline}\label{eq: y to y bar}
       \|\iy - \overline{\iy}\| =
       \|\psi_k(\pi_k(\widehat{\iy}^\varphi)) - \psi_k(\pi_k(\overline{\iy}^\prime))\|
       \leq \|\pi_k(\widehat{\iy}^{\varphi})- \pi_k(\overline{\iy}^\prime)\| 
       \\ 
       \;\leq\; \bd_k(\CT(\mathbb{B}_R \times \Delta_K^m)_{2t}) \;\leq\; \dfrac{\Gamma(\mathbb{B}_R \times \Delta_K^m,k)2\gamma(n+m,R+mK^2,k)}{(t-2)^2} =: \varepsilon.
     \end{multline}
     
     \textbf{Evaluating $\|\overline{\iy}-\widetilde{\iy}\|$:} We use the same argument as in the proof of Theorem~\ref{thm: error on variety} to derive the following inequality:
     \begin{displaymath}
          \sum_{s=1}^N w_s \max\left\{ |h_i(\overline{\bx}_s)|:\ i \in [p]\right\} \leq \varepsilon^{1/2} \sum_{i=1}^p \|h_i\|_1.
     \end{displaymath}
    
     For $j \in [m]$, we consider $q_j(\bz) = (g_j(\bx)-u_j)^2$. Then we have
     \begin{align*}
         \ell_{\widehat{\iy}^{\varphi}}(q_j) = \ell_{\iy}((g_j(\bx)-g_j(\bx))^2) =0,\quad \mbox{and} \quad \|q_j\|_1\leq (1+\|g_j\|_1)^2.
     \end{align*}
     Note that ${\rm deg}(q_j) \leq 4d \leq k.$ Next, using the Cauchy–Schwarz inequality and the fact that 
     $\overline{\bu}_s=(\overline{u}_{s1},\ldots,\overline{u}_{sm}) \geq 0\; \forall\; s\in [N]$, we have 
     \begin{align*}
     &\sum_{s=1}^Nw_sq_j(\overline{\bz}_s) = 
     \sum_{s=1}^Nw_sq_j(\overline{\bz}_s)-\ell_{\widehat{\iy}^{\varphi}}(q_j) 
     = \left\langle {\bf q}_j,\overline{\iy}^\prime - \widehat{\iy}^{\varphi} \right\rangle
     \leq \varepsilon  \|q_j\|_1 \quad \mbox{(by \eqref{eq: y to y bar})}
     \\
     \Rightarrow & 
     \sum_{s=1}^N w_s(g_j(\overline{\bx}_s)-\overline{u}_{sj})^2 \leq \varepsilon \|q_j\|_1 
     \quad \Rightarrow \quad 
     \sum_{s=1}^Nw_s \max \{0, -g_j(\overline{\bx}_s)\}^2 \leq \varepsilon \|q_j\|_1\\
     \Rightarrow & \sum_{s=1}^N w_s\max \{0, -g_j(\overline{\bx}_s)\} \leq \left( \sum_{s=1}^Nw_s\right)^{1/2}\left( \sum_{s=1}^Nw_s \max \{0, -g_j(\overline{\bx}_s)\}^2\right)^{1/2} \leq \varepsilon^{1/2}(1+\|g_j\|_1)\\
     \Rightarrow & \sum_{s=1}^Nw_s \max \left\{-g_j(\overline{\bx}_s),|h_i(\overline{\bx}_s)|:\ i\in [p],\ j\in [m]\right\}\leq \varepsilon^{1/2}\left(m+\sum_{j=1}^m\|g_j\|_1+\sum_{i=1}^p\|h_i\|_1\right).
     \end{align*}
     
     Furthermore, the Łojasiewicz inequality~\eqref{Lojaciewicz inequality} implies that
    \begin{displaymath}
        \|\bx_s -\overline{\bx}_s\| = \bd(\overline{\bx}_s,\CX) \leq c\max \left\{-g_j(\overline{\bx}_s),|h_i(\overline{\bx}_s)|:\ i\in [p],\ j\in [m]\right\}^{\Lo}.
    \end{displaymath}
    
    Let  $L(R,k)>0$ be the Lipschitz number of $\bv_k(\bx)$ on the ball $\mathbb{B}_R$. Then it directly leads to the inequality:
    \begin{displaymath}
        \|\bv_k(\bx_s) - \bv_k(\overline{\bx}_s)\| \leq L(R,k) \|\bx_s -\overline{\bx}_s\|  \leq cL(R,k) \max \left\{-g_j(\overline{\bx}_s),|h_i(\overline{\bx}_s)|:\ i\in [p],\ j\in [m]\right\}^{\Lo}.
    \end{displaymath}
    Hence, we obtain that 
    \begin{eqnarray*}
        \|\overline{\iy}-\widetilde{\iy} \| &=&
        \left\|\sum_{s=1}^N w_s\bv_k(\overline{\bx}_s) -
        \sum_{s=1}^N w_s\bv_k(\bx_s)\right\|
        \\  
        &\leq& cL(R,k)\sum_{s=1}^N w_s\max \left\{-g_j(\overline{\bx}_s),|h_i(\overline{\bx}_s)|:\ i\in [p],\ j\in [m]\right\}^{\Lo}.
    \end{eqnarray*} 
    Recall that $\Lo \leq 1$, applying Jensen's inequality gives us the following inequality:
        \begin{align}\label{eq: y bar to y tidle}
            &\sum_{s=1}^N w_s\max \left\{-g_j(\overline{\bx}_s),|h_i(\overline{\bx}_s)|:\ i\in [p],\ j\in [m]\right\}^{\Lo}\nonumber\\
            \leq &\ \left(\sum_{s=1}^N w_s\max \left\{-g_j(\overline{\bx}_s),|h_i(\overline{\bx}_s)|:\ i\in [p],\ j\in [m]\right\}\right)^{\Lo}\nonumber\\
            \Rightarrow &\ \sum_{s=1}^N w_s\max \left\{-g_j(\overline{\bx}_s),|h_i(\overline{\bx}_s)|:\ i\in [p],\ j\in [m]\right\}^{\Lo} \leq\varepsilon^{\Lo/2}\left(m+\sum_{j=1}^m\|g_j\|_1+\sum_{i=1}^p\|h_i\|_1\right)^{\Lo}\nonumber\\
            \Rightarrow&\ \|\overline{\iy}-\widetilde{\iy} \| 
            \leq cL(R,k)\varepsilon^{\Lo/2}\left(m+\sum_{j=1}^m\|g_j\|_1+\sum_{i=1}^p\|h_i\|_1\right)^{\Lo}.
        \end{align}
        We substitute \eqref{eq: y to y bar} and \eqref{eq: y bar to y tidle} back to \eqref{eq: inequality of distance with compactness} to obtain that
        \begin{align*}
            &\bd(\iy,\CM_k(\CX)) \leq \varepsilon+cL(R,k)
            \varepsilon^{\Lo/2}\left(m+\sum_{j=1}^m\|g_j\|_1+\sum_{i=1}^p\|h_i\|_1\right)^{\Lo} \quad 
            \forall \iy \in \CM_k(\CR(\CX)_{2r})
            \\
            \Rightarrow& \ \bd_k(\CR(\CX)_{2r}) \leq \varepsilon+cL(R,k)
            \varepsilon^{\Lo/2}\left(m+\sum_{j=1}^m\|g_j\|_1+\sum_{i=1}^p\|h_i\|_1\right)^{\Lo}.
        \end{align*}
        By the definition $t$, we have $t \geq \frac{r-2d}{2d}$. Hence, we obtain that
        \begin{multline*}
            \bd_k(\CR(\CX)_{2r}) \leq \Gamma(\mathbb{B}_R \times \Delta_K^m,k)2\gamma(n+m,R+mK^2,k)\dfrac{(2d)^2}{(r-4d)^2}
            \\
            +\frac{cL(R,k)(2d)^{\Lo}}{(r-4d)^{\Lo}}
            [\Gamma(\mathbb{B}_R \times \Delta_K^m,k)2\gamma(n+m,R+mK^2,k)]^{\Lo/2}
        \left(m+\sum_{j=1}^m\|g_j\|_1+\sum_{i=1}^p\|h_i\|_1\right)^{\Lo}.
        \end{multline*}
        This completes the proof.
 \end{proof}
 
 \begin{corollary}\label{cor: convergence rate for compact semi-algebraic set}
     Let $\CX$ be the set and $k,r$ be the numbers defined as in Theorem~\ref{thm: error on general}. Then the following inequality holds:
     \begin{multline*}
          \fmin-\lb(f,\CR(\CX))_r \leq \Biggl[\Gamma(\mathbb{B}_R \times \Delta_K^m,k)2\gamma(n+m,R+mK^2,k)\dfrac{(2d)^2}{(r-4d)^2}\\+\frac{cL(R,k)(2d)^{\Lo}}{(r-4d)^{\Lo}}
            [\Gamma(\mathbb{B}_R \times \Delta_K^m,k)2\gamma(n+m,R+mK^2,k)]^{\Lo/2}\left(m+\sum_{j=1}^m\|g_j\|_1+\sum_{i=1}^p\|h_i\|_1\right)^{\Lo}\Biggr]\cdot \|f\|_1.
     \end{multline*}
     In conclusion, it is shown that the error for Schmüdgen-type truncated pseudo-moment sequences on a compact basic semi-algebraic $\CX$ is $\mathrm{O}(1/r^{\Lo})$, where $\Lo$ is the Łojasiewicz exponent depending on the polynomial inequalities defining $\CX$.
 \end{corollary}
 \begin{proof}
     The result follows straightforwardly from Theorem~\ref{thm: error on general} and Lemma~\ref{lemma: distance to convergence rate}.
 \end{proof}
 
 \subsection{Error estimation under the Polyak-Łojasiewicz condition}
 The last subsection has emphasized the connection between the error of pseudo-moment sequences and the Łojasiewicz exponent of the domain. In this subsection, we review the Polyak-Łojasiewicz inequality, which plays a significant role in the study of analytic gradient flows, to sharpen the Łojasiewicz exponent. We first set up our problem with additional assumptions. For simplicity, let $\CX$ be a compact basic semi-algebraic set
 defined as follows: 
 \begin{displaymath}
     \CX = \left\{\bx \in \RR^n:\ g_i(\bx) \leq 0 \ \forall i \in [m] \right\}.
 \end{displaymath}
 
 It is clear that the Hausdorff distance $
 \bd_k(\CR(\CX)_{2r})$ is $\mathrm{O}(1/r^{\Lo})$. Thus sharpening the exponent $\Lo$ would lead to a better bound on the error. To do so, we define the violating function $g$, which indicates how much the inequalities $g_j(\bx) \leq 0$ are violated at the point $\bx\in\RR^n$, i.e, 
 \begin{displaymath}
     g(\bx) = \max\{0,\ g_i(\bx):\ i \in [m]\} \quad \Rightarrow \quad\CX = \{\bx \in \RR^n:\ g(\bx) = 0\}.
 \end{displaymath}
We note that the function $g$ is non-smooth in general. To state the 
Polyak-Łojasiewicz (PŁ) condition, we first review the limiting subdifferential
of a nonsmooth function (see e.g., \cite[Chapter 8, 10]{variational_analysis}). In particular, the Fréchet subdifferential of $g$ at $\bx$, denoted by $\partial^Fg(\bx)$ is the set of vectors $\bw$ satisfying the following condition:
\begin{displaymath}
    \bw \in \partial^Fg(\bx) \ \Leftrightarrow \ \liminf_{\by \to \bx}\dfrac{g(\by) - g(\bx)-\langle \bw,\by-\bx \rangle}{\|\by -\bx\|} \geq 0.
\end{displaymath}
The limiting subdifferential of $g$ at $\bx$, denoted by $\partial^Lg(\bx)$, consists of vectors $\bw \in \partial^Fg(\bx)$ such that
\begin{displaymath}
    \exists \bx_n \to \bx,\ \exists \bw_n \to \bw, \quad \mbox{satisfying} \quad g(\bx_n) \to g(\bx)
     \;\;\mbox{and}\;\; 
    \bw_n \in \partial^Fg(\bx_n).
\end{displaymath}
We said that the function $g$ is globally $\mu-$PŁ for a positive number $\mu$ if 
\begin{equation}\label{def LP}
    \forall \bx \in \RR^n, \quad g(\bx) -\inf g \leq \dfrac{1}{2\mu}\|\bw\|^2\quad \forall\; \bw \in \partial^Lg(\bx).
\end{equation}
When the global PŁ condition is met, the Łojasiewicz exponent and the Łojasiewicz constant are explicitly defined in the following lemma (see e.g., \cite[Corollary 12]{Polyak-Łojasiewicz-for-nonsmooth}).

\begin{lemma}\label{lemma: PL }
    Let $q: \RR^n \to \RR$ be a continuous function whose set of global minimizers $\mathrm{argmin}\, q$ is nonempty. Assume that $q$ is globally $\mu-$PŁ with constant $\mu >0$. Then the Łojasiewicz exponent $\Lo = \frac{1}{2}$ and the Łojasiewicz constant is $\sqrt{\frac{2}{\mu}}$ for the distance to the set $\mathrm{argmin}\, q$, i.e,
    \begin{displaymath}
        \bd(\bx, \mathrm{argmin}\ q) \leq \sqrt{\dfrac{2}{\mu}} \big(q(\bx) -\inf_{\RR^n} q\big)^{1/2}, \quad \forall \bx \in \RR^n.
    \end{displaymath}
\end{lemma}
Therefore, if the violating function $g$ satisfies the P\L\ condition, i.e., the inequality \eqref{def LP} holds for some positive constant $\mu$, then the error for the set of truncated pseudo-moment sequences is $\mathrm{O}(1/\sqrt{r})$. However, the inequality 
\eqref{def LP} is challenging to check in practice. Thus, we relax the P\L\ condition by the strong convexity of the defining polynomials $\{g_i\mid i\in[m]\}$.

\begin{assume}\label{assumption: convexisty}
    We assume that $\CX$ is defined by the polynomial inequalities $g_i(\bx) \leq 0$ for all $i\in[m]$, where $g_i$'s are locally strongly convex function with the constant $\mu_i >0$, i.e., there exists a compact convex set $\Omega$ such that $\CX \subset \Omega$ and for any $\bx, \by \in \Omega$ and $i \in [m]$, the following inequality holds
    \begin{displaymath}
        g_i(\by) \geq g_i(\bx) + \left \langle \nabla g_i(\bx),\by -\bx \right \rangle + \dfrac{\mu_i}{2}\|\by-\bx\|^2.
    \end{displaymath}
    In addition, we set $\mu = \min\{ \mu_i,\ i \in [m]\} >0$. Then the inequality 
    \begin{displaymath}
        g_i(\by) \geq g_i(\bx) + \left \langle \nabla g_i(\bx),\by -\bx \right \rangle + \dfrac{\mu}{2}\|\by-\bx\|^2
    \end{displaymath}
    holds for all $\bx, \by \in \Omega$ and $i \in [m]$. In this case, we call the set $\CX$ a $\mu$-strongly convex semi-algebraic set.
\end{assume}

The work \cite{strong-convex-to-PL} implies that global strong convexity induces the global P\L\ condition. However, the Łojasiewicz inequality that has been used throughout this paper is local, i.e., the inequality holds true on a compact domain. Therefore, to sharpen the Łojasiewicz exponent, we analyse the connection between local strong convexity (LSC) and the local Polyak-Łojasiewicz condition (LP\L), which leads to the local Łojasiewicz inequality with explicit exponent and constant (LLI), i.e., a local version of Lemma~\ref{lemma: PL }. In summary, we aim to prove that 
\begin{displaymath}
    \textrm{(LSC)} \quad \xrightarrow{\text{Lemma }\ref{lemma: (LSC) to (LPL)}} \quad \textrm{(LP\L)} \quad\Longrightarrow \quad \textrm{(LLI)}.
\end{displaymath}

We first show the first connection where the (LSC) condition of $g_i$'s implies the (LP\L) property of the violating function $g$.
\begin{lemma}\label{lemma: (LSC) to (LPL)}
    Let $\CX$ be a non-empty, $\mu-$strongly convex semi-algebraic set contained in a compact convex set $\Omega$, i.e., $\CX \subset \Omega$,
    and for any $\bx, \by \in \Omega$ and $i \in [m]$, the following inequality holds for some positive parameter $\mu$:
    \begin{displaymath}
        g_i(\by) \geq g_i(\bx) + \left \langle \nabla g_i(\bx),\by -\bx \right \rangle + \dfrac{\mu}{2}\|\by-\bx\|^2.
    \end{displaymath}
    Then the function $g$ is $\mu-$P$\Lo$ on $\Omega$, i.e., for any $\bx \in \Omega$ and $\bw \in \partial^Lg(\bx)$, we have that 
    \begin{displaymath}
        g(\bx) -\inf_{\Omega} g \leq \dfrac{1}{2\mu}\|\bw\|^2.
    \end{displaymath}
\end{lemma}
\begin{proof}
    Since the functions $0, g_1, \dots, g_m$ are convex, $g$ is convex. Thus, the Fréchet subdifferential and the limiting subdifferential of $g$ are both equal to the classical subdifferential for a convex function~\cite[Proposition 8.12]{variational_analysis}, i.e.,
    \begin{displaymath}
        \partial^Fg(\bx) = \partial^Lg(\bx) = \partial g(\bx),
    \end{displaymath}
    where $\partial g(\bx)$ denotes the subdifferential of $g$. We next consider our desired inequality. For any $\bx \in \mbox{argmin}(g)$, i.e., $g(\bx)= \inf_{\Omega}g$, it is obvious that the following inequality
    \begin{equation*}
        g(\bx) - \inf_{\Omega}g = 0 \leq\frac{1}{2\mu}\|\bw\|^2, \quad \forall \ \bw \in \partial^Lg(\bx)
    \end{equation*} 
    holds true. For $\bx \in \Omega \backslash \mbox{argmin}(g)$, $g(\bx) > 0$. Thus, the set $A(\bx)$ of active indices at $\bx$, defined as
    \begin{displaymath}
         A(\bx) = \{i \in [m]:\ g(\bx) = g_i(\bx)\} \neq \emptyset.
    \end{displaymath}
     Additionally, we can combine it with \cite[Example 5.4.5]{Dimitri-convex} to obtain the subdifferential of the maximum of differentiable functions as
     \begin{displaymath}
         \partial^Fg(\bx) = \partial^Lg(\bx) = \partial g(\bx) = \conv\{\nabla g_i(\bx)\; :\; i \in A(\bx)\}.
     \end{displaymath}
     For any $\bw \in \partial g(\bx)$, we set $\bw = \sum_{i \in  A(\bx)} w_i\nabla g_i(\bx)$ where $0 \leq w_i \leq 1$ for all $i \in A(\bx)$ and $\sum_{i \in A(\bx)}w_i =1$. Then we have that 
    \begin{align*}
        &g_i(\by) \geq g_i(\bx) + \left \langle \nabla g_i(\bx),\by -\bx \right \rangle + \dfrac{\mu}{2}\|\by-\bx\|^2\\
        \Rightarrow & \  \sum_{i \in A(\bx)}w_ig_i(\by) \geq \sum_{i \in A(\bx)}w_ig_i(\bx) + \left \langle  \sum_{i \in A(\bx)}w_i\nabla g_i(\bx),\by -\bx \right \rangle + \sum_{i \in A(\bx)}w_i\dfrac{\mu}{2}\|\by-\bx\|^2\\
        \Rightarrow & \  g(\by) \geq g(\bx) + \left\langle \bw,\by-\bx \right\rangle + \dfrac{\mu}{2}\|\by-\bx\|^2, \ \forall \bx, \by \in \Omega,\ \bw \in \partial g(\bx).
    \end{align*}
    Since $\CX$ is a non-empty set contained in $\Omega$, there exists a point $\bx^* \in \mbox{argmin}(g) \subset \Omega$ . Applying the last inequality and the Cauchy–Schwarz inequality, we  obtain that 
    \begin{displaymath}
        \inf_{\Omega}\ g = g(\bx^*) \geq g(\bx) + \left\langle \bw,\bx^*-\bx \right\rangle + \dfrac{\mu}{2}\|\bx^*-\bx\|^2 \geq g(\bx) - \dfrac{1}{2\mu}\|\bw\|^2.
    \end{displaymath}
    Hence, we obtain that $g(\bx) - \inf_\Omega g \leq \frac{1}{2\mu}\|\bw\|^2 \ \forall \bw \in \partial g(\bx)= \partial^L g(\bx),\ \forall \bx \in \Omega$, which is the local version of the P\L\ condition \eqref{def LP} on $\Omega$.
 \end{proof}
 
 We next prove the local Łojasiewicz inequality of $g$ on $\Omega$ from the local 
 P\L\ condition by adopting the proof in \cite[Theorem 11]{Polyak-Łojasiewicz-for-nonsmooth} with some modifications for our case.

 \begin{lemma}\label{lemma: (LPL) to (LLI)}
     Let $\CX$ satisfies Assumption \ref{assumption: convexisty}. 
     Then the following Łojasiewicz inequality holds:
     \begin{displaymath}
         \bd(\bx,\CX) \leq \sqrt{\dfrac{2}{\mu}}
         g(\bx)^{1/2}\quad \forall\; \bx \in \Omega.
     \end{displaymath}
 \end{lemma}
 \begin{proof}
     According to Lemma~\ref{lemma: (LSC) to (LPL)}, the P\L\ condition holds on $\Omega$, i.e.,
     \begin{displaymath}
         g(\bx) - \inf_\Omega g \leq \frac{1}{2\mu}\|\bw\|^2 \quad \forall \bw \in \partial g(\bx),\quad \forall \bx \in \Omega.
     \end{displaymath}
     Now consider the function
     \begin{displaymath}
          f(\bx) = 
            \begin{cases} 
                g(\bx) & \text{if } \bx \in \Omega, \\
                +\infty & \text{if } \bx \notin \Omega,
            \end{cases}
     \end{displaymath}
     which is a convex lower semi-continuous function on $\RR^n$ with $\textrm{dom}\ f = \Omega$ and satisfies the P\L\ condition on $\Omega$. Observe that $\inf_{\RR^n} f = \inf_{\Omega} g =0$, and $\textrm{argmin}\ f = \CX$.
     For any $\widehat{\bx} \in \Omega$, we construct a sequence $(\bx_k)_{k \in \NN}$ as follows: $\bx_0 =\widehat{\bx}$, and 
     \begin{displaymath}
         \bx_{k+1} \in \underset{\bx \in \RR^n}{\textrm{argmin }} f(\bx) + \dfrac{1}{2}\|\bx-\bx_k\|^2, \quad S_{k+1} = \{\bx \in \RR^n:\ f(\bx) \leq f(\bx_{k+1})\}.
        \end{displaymath}
        Then the level sets $S_k \subset \Omega \ \forall k \in \NN$, and $f(\bx_{k+1}) + \dfrac{1}{2}\|\bx_{k+1} - \bx_k\|^2 \leq f(\bx) + \dfrac{1}{2}\|\bx - \bx_k\|^2\ \forall \bx \in \RR^n$. Hence, we obtain that 
        \begin{displaymath}
            \|\bx_{k+1}-\bx_{k}\|^2 \leq \|\bx-\bx_{k}\|^2 \;\; \forall\; \bx \in S_{k+1}\quad \Rightarrow \quad \|\bx_{k+1}-\bx_{k}\| = \bd(\bx_k,S_{k+1}).
        \end{displaymath}
        
        We next consider the function $\varphi(t) = \sqrt{\dfrac{2}{\mu}}t^{1/2}$. Then $\varphi'(t) =\dfrac{1}{\sqrt{2\mu}}t^{-1/2}$ and $\varphi^{-1}(t) = \dfrac{\mu}{2}t^2$. Since $\inf_{\RR^n} f =0$, the P\L\ condition can be written equivalently as
        \begin{displaymath}
            f(\bx) -\inf_{\RR^n} f \leq \dfrac{1}{2\mu}\|\bw\|^2 \quad \Leftrightarrow \quad 1 \leq \varphi'(f(\bx)) \|\bw\| \quad \forall \bx \in \Omega,\ \forall\; \bw \in \partial f(\bx).
        \end{displaymath}
        We use the chain rule to obtain that 
        \begin{displaymath}
            \partial(\varphi \circ f)(\bx) \subset \varphi'(f(\bx))\partial f(\bx) \quad \Rightarrow \quad 1 \leq \|\bv\| \;\; \forall \bv \in \partial(\varphi \circ f)(\bx),\quad \forall\; \bx\in \Omega.
        \end{displaymath}
        The above  result can be combined with the condition of \cite[Proposition 4.6]{curves-of-descent} to show that the limiting slope of $f$ at $\bx$ (see e.g., \cite{curves-of-descent,Characterizations-of-LP}), denoted by $\left|\overline{\nabla }f(\bx)\right|$, satisfies that
        \begin{displaymath}
            \left|\overline{\nabla }f(\bx)\right| = \bd(0, \partial f(\bx)) = \inf \left\{\|\bv\|,\ \bv \in \partial f(\bx) \right\} \geq 1 \ \forall\; \bx \in \Omega.
        \end{displaymath}
        
         We note that a convex function is subdifferential regular at every point of its effective domain (see e.g., \cite[Proposition 8.21]{variational_analysis}). Then the slope and the limiting slope coincides. Thus, the assumption on the K\L-inequality and 
         sub-level set mapping in \cite[Corollary 4]{Characterizations-of-LP} is satisfied, which when combines with the definition of $S_k = \{ \bx \in \RR^n:\ 0 \leq f(\bx) \leq f(\bx_k)\} \subset \Omega$ leads to the following inequality:
        \begin{displaymath}
            \|\bx_{k+1}-\bx_{k}\|= \bd(\bx_k,S_{k+1}) \leq \varphi(f(\bx_k)) -\varphi(f(\bx_{k+1})). 
        \end{displaymath}
        
        We note that $\bx_0 =\widehat{\bx}$. The above inequality gives
        \begin{equation} \label{eq-PL}
            \|\bx_{k+1}-\widehat{\bx}\| \leq \sum_{i=0}^k\|\bx_{i+1}-\bx_i\| \leq \sum_{i=0}^k 
            \varphi(f(\bx_{i}))-\varphi(f(\bx_{i+1})) \leq \varphi(f(\widehat{\bx})).
        \end{equation}
        Furthermore, since $(\bx_k)_{k \in \NN} \subset \Omega$,  the compactness of $\Omega$ ensures the existence of a limit point $\overline{\bx}$ of $\{\bx_k: k \in \NN\}$. Without loss of generality, we can assume that $\lim_{k \to \infty} \bx_k = \overline{\bx}$. The first-order optimality condition and the local P\L\ condition for $f$ imply that 
        \begin{align*}
            &\bx_k - \bx_{k+1} \in \partial f(\bx_{k+1}) \ \Rightarrow \ f(\bx_{k+1}) \leq \dfrac{1}{2\mu}\|\bx_k-\bx_{k+1}\|^2 \\ \Rightarrow &\ f(\overline{\bx}) = \lim_{k \to \infty}f(\bx_k) \leq \frac{1}{2\mu}\lim_{k \to \infty}\|\bx_k-\bx_{k+1}\|^2 =0.
        \end{align*}
        Hence, $\overline{\bx} \in \textrm{argmin\;}f$. Moreover, \eqref{eq-PL} implies that 
        $$
         \|\widehat{\bx}-\overline{\bx}\| \leq \varphi(f(\widehat{\bx})). 
        $$
        Now the increasing property of $\varphi^{-1}$ implies the following inequality:
        \begin{displaymath}
            \varphi^{-1}(\bd(\widehat{\bx}, \textrm{argmin }f)) \leq \varphi^{-1}(\|\widehat{\bx}-\overline{\bx}\|) \leq f(\widehat{\bx})\quad \Rightarrow \quad \bd(\widehat{\bx},\CX) \leq \sqrt{\dfrac{2}{\mu}}g(\widehat{\bx})^{1/2}.
        \end{displaymath}
        Since $\widehat{\bx}\in \Omega$ is arbitrary, the last inequality
        holds for all $\widehat{\bx}\in \Omega.$
        This completes the proof.
 \end{proof}
 \begin{theorem}\label{thm: error and rate under the strong convexity}
            Let $\CX \subset \interior{\Omega}$ be either a strongly convex semi-algebraic set as in Assumption \ref{assumption: convexisty} or the violating function $g$ satisfies the Polyak-Łojasiewicz condition, where $\interior{\Omega}$ is 
            the interior of the convex set $\Omega$. Then the error for the set of truncated pseudo-moment sequences and the convergence rate of the Schm\"udgen-type moment-SOS hierarchy
            $\CM_k(\CR(\CX)_{2r})$ is $\mathrm{O}(1/\sqrt{r})$. 
\end{theorem}
\begin{proof}
    We recall that under either the strong convexity condition in Assumption \ref{assumption: convexisty} or the Polyak-Łojasiewicz condition, the Łojasiewicz exponent of $\CX$ on $\Omega$ is $1/2$. To apply Theorem~\ref{thm: error on general} with $\Lo = 1/2$, it is necessary to show that the Łojasiewicz exponent of $\CX$ on $\mathbb{B}_R$ is also $1/2$ for some radius $R$ such that $\CX \subset \mathbb{B}_R$. This relies on the containment of $\CX$ within the interior of $\Omega$, which allows the Łojasiewicz exponent of $1/2$ to be extended to $\mathbb{B}_R$. This concept is illustrated in Figure~\ref{fig: exponent}.
    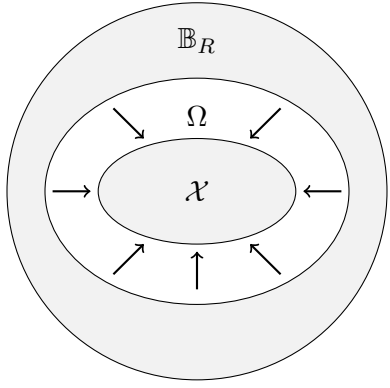
\begin{figure*}
    \centering
    \hfill\\
    \begin{tikzpicture}
    \filldraw[fill=gray!10, draw=black] (0,0) circle (2.5cm);
    \filldraw[fill=white!10, draw=black] (0,0) ellipse (2cm and 1.5cm);
    \filldraw[fill=gray!10, draw=black] (0,0) ellipse (1.3cm and 0.7cm);
    
    \draw[->, thick] (-1.9,0) -- (-1.4,0); 
    \draw[->, thick] (1.9,0) -- (1.4,0); 
    \draw[->, thick] (-1.1,-1.1) -- (-0.7,-0.7); 
    \draw[->, thick] (-1.1,1.1) -- (-0.7,0.7); 
    \draw[->, thick] (1.1,1.1) -- (0.7,0.7); 
    \draw[->, thick] (1.1,-1.1) -- (0.7,-0.7); 
    \draw[->, thick] (0,-1.3) -- (0,-0.8); 
    \node at (0,0) {$\mathcal{X}$};
    \node at (0,1) {$\Omega$};
    \node at (0,2) {$\mathbb{B}_R$};
    \end{tikzpicture}
    \caption{Proof of Theorem~\ref{thm: error and rate under the strong convexity}.}
    \label{fig: exponent}
    \end{figure*}
    Note that from Lemmas \ref{lemma: (LSC) to (LPL)} and \ref{lemma: (LPL) to (LLI)}, there exist the the Łojasiewicz exponent $1/2$ and constant $c$ such that 
    \begin{equation*}
        \bd(\bx,\CX) \leq c g(\bx)^{1/2} \quad \forall\; \bx \in \Omega.
    \end{equation*}
    Since $\CX \subset \interior{\Omega}$, $g(\bx)$ is continuous on $\mathbb{B}_R$, and $\CX = g^{-1}(0)$,  we obtain that
    \begin{displaymath}
    \min_{\bx \in \mathbb{B}_R \backslash \interior{\Omega}}\frac{cg(\bx)^{1/2}}{\bd(\bx,\CX)} > 0, \quad \eta := \min \left\{ 1,\min_{\bx \in \mathbb{B}_R \backslash \interior{\Omega}}\frac{cg(\bx)^{1/2}}{\bd(\bx,\CX)}\right\} > 0.
    \end{displaymath}
    Then we can construct a version of the Łojasiewicz inequality of $\CX$ over $\mathbb{B}_R$ as follows: 
    \begin{align*}
        &\bd(\bx,\CX) \leq c g(\bx)^{1/2} \leq \frac{c}{\eta} g(\bx)^{1/2} \quad \forall\; \bx \in \Omega, \quad \bd(\bx,\CX) \leq \frac{c}{\eta} g(\bx)^{1/2} \quad \forall\; \bx \in \mathbb{B}_R\backslash \interior{\Omega}\\
        \Rightarrow \;& \bd(\bx,\CX) \leq \frac{c}{\eta} g(\bx)^{1/2} \quad \forall \; \bx \in \mathbb{B}_R.
    \end{align*} 
    This implies that the Łojasiewicz exponent of $\CX$ on $\mathbb{B}_R$ is $1/2$, which is what we desire. 
\end{proof}

\subsection{Convergence rates for some special cases}

First, we consider the case when $\CX$ is a polytope. The following lemma identify the Łojasiewicz exponent of $\CX$.

\begin{lemma}[\cite{bergthaller1992distance},Theorem 0.1]\label{lemma: bound on distance to a polyhedron}
     Let $A$ be an $m \times n$ matrix, and  $\beta$ be the least number such that for each non-singular sub-matrix $B$ of $A$, all entries of $B^{-1}$ are at most $\beta$ in absolute value. Consider the polyhedron $P$ defined by a system of linear inequalities $A\bx^\prime \leq \bb^\prime$. Then for each $\bb^0 \in \RR^m$ and $\bx^0 \in \RR^n$ such that $A\bx^0 \leq \bb^0$, there exists $\bx^\prime \in \RR^n$ satisfying 
     \begin{align*}
         &A\bx^\prime \leq \bb^\prime,\quad \mbox{and}\quad \|\bx^0-\bx^\prime\|_{\infty} \leq n\beta \|\bb^0-\bb^\prime\|_{\infty}.
     \end{align*}
     Consequently, the Łojasiewicz exponent of $P$ is $1$. Here, $\|\cdot\|_{\infty}$ denotes the infinity norm of a vector.
 \end{lemma}
 \begin{theorem}\label{cor: error and rate for polytope}
      Let $\CX$ be a polytope. Then the Łojasiewicz exponent of $\CX$ is $1$. Consequently, we obtain the followings.
      \begin{enumerate}
            \item For any fixed positive integer $k$, the error $\bd_k(\CT(\CX)_{2r})$ of pseudo-moment sequences on $\CX$ is of
            the order $\mathrm{O}(1/r)$. 
            \item The convergence rate of the reduced moment-SOS hierarchy \eqref{hierarchy: reduced moment},~\eqref{hierarchy: reduced SOS}, and the Schmüdgen-type moment-SOS hierarchy \eqref{hierarchy: moment Schmudgen}, \eqref{hierarchy: SOS Schmudegen} are of the order $\mathrm{O}(1/r)$.
      \end{enumerate}
 \end{theorem}
 \begin{proof}
     The theorem is a corollary of Theorem~\ref{thm: error on general} and Corollary~\ref{cor: convergence rate for compact semi-algebraic set} when the Łojasiewicz exponent is $1$.
 \end{proof}

In what follows, we extend the results of \cite{fang2021sum} on the convergence rate of
$O(1/r^2)$ of the moment-SOS hierarchy over the sphere $ S_R^{n-1} = \{\bx \in \RR^n : \|\bx\| = R\} $. The established convergence rate of $\mathrm{O}(1/r^2)$ 
in \cite{fang2021sum} applies to polynomial optimization problems  
on  $S_R^{n-1}$  with homogeneous objective functions. Our goal here is to generalize this result to POPs with possibly non-homogeneous  polynomial objective functions. To achieve this, we refine the bounds on the pseudo-moment sequences for  $S_R^{n-1}$. Observe that for any positive integer $r$,
\begin{equation*}
    \CM(\CT(S_R^{n-1})_{2r})=\CM(\CQ(S_R^{n-1})_{2r})= \left\{\iy \in \RR^{s(n,2r)}:\; \bM_r(\iy) \succeq 0,\; 
    \bM_{r-1}((R^2-\|\bx\|^2)\iy)=0 \right\}.
\end{equation*}
The error for the pseudo-moment sequences on $S_R^{n-1}$ is sharpened in the following theorem.
\begin{theorem}\label{thm: error on sphere}
    Let $k,l$ and $r$ be positive integers such that $k=2l$ and $2r\geq k$. Then the Hausdorff distance $\bd_k(\CQ(S_R^{n-1})_{2r}) =\bd_k(\CT(S_R^{n-1})_{2r})$ admits the following upper bound:
    \begin{equation}\label{error on sphere}
        \bd_k(\CT(S_R^{n-1})_{2r}) \leq   
        \left(1+ \frac{\sqrt{n}L(R,k)}{R} \right)\frac{2\gamma(R,n,k)\Gamma(\mathbb{B}_R,k)}{r^2}.
    \end{equation}
\end{theorem}
\begin{proof}
    We first proceed as in the proof of Theorem~\ref{thm: error on variety}. Let $\iy \in \CM_k(\CT(S_R^{n-1})_{2r})$. Then $\iy \in \CM_k(\CT(\mathbb{B}_R)_{2r})$, and there exists the projection $\overline{\iy}$ of $\iy$ onto $\CM_k(\mathbb{B}_R)$ such that 
    \begin{displaymath}
        \|\iy -\overline{\iy}\| \;\leq\; \frac{2\gamma(R,n,k) \Gamma(\mathbb{B}_R,k)}{r^2}=: \varepsilon.
    \end{displaymath}
    Moreover, there exists $N \leq  s(n,k) $ points $\{\overline{\bx}_s:\; s\in[N]\}$ in $\mathbb{B}_R$ and positive weight $\{w_s:\; s\in [N]\}$ with $\sum_{s=1}^Nw_s =1$ satisfying that 
    \begin{displaymath}
        \overline{\iy}= \sum_{s=1}^Nw_s\bv_k(\overline{\bx}_s).
    \end{displaymath}
    Applying the Cauchy–Schwarz inequality to $\left|\ell_{\iy}(R^2-\|\bx\|^2)-\ell_{\overline{\iy}}(R^2-\|\bx\|^2)\right|$, we obtain the following:
    \begin{align*}
        &\left|\ell_{\iy}(R^2-\|\bx\|^2)-\ell_{\overline{\iy}}(R^2-\|\bx\|^2)\right|  
        \leq  \sum_{i=1}^n |y_{2e_i} - \overline{y}_{2e_i}|
        \leq \sqrt{n}\|y-\overline{y}\|      
        \\
        \Rightarrow\;& \left|\ell_{\overline{\iy}}(R^2-\|\bx\|^2) \right| 
        \leq \sqrt{n}\varepsilon \quad \text{(since $\ell_{\iy}(R^2 - \|\bx\|^2) = 0$)}\\
        \Rightarrow\;& \sum_{s=1}^Nw_s(R^2-\|\overline{\bx}_s\|^2) \leq \sqrt{n}\varepsilon.
    \end{align*}
    Define $\widetilde{\bx}_s = R\,{\overline{\bx}_s}/{\|\overline{\bx}_s\|} \in S_R^{n-1}$ to be the projection of $\overline{\bx}_s$ onto $S_R^{n-1}$.
    Then we can bound the distance between $\overline{\bx}_s$ and 
    $\widetilde{\bx}_s$ as follows:
    \begin{displaymath}
        \|\widetilde{\bx}_s-\overline{\bx}_s\| = \left\|\overline{\bx}_s\left(\frac{R}{\|\overline{\bx}_s\|}-1 \right)\right\|= R-\|\overline{\bx}_s\| 
        \leq \frac{(R+\|\overline{\bx}_s\|)(R-\|\overline{\bx}_s\|)}{R}= \frac{R^2-\|\overline{\bx}_s\|^2}{R}.
    \end{displaymath}
    We set $\widetilde{\iy}= \sum_{s=1}^Nw_s\bv_k(\widetilde{\bx}_s)$. Since $\widetilde{\bx}_s \in S_R^{n-1}\; \forall s \in [N]$, $\widetilde{\iy} \in \CM_k(S_R^{n-1})$. Moreover, let $L(R,k)$ be the Lipschitz number of $\bv_k(\bx)$ on the ball $\mathbb{B}_R$. Then we can perform the following evaluation on the distance between $\overline{\iy}$ and $\widetilde{\iy}$ as follows:
    \begin{align*}
        &\|\overline{\iy}-\widetilde{\iy}\| \leq \sum_{s=1}^Nw_s\|\bv_k(\overline{\bx}_s)-\bv_k(\widetilde{\bx}_s)\| \leq \sum_{s=1}^Nw_sL(R,k)\|\overline{\bx}_s-\widetilde{\bx}_s\|\\
        \leq\; & \sum_{s=1}^Nw_sL(R,k)\frac{R^2-\|\overline{\bx}_s\|^2}{R} 
        = \frac{L(R,k)}{R}\sum_{s=1}^Nw_s(R^2-\|\overline{\bx}_s\|^2) \leq \frac{\sqrt{n}L(R,k)}{R}\varepsilon.
    \end{align*}
    Hence, we can bound the distance from $\iy$ to $\CM_k(S_R^{n-1})$ as 
    \begin{displaymath}
        \bd_k(\CT(S_R^{n-1})_{2r}) \leq \|\iy -\widetilde{\iy}\| \leq \|\iy - \overline{\iy}\| + \| \overline{\iy}-\widetilde{\iy}\| \leq \varepsilon\left(1+ \frac{\sqrt{n}L(R,k)}{R} \right). 
    \end{displaymath}
    In conclusion, we obtain the desired bound.
    This completes the proof.
\end{proof}
\begin{corollary}\label{cor: rate on sphere}
    Consider the following POP on $S_R^{n-1}$:
    \begin{displaymath}
        \min_{\bx \in S_R^{n-1}}f(\bx),
    \end{displaymath}
    where $f(x)$ is a polynomial of degree $k$. Then the convergence rate of the moment-SOS hierarchy $\mlb(f,\CQ(\CX))_r$ and $\lb(f,\CQ(\CX))_r$ is of the order $\mathrm{O}(1/r^2)$.
\end{corollary}
\begin{proof}
    The proof is similar to the proof of Lemma~\ref{lemma: distance to convergence rate}. We can obtain that 
    \begin{eqnarray*}
        | \fmin - \mlb(f,\CQ(\CX))_r| &\leq& \|f\|_1\,
        \bd_k(\CQ(S_R^{n-1})_{2r})\\
        &\leq& \|f\|_1\left(1+ \frac{\sqrt{n}L(R,k)}{R} \right)\frac{2\gamma(R,n,k)\Gamma(\mathbb{B}_R,k)}{r^2}.
    \end{eqnarray*}
    The convergence rates of $\mlb(f,\CQ(\CX))_r$ and $\lb(f,\CQ(\CX))_r$ are since strong duality holds under the Archimedean condition, which is satisfied in this case. The proof is completed.
\end{proof}

In the paper~\cite{Baldi_2025}, Baldi shows that the Łojasiewicz exponent 
of $\CX$ is $1$ under the Constraint Qualification Condition (CQC), which is stated next. Hence, we can sharpen our results under that case.
\begin{proposition}\text{\cite[Theorem 2.10 and Theorem 2.14]{Baldi_2025}}\label{thm: CQC}
    Consider the domain 
    \begin{displaymath}
        \CX = \{\bx \in \RR^n:\; g_i(\bx) \geq 0,\; i\in [m] \}. 
    \end{displaymath}
    For $\bx \in \CX$, let $I(\bx)$ be the set of active indices (the index $i$ is said to be active if $g_i(x) = 0$). We say that the Constraint Qualification Condition (CQC) holds at $\bx$ if $\{ \nabla g_i(x):\; i \in I(\bx)\} $ are linearly independent. 

    We say that $\CX$ satisfy the CQC if the CQC holds at every point of $\CX$. In this case, the Łojasiewicz exponent is equal to $1$.
\end{proposition}
\begin{corollary}\label{cor: error and rate under CQC}
    If $\CX$ satisfies the CQC, then the error $\bd_k(\CT(\CX)_{2r})$ of pseudo-moment sequences and the convergence rate of the reduced moment-SOS hierarchy~\eqref{hierarchy: reduced moment},~\eqref{hierarchy: reduced SOS}, and the Schm\"udgen-type moment-SOS hierarchy \eqref{hierarchy: moment Schmudgen}, \eqref{hierarchy: SOS Schmudegen} are of the order $\mathrm{O}(1/r)$.
\end{corollary}

\section*{Conclusion} Our work has demonstrated a strong connection between the convergence rate of the moment-SOS hierarchy for a compact semi-algebraic set and the \L{}ojasiewicz exponent of the domain. This insight provides a novel framework for analyzing the behavior of the moment-SOS hierarchy. However, the methodology appears to be limited to the Schm\"udgen-type hierarchy. An intriguing direction for future research would be to extend this approach to analyze the convergence rate of the Putinar-type hierarchy. Additionally, we have shown in this paper that the Schm\"udgen-type hierarchy can be simplified without affecting its theoretical convergence rate. This raises an important question for future investigation: What kind of reductions can be applied to the moment-SOS hierarchy to further lower the computational cost of the SDP relaxation?

\bibliographystyle{alpha}
\bibliography{references}

\appendix

\section{Convergence rate under a linear transformation}\label{appendix: linear transformation}
    
    In this section, we prove that the error of the pseudo-moment sequences and the convergence rate of the moment-SOS hierarchy are asymptotically invariant under a linear transformation. Equivalently, if either one of the Hausdorff distances $\bd_k(\CT(\CX)_{2r})$, $\bd_k(\CQ(\CX)_{2r})$, $\bd_k(\CR(\CX)_{2r})$ or one of the convergence rates of the hierarchies \eqref{hierarchy: moment Schmudgen}, \eqref{hierarchy: moment Putinar},    
    \eqref{hierarchy: reduced moment},
    is of the order $\mathrm{O}(1/r^c)$ for some constant $c$ over a domain $\CX$, then it is also valid for any image of $\CX$ under an invertible linear transformation. The next theorem captures this fact.
    \begin{theorem}\label{thm: error on image}
        Let $\CX$ be a product of simple sets defined by 
        \begin{displaymath}
        \CX := \{ \bx \in \RR^n:\ g_j(\bx) \geq 0\ \forall j \in [m]\}.
    \end{displaymath}
     Let $k= 2l$ be a positive even integer, $A$ be an invertible matrix in $\RR^{n\times n}$. Then the image of $\CX$ via $A$ remains a basic semi-algebraic set defined by
    \begin{displaymath}
        A(\CX) = \{\bx \in \RR^n:\  g_j(A^{-1}\bx) \geq 0\ \forall j \in [m]\}.
    \end{displaymath}
    We have the following upper bound on the error of truncated pseudo-moment sequences:
    \begin{displaymath}
        \bd_k(\CT(A(\CX)))_r \leq \Gamma(A(\CX),k)\dfrac{\gamma(R,n,k)}{r^2},
    \end{displaymath}
    where $\Gamma(A(\CX)),k)$ is a polynomial in the dimension $n$, the norm of  $A$ and $k$. The result also holds true for $\CQ(A(\CX))_{2r}$ and $\CR(A(\CX))_{2r}$.
    \end{theorem}

    \begin{remark}
        Theorem~\ref{thm: error on image}  can be proved similarly as in Theorem~\ref{thm: error on product of simple sets} when $\CX$ is a product of simple sets by constructing a push-forward measure on 
        $A(\CX)$. In particular, let $\mu $ be the measure on $\CX$ used in the proof of Theorem~\ref{thm: error on product of simple sets}. We consider the push-forward measure $\mu \circ A^{-1}$, and note that $A$ is invertible. Then the CD kernel on $A(\CX)$ can be reconstructed, for which the same arguments as in the proof of Theorem~\ref{thm: error on product of simple sets} remain valid. However, we propose here another proof for a general set $\CX$ that also explains how a pseudo-moment sequence is transformed under the action of an invertible linear transformation on the domain.
    \end{remark}
\begin{proof}
        Denote the rows of $A$ by $a_i^\top$ for $i \in [n]$, and we define an isomorphism $\overline{A}$ on the $\RR[\bx]-$algebra induced by 
        $A$ as follows:
    \begin{align*}
        \overline{A}: \RR[\bx] \to \RR[\bx],  
        \quad 
        \overline{A}\bx^\alpha = (A\bx)^\alpha = \prod_{i=1}^n \langle a_i,\bx \rangle^{\alpha_i}.
    \end{align*}
    Since $A$ is invertible,  for all $\alpha \in \NN^n$, the degree of the polynomial $(A\bx)^{\alpha}=\prod_{i=1}^n\langle a_i, \bx \rangle^{\alpha_i}$ remains to be $|\alpha|$. We denote the restriction of $\overline{A}$ on $\RR[\bx]_{2r}$ by $\overline{A}_{2r}$. Fix $\bv_{2r}(\bx) = (\bx^{\alpha})_{\alpha \in \NN^n_{2r}}$ to be a basis of the linear space $\RR[\bx]_{2r}$. Then for any $p(\bx) \in \RR[\bx]_{2r}$, there exists a unique vector $\mathbf{p} \in \RR^{s(n,2r)}$ such that 
    \begin{displaymath}
        p(\bx)= \sum_{\alpha \in \NN^n_{2r}}p_{\alpha}\bx^{\alpha} = \langle \mathbf{p},\bv_{2r}(\bx)\rangle, \quad \mathbf{p}:= (p_{\alpha})_{\alpha \in \NN^n_{2r}}.
    \end{displaymath}
    The mapping $\varphi_{2r}: p(\bx) \mapsto \mathbf{p}$ is an isomorphism between $\RR[\bx]_{2r}$ and $\RR^{s(n,2r)}$. We next define the linear transformation $\mathcal{A}_{2r}:\RR^{s(n,2r)}\to \RR^{s(n,2r)}$ as follows:
    \begin{equation}\label{define transformation}
        \mathcal{A}_{2r}= \varphi_{2r}\circ \overline{A}_{2r}\circ \varphi_{2r}^{-1}.
    \end{equation}
     Let 
     $\{e_{\alpha}\}_{\alpha \in \NN^n_{2r}}$ be the standard basis of $\RR^{s(n,2r)}$.
    In particular, $\mathcal{A}_{2r}(e_\alpha) = \varphi_{2r}((A\bx)^\alpha)$ for all $\alpha\in \NN^n_{2r}$.
    
     We claim that $\mathcal{A}_{2r}(\CM(\CT(\CX)_{2r}))= \CM(\CT(A(\CX))_{2r})\ \forall r \in \NN$. Indeed, let $\iy \in \CM(\CT(\CX)_{2r})$, i.e., $\iy$ satisfies the following conditions
     \begin{displaymath}
         \iy_0 = 1, \quad \bM_r(\iy) \succeq 0,\quad \bM_{r-\degg{g_J}}(g_J\iy) \succeq 0 \quad \forall J \subset [m],\; 
         \degg{g_J} \leq r.
     \end{displaymath}
    Set $\overline{\iy}= \mathcal{A}_{2r}(\iy)$. To prove that $\overline{\iy} \in \CM(\CT(A(\CX))_{2r})$, we check the following constraints on $\overline{\iy}$ 
    \begin{equation}\label{condition 1}
       \overline{\iy}_0 = 1, \quad  \bM_r(\overline{\iy}) \succeq 0,\quad \bM_{r-\degg{g_J}}((g_J\circ A^{-1})\overline{\iy}) \succeq 0 \quad \forall J \subset [m],\; \degg{g_J} \leq r.
     \end{equation}
     Let $[\overline{A}_{2r}]$ be the representation matrix of $\overline{A}_{2r}$ on $\RR[\bx]_{2r}$ with respect to the basis $\{\bx^{\alpha}\}_{\alpha \in \NN^n_{2r}}$.
     Then for any $\alpha \in \NN^n_{2r}$, we have that 
     \begin{align*}
         \mathcal{A}_{2r} (e_{\alpha}) &\;= \varphi_{2r} \circ \overline{A}_{2r} \circ \varphi_{2r}^{-1}(e_{\alpha})
         \;=\; \varphi_{2r} \circ \overline{A}_{2r}(\bx^{\alpha})
        \;=\; \varphi_{2r} ( (A\bx)^\alpha)
         \;=\; [\overline{A}_{2r}] e_{\alpha}
     \end{align*}
     and hence $\overline{y}= \mathcal{A}_{2r}(y) = [\overline{A}_{2r}] y.$
     So $[\overline{A}_{2r}]$ is also the representation matrix of $\mathcal{A}_{2r}$ on $\RR^{s(n,2r)}$ with respect to the standard basis.
     Moreover, for any $\alpha,\, \beta \in \NN^n_r$, the $(\alpha,\beta)$-entry of $[\overline{A}_{2r}]$ is $e_{\alpha}^{\top}[\overline{A}_{2r}]e_{\beta} = e_\beta^\top [\overline{A}_{2r}]^\top e_\alpha $, which is the coefficients of $\bx^{\beta}$ in the polynomial $(A\bx)^{\alpha}$.
     Therefore, we have 
     \begin{displaymath}
         \overline{\iy}_0 =\langle e_0, \overline{\iy} \rangle= \langle  e_0, [\overline{A}_{2r}] y \rangle
         = \sum_{\alpha} y_\alpha 
         \langle e_0, \varphi_{2r}((A\bx)^\alpha) \rangle
         = y_0 =1.
     \end{displaymath}
     For the rest of the conditions in  \eqref{condition 1}, we adopt the convention that $g_{\emptyset} =1$ and $d_{\emptyset} =0$. 
     Then $\bM_r(\overline{\iy}) = \bM_{r-d_{\emptyset}}( (g_{\emptyset}\circ A^{-1})\overline{\iy})$. 
     Hence, it suffices to prove that for any $J \subset [m]$ and 
     $\degg{g_J} \leq r$, $\bM_{r-\degg{g_J}}((g_J\circ A^{-1})\overline{\iy}) \succeq 0$.
     Let $t=r-\degg{g_J}$.  We will prove that $\bM_{t}((g_J\circ A^{-1})\overline{\iy})= [\overline{A}_{t}]^{\top}\bM_{t}(g_J\iy)[\overline{A}_{t}]$, and the semidefiniteness of $\bM_{t}((g_J\circ A^{-1})\overline{\iy})$ will then follow from that of $\bM_{t}(g_J\iy)$. Here $\overline{A}_{t}$ is the restriction of $\overline{A}$ on $\RR[\bx]_{t}$, whose representation matrix with respect to the basis $\{\bx^{\alpha}\}_{\alpha \in \NN^n_{t}}$ is $[\overline{A}_{t}]$. 
     Let $d_J = \degg{g_J}.$ For any $\alpha,\beta \in \NN^n_t$, 
     we observe that 
     \begin{align*}
         &\bM_{t}((g_J\circ A^{-1})\overline{\iy})(\alpha,\beta)= \sum_{\gamma \in \NN^n_{2d_J}}(g_J \circ A^{-1})_{\gamma}\overline{\iy}_{\alpha + \beta+\gamma}
         \\
         =\;&\sum_{\gamma \in \NN^n_{2d_J}}(g_J \circ A^{-1})_{\gamma}\left\langle e_{\alpha+\beta+\gamma},[\overline{A}_{2r}]\iy \right\rangle
         \\
         =\;& \sum_{\tau \in \NN^n_{2r}}\sum_{\gamma \in \NN^n_{2d_J}}(g_J \circ A^{-1})_{\gamma} [\overline{A}_{2r}]_{\alpha+\beta+\gamma, \tau}\cdot\iy_{\tau}.
     \end{align*}
     We consider the coefficient of $\iy_{\tau}$ in the last sum
     \begin{align*}
         &\sum_{\gamma \in \NN^n_{2d_J}}(g_J \circ A^{-1})_{\gamma} [\overline{A}_{2r}]_{\alpha+\beta+\gamma, \tau}
         = \left \langle e_{\tau},\, \sum_{\gamma \in \NN^n_{2d_J}}(g_J \circ A^{-1})_{\gamma} [\overline{A}_{2r}]^\top e_{\alpha+\beta+\gamma}\right \rangle
         \\
         =&\; \left \langle e_{\tau},\, \varphi_{2r}\Big(\sum_{\gamma \in \NN^n_{2d_J}}(g_J \circ A^{-1})_{\gamma} (A\bx)^{\alpha+\beta+\gamma}
         \Big)\right \rangle
         =\left \langle e_{\tau},\, \varphi_{2r}\left((A\bx)^{\alpha+\beta}\cdot 
         (g_J \circ A^{-1}) (A\bx) \right)\right \rangle\\
         =&\; \left \langle  e_{\tau},
         \varphi_{2r} ((A\bx)^{\alpha+\beta}\cdot g_J(\bx)) \right \rangle,
     \end{align*}
     which is the coefficient of $\bx^{\tau}$ in the polynomial $(A\bx)^{\alpha+\beta}\cdot g_J(\bx)$.
     
     We next consider the $(\alpha,\beta)-$entry of $[\overline{A}_{t}]^{\top}\bM_{t}(g_J\iy)[\overline{A}_{t}]$ given as follows:
     \begin{align*}
         &[\overline{A}_{t}]^{\top}\bM_{t}(g_J\iy)[\overline{A}_{t}](\alpha,\beta)= ([\overline{A}_{t}]e_{\alpha})^{\top}\bM_{t}(g_J\iy)([\overline{A}_{t}]e_{\beta})
         \;=\;\sum_{\sigma, \theta \in \NN^n_t}[\overline{A}_t]_{\sigma,\alpha}[\overline{A}_t]_{\theta,\beta}\sum_{\gamma \in \NN^n_{2d_J}}(g_J)_{\gamma}\iy_{\sigma+\theta+\gamma}.
     \end{align*}
     Then for any $\tau \in \NN^n_{2r}$, the coefficient of $\iy_{\tau}$ in the last sum is 
     \begin{align*}
         \sum_{\substack{\sigma, \theta \in \NN^n_r, \gamma \in \NN^n_{2d_J}\\ \sigma+\theta +\gamma= \tau}}[\overline{A}_t]_{\sigma,\alpha}[\overline{A}_t]_{\theta,\beta}(g_J)_{\gamma},
     \end{align*}
     which is the coefficient of $\bx^{\tau}$ in the product $g_J(\bx)(A\bx)^{\alpha}(A\bx)^{\beta}=(A\bx)^{\alpha+\beta}\cdot g_J(\bx) $. Hence, the $(\alpha,\beta)-$entries of $\bM_{t}((g_J\circ A^{-1})\overline{\iy})$ and $ [\overline{A}_{t}]^{\top}\bM_{r-\degg{g_J}}(g_J\iy)[\overline{A}_{t}]$ coincide. Consequently, we have 
     \begin{displaymath}
         \mathcal{A}_{2r}(\CM(\CT(\CX)_{2r}))\subset \CM(\CT(A(\CX))_{2r}).
     \end{displaymath}
     
     Since $A$ is invertible, we use the same argument to prove that 
     \begin{displaymath}
         \CM(\CT(A(\CX))_{2r})=\mathcal{A}^{-1}_{2r}(\CM(\CT(A(\CX))_{2r}))\subset \CM(\CT(\CX)_{2r}),
     \end{displaymath}
     where $\mathcal{A}^{-1}_{2r} =\varphi_{2r}\circ 
     (\overline{A}_{2r})^{-1}\circ \varphi_{2r}^{-1}$, and $(\overline{A}_{2r})^{-1}$ is the inverse of $\overline{A}_{2r}$.
     Applying $\mathcal{A}_{2r}$ to the above inclusion, we get 
     \begin{displaymath}
         \CM(\CT(A(\CX))_{2r})\subset \mathcal{A}_{2r}
         (\CM(\CT(\CX)_{2r})).
     \end{displaymath}
     Thus we have proved that $ \mathcal{A}_{2r}
         (\CM(\CT(\CX)_{2r}))=\CM(\CT(A(\CX))_{2r}).$

     Now consider any $2r-$truncated moment sequence $y\in \CM_{2r}(\CX)$.
    By Tchakaloff's theorem, 
     there exist at most $N = s(n,2r)$ points $\{\bx_1,\ldots,\bx_N\}$ in $\CX$ and positive scalars 
     $\{w_1,\ldots,w_N\}$ satisfying $\sum_{i=1}^Nw_i=1$ such that 
     \begin{displaymath}
         \iy = \sum_{i=1}^Nw_i\bv_{2r}(\bx_i).
     \end{displaymath}
     Then  $\mathcal{A}_{2r}(\iy) =  \sum_{i=1}^Nw_i \varphi_{2r}(\bv_{2r}(A\bx_i))$, and $A\bx_i \in A(\CX) \; \forall i \in [N]$.
     Hence, we also have $\mathcal{A}_{2r}(\CM_{2r}(\CX)) = \CM_{2r}(A(\CX))$. Combining this with the result in the last paragraph, we obtain that 
     \begin{displaymath}
         \bd_{k}(\CT(A(\CX))_{2r}) \leq \|\pi_k \circ \mathcal{A}_{2r}\|_2 \,\bd_k(\CT(\CX)_{2r})\quad \forall k \leq 2r.
     \end{displaymath}
     In conclusion, the error estimation for the $k-$truncated 
pseudo-moment problem on $A(\CX)$ is conveyed from that on $\CX$.
The same holds true for the convergence rate of the moment-SOS hierarchy
on $A(\CX)$.
    \end{proof}
    \begin{remark}\label{rem: scalling}
        Notice that $\mathbb{B}_R$ is the image of the unit ball $B_n$ via the scalling: 
    \begin{align*}
       & A: \RR^n \to \RR^n: \quad \bx \mapsto A\bx = R\bx\\
       \Rightarrow\; & \mathcal{A}_{2r}(e_{\alpha}) = R^{|\alpha|}e_{\alpha} \quad \Rightarrow \;\|\pi_k \circ \mathcal{A}_{2r}\|_2 = R^{k} \quad \forall R \geq 1.
    \end{align*}
    This implies that 
    $\Gamma(\mathbb{B}_R,k) = R^{k}\cdot\Gamma(B_n,k).$
    A similar result also holds true for the simplex $\Delta^m_K$ and the product $\mathbb{B}_R \times \Delta^m_K$, i.e., we have 
    \begin{eqnarray*}
        & \Gamma(\Delta_K^m,k) = K^{k}\cdot\Gamma(\Delta_m,k) \quad \forall K \geq 1,\\
        & \Gamma(\mathbb{B}_R \times \Delta^m_K) = \max\{R,K\}^{k}\cdot \Gamma(B_n \times \Delta_m,k) \quad \forall R, K \geq 1.
    \end{eqnarray*}
    \end{remark}

    \section{The harmonic constant}\label{appendix: harmonic constant}
    
    In this appendix, we provide a quantitative analysis 
on the harmonic constant 
defined in \eqref{eq-Lambda}
in the following proposition.
    \begin{proposition}
        Let $\Lambda(\CX,k)$ be the harmonic constant defined as in \eqref{eq-Lambda} on the product $\CX$ of  simple sets $\CX_i\subset\RR^{n_i}$ for $i \in [m]$. Then $\Lambda(\CX,k)$ depends polynomially on $k$ (for fixed $\CX$) and polynomially on $(n_1,\dots,n_m)$ (for fixed $k$).
    \end{proposition}
    \begin{proof}
       Let $p \in \RR[\bx]_k$ be a polynomial of degree $k$ and we assume that $\|p\|_{\CX} =1 $. Recall that $\mu = \otimes_{i=1}^m\mu_i$, where $\mu_i$ is the measure corresponding to the simple set $\CX_i$ as in Table~\ref{tab:measure}. We know from the proof of Lemma~\ref{lem: equivalent kernel} that for $j_1,\dots,j_m$ such that $j_1 + \dots + j_m \leq k$,
       \begin{equation*}
           p_{j_1,\dots,j_m}(\bx) = \int_{\CX}\left[\prod_{i=1}^m C^{(j_i)}[\CX_i,\mu_i](\bx^{(i)}, \overline{\bx}^{(i)})\right]p(\overline{\bx})d\mu(\overline{\bx}) \quad \forall\; \bx \in \CX.
       \end{equation*}
       We use the Cauchy–Schwarz inequality and the fact that $\|p\|_{\CX} = 1$ to obtain that 
       \begin{align*}
           |p_{j_1,\dots,j_m}(\bx)|^2 &= \left|\int_{\CX}\left[\prod_{i=1}^m C^{(j_i)}[\CX_i,\mu_i](\bx^{(i)}, \overline{\bx}^{(i)})\right]p(\overline{\bx})d\mu(\overline{\bx}) \right|^2\\
           &\leq \int_{\CX}\left[\prod_{i=1}^m C^{(j_i)}[\CX_i,\mu_i](\bx^{(i)}, \overline{\bx}^{(i)})\right]^2d\mu(\overline{\bx}) \cdot \int_{\CX}p(\overline{\bx})^2d\mu(\overline{\bx})\\
           &\leq \prod_{i=1}^m \int_{\CX_i}C^{(j_i)}[\CX_i,\mu_i](\bx^{(i)}, \overline{\bx}^{(i)})^2d\mu_i(\overline{\bx}^{(i)}).
       \end{align*}
       Using the property of the CD kernel, we have:
       \begin{equation*}
           \int_{\CX_i}C^{(j_i)}[\CX_i,\mu_i](\bx^{(i)}, \overline{\bx}^{(i)})^2d\mu_i(\overline{\bx}^{(i)}) = C^{(j_i)}[\CX_i,\mu_i](\bx^{(i)},\bx^{(i)}).
       \end{equation*}
       Then it follows that 
       \begin{align*}
           \Lambda(\CX,k)^2 &= \max_{p \in \RR[\bx]_k}\max_{j_1 + \dots + j_m \leq k} \frac{\|p_{j_1,\dots,j_m}\|_{\CX}^2}{\|p\|_{\CX}^2} \leq \max_{j_1 + \dots + j_m \leq k} \max_{\bx \in \CX}\ \prod_{i=1}^mC^{(j_i)}[\CX_i,\mu_i](\bx^{(i)},\bx^{(i)})\\
           &\leq \max_{j_1 + \dots + j_m \leq k} \prod_{i=1}^m \max_{\bx^{(i)} \in \CX_i} C^{(j_i)}[\CX_i,\mu_i](\bx^{(i)},\bx^{(i)})\\
           &\leq \prod_{i=1}^m \tau(\CX_i,k), \quad \mbox{where} \quad \tau(\CX_i,k):=\max_{0 \leq j_i \leq k}\max_{\bx^{(i)} \in \CX_i} C^{(j_i)}[\CX_i,\mu_i](\bx^{(i)},\bx^{(i)}).
        \end{align*}
        We recall from the works \cite{slot2111sum} and \cite{laurent2023effective} that $\tau(\CX_i,k)$ depends polynomially on $n_i$ for fixed $k$ and polynomially on $k$ for fixed $n_i$. This leads to our required result. We refer to \cite{slot2111sum} and \cite{laurent2023effective} for the explicit bound on each $\tau(\CX_i,k)$ that can be used  to derive an explicit bound on $\Lambda(\CX,k)$ in terms of $k$ and $n_i$'s.
    \end{proof}

\end{document}